\documentclass[12pt,a4paper]{article}
\usepackage[latin1]{inputenc}
\usepackage[T1]{fontenc}
\usepackage{amsmath}
\usepackage{amssymb}
\usepackage{textcomp}
\usepackage[english]{babel}
\usepackage[all,dvips]{xy}

\newtheorem{thm}{Theorem}[section]
\newtheorem{lem}[thm]{Lemma}
\newtheorem{prop}[thm]{Proposition}
\newtheorem{cor}[thm]{Corollary}

\newcommand\dem{\textbf{Proof: }}
\newcommand\Def{\textbf{Definition: }}
\newcommand\Defs{\textbf{Definitions: }}
\newcommand\rem{\textbf{Remark: }}
\newtheorem{soulem}[thm]{Sublemma}
\newtheorem{blank}{}

\title{Strong almost reducibility for analytic and Gevrey quasi-periodic cocycles}
\author{Claire Chavaudret\\
Institut de Mathématiques de Jussieu \\
175 rue du Chevaleret, 75013 Paris, France }
\date{}

\begin{document}



\maketitle

Abstract: This paper is about almost reducibility of quasi-periodic cocycles with a diophantine frequency which are sufficiently close to a constant. Generalizing previous works by L.H.Eliasson, we show a strong version of almost reducibility for analytic and Gevrey cocycles, that is to say, almost reducibility where the change of variables is in an analytic or Gevrey class which is independent of how close to a constant the initial cocycle is conjugated. This implies a result of density, or quasi-density, of reducible cocycles near a constant. Some algebraic structure can also be preserved, by doubling the period if needed. 

%

\section{Introduction}


We are concerned with quasi-periodic cocycles, that is, solutions of equations of the form

\begin{equation}\label{defducc}\forall (\theta,t)\in 2\mathbb{T}^d\times \mathbb{R},\ \frac{d}{dt}X^t(\theta)=A(\theta+t\omega)X^t(\theta); \ X^0(\theta)=Id
\end{equation}

\noindent 
where $A\in C^0(2\mathbb{T}^d, \mathcal{G})$ and $\mathcal{G}$ is a linear Lie algebra. Here $\mathbb{T}^d=\mathbb{R}^d/\mathbb{Z}^d$ stands for the $d$-torus, $d\geq 1$,  and  $2\mathbb{T}^d=\mathbb{R}^d/(2\mathbb{Z}^d)$ stands for 
the double torus. We will assume in this article that $\omega\in\mathbb{R}^d$ satisfies some diophantine conditions. 
The solution of \eqref{defducc} is called the quasi-periodic cocycle associated to $A$ and is defined on $2\mathbb{T}^d\times \mathbb{R}$ with values in the connected component of the identity of a Lie group $G$ whose associated Lie algebra is $\mathcal{G}$. Terminology is explained by the fact that $A$ is the envelope of a quasi-periodic function, since $t\mapsto A(\theta+t\omega ) $ is a quasi-periodic function for all $\theta\in 2\mathbb{T}^d$. We say $X$ is a constant cocycle if $A$ is constant. A constant cocycle is always of the form $t\mapsto e^{tA}$. 

\bigskip
\noindent A cocycle is said to be reducible if it is conjugated to a constant cocycle, in a sense that will be defined later on. The problem of reducibility of cocycles has been thoroughly studied and is of interest because the dynamics of reducible cocycles is well understood and because this problem has links with the spectral theory of Schrödinger cocycles and with the problem of lower dimensional invariant tori in hamiltonian systems. In the periodic case ($d=1$), Floquet theory tells that every cocycle is reducible modulo a loss of periodicity. However, the problem is far more difficult if $d$ is greater than 1 and it is not true that every cocycle is then reducible. The question becomes whether every cocycle is close, up to a conjugacy, to a reducible one; from this question comes the notion of almost-reducibility. A cocycle is said to be almost-reducible, roughly speaking, if it can be conjugated to a cocycle which is arbitrarily close to a reducible one. Reductibility implies almost reducibility, however the reverse is not true: there are non reducible cocycles even close to a constant cocycle (see \cite{El92}). Almost reducibility is an interesting notion since the dynamics of an almost reducible cocycle are quite well known on a very long time.

\bigskip
\noindent We first focus on cocycles generated by functions which are analytic on a neighbourhood of the torus, i.e real analytic functions which are periodic in the direction of the real axis (recall that they are matrix-valued). For such a function $F$, we will let 

$$\mid F\mid_r=\sup_{\mid \operatorname{Im}\theta\mid \leq r} \mid \mid F(\theta)\mid \mid$$
 
where $\mid \mid .\mid \mid$ stands for the operator norm. 

\bigskip
\noindent The aim of this paper is to show that for

$$G=GL(n,\mathbb{C}),GL(n,\mathbb{R}),SL(2,\mathbb{C}),SL(n,\mathbb{R}),Sp(n,\mathbb{R}),O(n),U(n) $$

\noindent in the neighbourhood of a constant cocycle, every cocycle which is analytic on an $r$-neighbourhood of the torus and $G$-valued
 is almost reducible in $C^\omega_{r'}(2\mathbb{T}^d,G)$ for all ${0<r'<r\leq \frac{1}{2}}$. The width of the neighbourhood only depends on the dimensions $n,d$, 
on the diophantine class of $\omega$, on the constant cocycle and on the loss of analyticity $r-r'$. 

\bigskip
\noindent More precisely, we shall prove the following theorem, for $G$ among the groups cited above and $\mathcal{G}$ the Lie algebra associated to $G$:

\begin{thm}\label{th1}Let $0<r'<r\leq \frac{1}{2}$, $A\in \mathcal{G}$, $F\in C^\omega_r(\mathbb{T}^d, \mathcal{G})$. 
There is $\epsilon_0<1$ 
depending only on $n,d,\omega,A,r-r'$ such that if

$$|F|_r\leq \epsilon_0$$ 

\noindent 
then for all $\epsilon>0$, there exists $\bar{A}_\epsilon,\bar{F}_\epsilon
\in C^\omega_{r'}
(2\mathbb{T}^d,\mathcal{G})$, 
$\Psi_\epsilon,Z_\epsilon\in C^\omega_{r'}(2\mathbb{T}^d, G)$ and $A_\epsilon\in \mathcal{G}$ 
such that for all $\theta\in 2\mathbb{T}^d$,

$$\partial_\omega Z_\epsilon(\theta)=(A+F(\theta))Z_\epsilon(\theta)-Z_\epsilon(\theta)(\bar{A}_\epsilon(\theta)
+\bar{F}_\epsilon(\theta))
$$

\noindent with
\begin{enumerate}
\item \label{redtoby}
$\partial_\omega \Psi_\epsilon=\bar{A}_\epsilon\Psi_\epsilon-\Psi_\epsilon A_\epsilon$,

\item $|\bar{F}_\epsilon|_{r'}\leq \epsilon$, 

\item $\mid \Psi_\epsilon \mid _{r'} \leq \epsilon^{-\frac{1}{8}}$,

\item and $|Z_\epsilon
-Id|_{r'}\leq 2\epsilon_0^\frac{1}{2}$. 
\end{enumerate}

\bigskip
\noindent Moreover, in dimension 2 or if $G=GL(n,\mathbb{C})$ or $U(n)$, $Z_\epsilon,\bar{A}_\epsilon,\bar{F}_\epsilon$ 
are continuous on $\mathbb{T}^d$. 

\end{thm}

\noindent Property \ref{redtoby} states the reducibility of $\bar{A}_\epsilon$. Theorem \ref{th1} immediately entails the following: 

\begin{thm}\label{cor1}
Let $0<r'<r\leq \frac{1}{2}$, $A\in \mathcal{G}$, $F\in C^\omega_r(\mathbb{T}^d, \mathcal{G})$. 
There is $\epsilon_0<1$ 
depending only on $n,d,\omega,A,r-r'$ such that if

$$|F|_r\leq \epsilon_0$$ 

\noindent 
then for all $\epsilon>0$, there exists ${F}_\epsilon
\in C^\omega_{r'}
(2\mathbb{T}^d,\mathcal{G})$, 
$Z_\epsilon\in C^\omega_{r'}(2\mathbb{T}^d, G)$ and $A_\epsilon\in \mathcal{G}$ 
such that for all $\theta\in 2\mathbb{T}^d$,

$$\partial_\omega Z_\epsilon(\theta)=(A+F(\theta))Z_\epsilon(\theta)-Z_\epsilon(\theta)({A}_\epsilon
+{F}_\epsilon(\theta))
$$

\noindent with $|{F}_\epsilon|_{r'}\leq \epsilon$. 

\end{thm}


\bigskip
\noindent Note that in Theorem \ref{cor1}, we do not have any good estimate of $Z_\epsilon$. Theorem \ref{th1} also holds if one chooses $F$ in a class which is bigger than $C^\omega_r(\mathbb{T}^d,
\mathcal{G})$, 
i.e the class of functions in $C^\omega_r(2\mathbb{T}^d,\mathcal{G})$ 
satisfying some 
"nice periodicity properties" with respect to the matrix $A$.

\bigskip 
\noindent There is a loss of analyticity in this result, but it is arbitrarily small. A result close to Theorem \ref{th1} in the case when $G=GL(n,\mathbb{R})$ had already been proven in \cite{E1} by L.H.Eliasson:

\begin{blank}
{Let $A\in gl(n,\mathbb{R})$ and $F\in C^\omega_r(\mathbb{T}^d, gl(n,\mathbb{R}))$. 
There is $\epsilon_0<1$ 
depending only on $n,d,\kappa,\tau,||A||,r$ such that if $|F|_r\leq \epsilon_0$, then for all $\epsilon>0$, there exists $0<r_\epsilon<r$, $Z_\epsilon\in C^\omega_{r_\epsilon}(2\mathbb{T}^d, GL(n,\mathbb{R}))$
such that for all $\theta\in 2\mathbb{T}^d$,

$$\partial_\omega Z_\epsilon(\theta)=(A+F(\theta))Z_\epsilon(\theta)-Z_\epsilon(\theta)({A}_\epsilon
+{F}_\epsilon(\theta))
$$

\noindent with
$A_\epsilon\in gl(n,\mathbb{R})$, 
${F}_\epsilon\in C^\omega_{r_\epsilon}(2\mathbb{T}^d,gl(n,\mathbb{R}))$ and
$|{F}_\epsilon|_{r_\epsilon}\leq \epsilon$.}
\end{blank}

\noindent Eliasson's theorem merely states almost reducibility in $C^\omega_0(2\mathbb{T}^d,GL(n,\mathbb{R}))$, 
since the sequence $(r_\epsilon)$ might well tend to 0. The achievement of Theorem \ref{th1} is to state almost reducibility in a more general algebraic framework, but also, and mostly, to show that almost reducibility holds in a fixed neighbourhood of a torus even when this torus has dimension greater than 1. This is almost reducibility in a strong sense.

\bigskip 
\noindent Note that, as was the case in \cite{E1}, one cannot avoid to lose periodicity in theorem \ref{th1} if $G$ is a real group with dimension greater than 2. The notion of "nice periodicity properties" that will be given aims at limiting this loss to a period doubling. In comparison with the real framework, the symplectic framework does not introduce any new constraints in the elimination of resonances; therefore there is no more loss of periodicity here than in the case when $G=GL(n,\mathbb{R})$. As before in \cite{C}, a single period doubling is sufficient in the case when $G$ is a real symplectic group.

\bigskip 
\noindent The second part of this paper is dedicated to showing that the same method gives an analogous result for cocycles which are in a Gevrey class (Theorem \ref{th1G'}). 


%






%



\bigskip
\noindent In dimension 2 or if $\mathcal{G}$ is $gl(n,\mathbb{C})$ or $u(n)$, these results can be rephrased as density of reducible cocycles in the neighbourhood of constant cocycles:

\begin{thm}\label{th3}Let $\mathcal{G}=gl(n,\mathbb{C}),u(n), gl(2,\mathbb{R}), sl(2,\mathbb{R})$ or $o(2)$. Let $0<r'<r\leq \frac{1}{2}$ and $A\in \mathcal{G}, F\in C^\omega_r(\mathbb{T}^d,\mathcal{G})$. There is 
$\epsilon_0$ depending only on $r-r',n,d,\omega,A$ such that if

$$|F|_r\leq \epsilon_0$$ 

\noindent then for all 
$\epsilon>0$ there exists $H\in C^\omega_{r'}(\mathbb{T}^d,\mathcal{G})$ which is reducible in $C^\omega_{r'}(\mathbb{T}^d,\mathcal{G})$ and such that

$$|A+F-H|_{r'}\leq \epsilon$$ 
\end{thm}


\noindent A similar result, for smooth cocycles with values in compact Lie groups, was obtained by R. Krikorian in \cite{Kr} (th.5.1.1). For cocycles over a rotation on the circle, analyticity is far better controlled (see for instance \cite{AK06}) since it is then possible to use global methods. In this article, we are considering the case of a torus of arbitrary  dimension. The KAM-type method that is being used here had already given way to full-measure reducibility results for cocycles with values in $SL(2,\mathbb{R})$ (\cite{El92}, \cite{SH06}).

\subsection*{Sketch of the proof and organization of the paper}

The proof of Theorems \ref{th1} and \ref{th3} is a refinement of the method in \cite{E1}; it is based on a KAM scheme. The central idea is to prove an inductive lemma where one conjugates a system which is close to a reducible one to another system which is even closer to something reducible. Iterating this lemma arbitrarily many times, one would then be able to conjugate the initial system to something which is arbitrarily close to a reducible one. An estimate on the reducing transformation would then imply almost reducibility. 
Now  consider a system close to a reducible one; if it is close to a system which can be reduced to a constant part satisfying some non-resonant conditions, then there exists a conjugation which is close to the identity in a good topology taking the first system to something closer to a reducible system. But the constant part might well be too resonant for such a conjugation to exist. In this case, it is possible to remove the resonances in the constant part, but then the conjugation will not stay very close to the identity except if one accepts to give up a lot of regularity. Now we want to avoid this loss of regularity in order to obtain a strong version of almost reducibility. So we will have to improve the step of removing the resonances and use the following two facts: when resonances have been removed up to some order $N$, firstly, the eigenvalues will be so close together that resonances are in fact removed up to an order $RN$ which is much greater than $N$; secondly, the eigenvalues are removed in a durable way, that is, one will not have to remove resonances again until a great number of conjugations is made that will take the cocycle to something much closer to a reducible one. The article is organized as follows:

\bigskip
Section \ref{C2} is dedicated to the proof of the theorem in the analytic case. Here are the main steps of the proof:

\begin{itemize}

\item Removing of the resonances by a map $\Phi$ called a reduction of the eigenvalues at order $R,\bar{N}$ (Proposition \ref{diophan})
for $R,N\in \mathbb{N}\setminus\{0\}$. 

\bigskip
\noindent In dimension 2, $\Phi$ will be such that for all $H$ continuous on 
$\mathbb{T}^d$, $\Phi H\Phi^{-1}$ is continuous on $\mathbb{T}^d$. 

\bigskip
\noindent This step is crucial in the obtention of strong almost reducibility. The reduction of the eigenvalues is defined in a way similar to \cite{E1}, 
however here it will remove resonances up to an order $R\bar{N}$ which is much greater than the value of the parameter $\bar{N}$ appearing in the estimates. 
The parameter 
$R$ will be used to define a map of reduction of the eigenvalues at order $R,\bar{N}$ where $\bar{N}$ does not depend on the loss of analyticity. 
This way, the map of reduction of the eigenvalues 
will stay under control on a neighbourhood of the torus which will not have to fade totally.

\bigskip

\item Resolution of the homological equation (Proposition \ref{homol}): if 
$\tilde{A}$ has a spectrum fulfilling some non-resonance conditions 
and $\tilde{F}$ is a function with nice periodicity properties with respect to $\tilde{A}$, then there exists a solution 
$\tilde{X}$ of equation 

$$\partial_\omega \tilde{X}=[\tilde{A},\tilde{X}]+\tilde{F}^{R\bar{N}};\ \hat{\tilde{X}}(0)=0$$

\noindent having the same periodicity properties as $\tilde{F}$; it takes its values in the same Lie algebra as does $\tilde{F}$. 
Moreover, it can be well controlled by losing some analyticity.

\bigskip

\item Inductive lemma (Proposition \ref{iter3}):
If $\tilde{F}\in C^\omega_{r}(2\mathbb{T}^d,\mathcal{G})$ has some periodicity properties (with respect to $
\tilde{A}$), 
if

$$\partial_\omega \Psi=\bar{A}\Psi-\Psi \tilde{A} $$

\noindent and $\bar{F}=\Psi \tilde{F} \Psi^{-1}$, then there exists
$Z\in C^\omega_{r'}(2\mathbb{T}^d,G)$ such that

\begin{equation}\partial_\omega Z =(\bar{A} +\bar{F} )Z -Z (\bar{A}' 
+\bar{F}' )
\end{equation}

\noindent with $\bar{A}'$ reducible, $\bar{F}'$ is much smaller than $\bar{F}$, 
$Z$ is close to the identity and $\Psi'^{-1}\bar{F}'\Psi'$ has periodicity properties with respect to $A'$ 
which are similar to the properties of $\tilde{F}$.

\bigskip
\noindent The estimate of $\bar{F}'$ depends on $\tilde{F}-\tilde{F}^{R\bar{N}}$, on the reduction of the eigenvalues $\Phi$, 
and on the solution $\tilde{X}$ of the homological equation. 

\bigskip

\item Iteration of the inductive lemma (Theorem \ref{PR}):
We shall iterate Lemma \ref{iter3} so as to obtain estimates of analytic functions on a sequence of neighbourhoods of the torus not tending to 0,
by means of a numerical lemma (Lemma \ref{eps'1}), 
to reduce the perturbation arbitrarily.

\end{itemize}

In section \ref{PRGev}, some lemmas are given (\ref{prelGev}) which show that it is possible to adapt the proof to the Gevrey case; namely, the estimates will be analogous to those that are obtained in the analytic case and so, by slightly modifying the parameters, the argument works in the same way: one obtains analogous reduction of the eigenvalues (\ref{redeigenG}), homological equation (\ref{homoleqG}) and inductive lemmas (\ref{indlemG}).

\subsection*{Notations, further definitions and a general assumption}

\noindent For a function $f\in C^1(2\mathbb{T}^d,gl(n,\mathbb{C}))$, for all $\theta\in 2\mathbb{T}^d$ we will denote by

\begin{equation}\partial_\omega f(\theta)=\frac{d}{dt} f(\theta+t\omega)_{\mid t=0}
\end{equation}

\noindent 
the derivative of $f$ in the direction $\omega$. Denote by $\langle .,. \rangle$ the complex euclidian scalar product, taking it antilinear in the second variable.
For a linear operator $M$, we shall call $M^*$ its adjoint, which is identical to the transpose of $M$ if $M$ is real. Also denote by $M_{\mathcal{N}}$ the nilpotent part of 
$M$, as follows:
let $M=PAP^{-1}$ with $A$ in Jordan normal form, let $A_D$ be the diagonal part of $A$, then $M_{\mathcal{N}}=P(A-A_D)P^{-1}$.
To simplify the writing, if $A: 2\mathbb{T}^d\rightarrow GL(n,\mathbb{C})$, we will 
denote by $A^{-1}$ the map ${\theta\mapsto A(\theta)^{-1}}$. For all $m=(m_1,\dots ,m_d)\in \frac{1}{2}\mathbb{Z}^d$, we shall denote $\mid m\mid = \mid m_1\mid +\dots +\mid m_d\mid$. 
The letter $J$ will stand for matrix $
J=\left(\begin{array}{cc}
0 & -Id\\
Id & 0\\
\end{array}\right)
$.

\bigskip
\noindent \Def A function $f$ is \textit{analytic on an $r$-neighbourhood of the torus} (resp. \textit{double torus)} 
if $f$ is holomorphic on 
$\{x=(x_1,\dots ,x_d)\in \mathbb{C}^d, \sup_j| \operatorname{Im} x_j|<r\}$ and $1$-periodic (resp. $2$-periodic) in $\operatorname{Re} x_j$ for all $1\leq j\leq d$.

\noindent  
For all subset $E$ of $gl(n,\mathbb{C})$, denote by 
$C^\omega_r(\mathbb{T}^d,E)$ the set of functions which are analytic on an 
$r$-neighbourhood of the torus and whose restriction to $\mathbb{R}^d$ takes its values in $E$; 
let $C^\omega_r(2\mathbb{T}^d,E)$ be the set of functions which are analytic on an 
$r$-neighbourhood of the double torus and whose restriction to $2\mathbb{T}^d$ takes its values in $E$. 
For all $f\in C^\omega_r(2\mathbb{T}^d,E)$, denote

\begin{equation}|f|_r=\sup_{|Im x|< r}||f(x)||\end{equation} 

\noindent where $||.||$ stands for the operator norm.

\bigskip 
Let $C^{G,\beta}_r$ be the class of \textit{Gevrey $\beta$ functions with parameter $r$}, i.e functions $f$ satisfying



$$\sum_{\alpha\in \mathbb{N}^d}\frac{r^{\beta \mid \alpha\mid}}{\alpha!^\beta}\sup_{\theta}\mid \mid \partial^\alpha F(\theta) \mid \mid <+\infty$$

\noindent 
Denote by $\mid \mid .\mid \mid_{\beta,r}$ the norm

$$\mid \mid F\mid \mid_{\beta,r}=\sum_{\alpha\in \mathbb{N}^d} \frac{r^{\beta \mid \alpha\mid}}{\alpha!^\beta} \sup_\theta \mid \mid \partial^\alpha F(\theta) \mid \mid$$

\bigskip
To formalize the notion of reducibility, 
we shall introduce an equivalence relation on cocycles. 

\bigskip 
\noindent 
\Def Let $G$ be a Lie group and $\mathcal{G}$ the Lie algebra associated to $G$. Let $r,r'>0$ and $A,B\in C^\omega_r(2\mathbb{T}^d,\mathcal{G})$. 
We say that $A$ and $B$ are \textit{conjugate in $C^\omega_{r'}(2\mathbb{T}^d,G)$} if there exists $Z\in C^\omega_{r'}(2\mathbb{T}^d,G)$ such that
 for all $\theta\in 
2\mathbb{T}^d$,

$$\partial_\omega Z(\theta)=A(\theta)Z(\theta)-Z(\theta)B(\theta)$$

\noindent where $\partial_\omega$ means the derivative in the direction $\omega$. If $B$ is constant in $\theta$, we say that $A$ is \textit{reducible in $C^\omega_{r'}(2\mathbb{T}^d,G)$},
or \textit{reducible by $Z$ to
$B$}. We will use an analogous definition with $C^{G,\beta}$ instead of $C^\omega$.

\bigskip Note that if $X$ is the quasi-periodic cocycle associated to $A$, then the map $A$ is reducible by $Z$ to $B$ if and only if 

\begin{equation}\forall (t,\theta),\ X^t(\theta)=Z(\theta+t\omega)^{-1}e^{tB}Z(\theta)
\end{equation}

\noindent Reducibility is also equivalent to the fact that the map from $2\mathbb{T}^d\times \mathbb{R}^n$ to itself:

\begin{equation}
\left(\begin{array}{c}
\theta\\
v\\
\end{array}\right)
\mapsto
\left(\begin{array}{c}
\theta+\omega\\
X^1(\theta)v\\
\end{array}\right)
\end{equation}

\noindent is conjugate to a map $\chi$ such that

\begin{equation}\frac{d\chi}{d\theta}  \left(\begin{array}{c}
\theta\\
v\\
\end{array}\right)
\equiv \left(\begin{array}{c}
\bar{1}\\
{0}\\
\end{array}\right)\end{equation}

\noindent

\bigskip
\textbf{Assumption:} The frequency $\omega$ is in the diophantine class $DC(\kappa,\tau)$, i.e

\begin{equation}
\forall \ m\in\mathbb{Z}^d\setminus\{0\},\ |\langle m,\omega\rangle|\geq \frac{\kappa}{|m|^\tau}\end{equation}

\noindent where $\kappa,\tau$ are fixed throughout the paper and $0<\kappa<1, \tau\geq \max(1,d-1)$.

\section{Strong almost reducibility for analytic quasi-periodic cocycles}\label{C2}

\subsection{Nice periodicity properties}

\noindent A few definitions will first be given. The notion of "triviality with respect to a decomposition" will make the 
construction of the map of reduction of the eigenvalues easier; the "nice periodicity properties" have been introduced in \cite{E1} and are used in the real case to make sure that only one period doubling will be needed in iterating the inductive lemma. 

\subsubsection{Invariant decompositions}

\noindent 
The set $\mathcal{L}=\{L_1,\dots ,L_R\}$ is called a \textit{decomposition of $\mathbb{C}^n$} if $\mathbb{C}^n=\bigoplus_j L_j$. If $\mathcal{L},\mathcal{L}'$ are decompositions of $\mathbb{C}^n$, then $\mathcal{L}$ is said to be 
\textit{finer than} 
$\mathcal{L}'$ if for all $L\in\mathcal{L}$, there is $L'\in\mathcal{L}'$ such that $L\subset L'$; $\mathcal{L}$ is said \textit{strictly finer than} 
$\mathcal{L}'$ if $\mathcal{L}$ is finer than 
$\mathcal{L}'$ and $\mathcal{L}\neq\mathcal{L}'$.

\bigskip
\noindent \Def Let $A\in gl(n,\mathbb{C})$; then $\mathcal{L}=\{L_1,\dots ,L_s\}$ is an \textit{$A$-decomposition}, or else \textit{$A$-invariant decomposition}, 
if it is a decomposition of $\mathbb{C}^n$ and for all $i$, $AL_i\subset L_i$. Subsets $L_i$ are called \textit{subspaces} of 
$\mathcal{L}$. 

\bigskip 
\noindent Let a Jordan decomposition for $A$ be an $A$-decomposition which is minimal (i.e no finer decomposition is an $A$-decomposition). Then
\begin{itemize} 
\item A matrix might have many Jordan decompositions. For instance, the identity has infinitely many Jordan decompositions. 
\item A decomposition is an $A$-decomposition if and only if it is less fine than some Jordan decomposition for $A$. Therefore, if operators $A$ and $A'$ have a common Jordan 
decomposition, then an $A$-decomposition which is less fine than this common Jordan decomposition is an $A'$-decomposition.
\end{itemize}

\bigskip
\noindent \textbf{Notation:} Let $\mathcal{L}$ be an $A$-decomposition. For all $L\in\mathcal{L}$, denote by 
$\sigma(A_{|L})$ the spectrum of the restriction of $A$ to subspace $L$.

\bigskip
\noindent \Def Let $\kappa'\geq 0$. Let $\mathcal{L}_{A,\kappa'}$ be the unique $A$-decomposition $\mathcal{L}$ such that for all $L\neq L'\in \mathcal{L}$, $\alpha\in\sigma(A_{|L})$ and 
$\beta\in\sigma(A_{|L'})$ 
$\Rightarrow |\alpha-\beta|>\kappa'$ and such that no $A$-decomposition strictly finer than $\mathcal{L}$ has this property.

\bigskip
\noindent \rem For $\kappa'\geq 0$, any Jordan decomposition is finer than $\mathcal{L}_{A,\kappa'}$.

%
%


\bigskip
\noindent \Def Let $\mathcal{L}$ be a decomposition of $\mathbb{C}^n$. For all $u\in\mathbb{C}^n$, there is a unique 
decomposition $u=\sum_{L\in\mathcal{L}}u_L$ such that $u_L\in L$ for all $L\in\mathcal{L}$. 
For all $L\in\mathcal{L}$, the \textit{projection 
on $L$ with respect to $\mathcal{L}$}, denoted by $P^\mathcal{L}_L$, is the map defined by $P^\mathcal{L}_Lu=u_L$.

\bigskip
\noindent 
\rem Let $A\in gl(n,\mathbb{C})$ and $\kappa'>0$. If $\mathcal{L}$ is an $A$-decomposition which is less fine than $\mathcal{L}_{A,\kappa'}$, then one has the following lemma,  which can be found in \cite{E1}, appendix, Lemma A\footnote{Lemma A from \cite{E1} gives in fact an estimate which depends on $\mid \mid A\mid \mid$, but the proof shows clearly that the estimate in fact only depends on $A_{\mathcal{N}}$.}:

\begin{lem}\label{c0}
There is a constant $C_0\geq 1$ depending only on $n$ such that for all subspace $L\in \mathcal{L}$,

\begin{equation}\mid \mid P^{\mathcal{L}}_L\mid \mid \leq C_0\left(\frac{1+\mid \mid A_{\mathcal{N}}\mid \mid }{\kappa'}\right)^{n(n+1)}
\end{equation}

\end{lem}

\noindent In what follows, $C_0$ will always stand for this constant fixed in Lemma \ref{c0}.

\bigskip
\noindent \Def An \textit{$(A,\kappa',\gamma)$-decomposition} is an $A$-decomposition $\mathcal{L}$ such that for all $L\in \mathcal{L}$, the projection on $L$ with respect to $\mathcal{L}$ 
satisfies

\begin{equation}\mid \mid P^{\mathcal{L}}_L\mid \mid \leq C_0\left(\frac{1+\mid \mid A_{\mathcal{N}}\mid \mid }{\kappa'}\right)^{\gamma}
\end{equation}

\bigskip
\noindent \rem 
For $A\in gl(n,\mathbb{C})$, one always has $A=\sum_{L,L'\in\mathcal{L}}P^\mathcal{L}_L AP^\mathcal{L}_{L'}$. 
In particular, if $\mathcal{L}$ is an $A$-decomposition, then $A=\sum_{L\in\mathcal{L}}P^\mathcal{L}_L A P^\mathcal{L}_L$.

\bigskip
\noindent \Defs Let $\mathcal{L}$ be a decomposition. We say that
\begin{itemize}
\item $\mathcal{L}$ is a \textit{real decomposition} if for all $L\in\mathcal{L}$, $\bar{L}\in\mathcal{L}$;

\item $\mathcal{L}$ is a \textit{symplectic decomposition} if it is a decomposition of $\mathbb{C}^{n}$ with even $n$ and for all $L\in \mathcal{L}$, there is a unique $L'
\in \mathcal{L}$ such that $\langle L,JL'\rangle \neq 0$; 
\item $\mathcal{L}$ is a \textit{unitary decomposition} if for all $L\neq L'\in \mathcal{L}$, 
$\langle L,L'\rangle =0$.

\end{itemize}

\bigskip
\noindent \rem 
\begin{itemize}
\item If $A$ is a real matrix, then for all $\kappa'\geq 0$,  $\mathcal{L}_{A,\kappa'}$ is a real decomposition. 
\item For all $L$, there is at least one $L'$ such that $\langle L,JL'\rangle \neq 0$. This comes from the fact that the symplectic form $\langle .,J.\rangle$ is non-degenerate.
\item If $A\in sp(n,\mathbb{R})$, then any $A$-decomposition $\mathcal{L}$ which is less fine than $\mathcal{L}_{A,0}$ is a real and symplectic decomposition. To see this, let $L,L'\in \mathcal{L}$ such that $\langle L,JL'\rangle \neq 0$; let $v\in L,v'\in L'$ be eigenvectors of $A$ such that $\langle v,Jv'\rangle \neq 0$ 
and $\lambda,\lambda'$ their associated eigenvalues. Then

\begin{equation}\nonumber\lambda \langle v,Jv'\rangle= \langle Av,Jv'\rangle =\langle v,A^*Jv'\rangle 
=-\langle v,JAv'\rangle =-\bar{\lambda}'\langle v,Jv'\rangle 
\end{equation}

\noindent and since $\langle v,Jv'\rangle \neq 0$, then $\lambda=-\bar{\lambda}'$.

\item If $A\in U(n)$, then any decomposition which is less fine than $\mathcal{L}_{A,0}$ is unitary. 

\item If $\mathcal{L}$ is unitary, then for every $L\in \mathcal{L}$, $P^{\mathcal{L}}_L$ is an orthogonal projection so 

\begin{equation}\nonumber\mid \mid P^{\mathcal{L}}_L \mid \mid \leq 1
\end{equation}

\end{itemize}

\subsubsection{Triviality and nice periodicity properties with respect to a decomposition}

\noindent \Def Let $\mathcal{L}$ be a decomposition of $\mathbb{C}^n$. 
We say a map $\Psi$ 
is \textit{trivial with respect to $\mathcal{L}$} if there exist 
${\{m_L,\ L\in \mathcal{L}\}\subset \frac{1}{2}\mathbb{Z}^d}$ such that for all $\theta\in 2\mathbb{T}^d$,

\begin{equation}
\Psi(\theta) =\sum_{L\in\mathcal{L}}e^{2i\pi\langle m_L,\theta\rangle}P_L^\mathcal{L}\end{equation}

\bigskip
\noindent We say that the function $\Psi$ is \textit{trivial} if there exists a decomposition $\mathcal{L}$ such that $\Psi$ is trivial with respect to $\mathcal{L}$. 

\bigskip
\noindent \rem 
\begin{itemize}

\item If $\Psi$ is trivial 
with respect to $\mathcal{L}$ and $\mathcal{L}'$ is finer than $\mathcal{L}$, 
then $\Psi$ is trivial with respect to $\mathcal{L}'$. 
\item If $\Phi,\Psi:2\mathbb{T}^d\rightarrow GL(n,\mathbb{C})$ are trivial with respect to 
$\mathcal{L}$, then the product $\Phi\Psi$ is trivial with respect to $\mathcal{L}$.
\item If $\Phi$ is trivial with respect to an $A$-decomposition $\mathcal{L}$, then for all $\theta\in
2\mathbb{T}^d$, $[A,\Phi(\theta)]=0$.

\end{itemize}

\begin{lem}\label{réel} Let $\mathcal{L}$ be a real decomposition of $\mathbb{C}^n$, ${\{m_L,\ L\in \mathcal{L}\}\subset \frac{1}{2}\mathbb{Z}^d}$ 
and $\Psi$ defined by

\begin{equation}
\label{psitriv}
\Psi(\theta) =\sum_{L\in\mathcal{L}}e^{2i\pi\langle m_L,\theta\rangle}P_L^\mathcal{L}\end{equation}

\noindent Then $\Psi$ is real if and only if for all $L$, $m_L=-m_{\bar{L}}$. Moreover, if $\Psi$ is real, then $\Psi$ takes its values in $SL(n,\mathbb{R})$.

\end{lem}

\noindent \dem Assume that for all $L\in \mathcal{L}$, $m_L=-m_{\bar{L}}$. Let $u\in \mathbb{R}^n$. Then

\begin{equation}\nonumber\overline{\Psi(\theta)u}=\sum_{L\in \mathcal{L}}e^{2i\pi \langle -m_L,\theta\rangle} \overline{P^{\mathcal{L}}_Lu}
=\sum_{L\in \mathcal{L}}e^{2i\pi \langle m_{\bar{L}},\theta\rangle} P^{\mathcal{L}}_{\bar{L}}u
=\Psi(\theta)u
\end{equation}

\noindent 
so $\Psi(\theta)$ is real. 

\bigskip
\noindent Now suppose that $\Psi$ is real. Then for all $\theta$,

\begin{equation}\nonumber
\sum_{L\in\mathcal{L}}e^{2i\pi\langle m_L,\theta\rangle}P_L^\mathcal{L}
=\sum_{L\in\mathcal{L}}e^{2i\pi\langle -m_L,\theta\rangle}\overline{P_L^\mathcal{L}}
=\sum_{L\in\mathcal{L}}e^{2i\pi\langle -m_L,\theta\rangle}P_{\bar{L}}^\mathcal{L}
\end{equation}

\noindent so $m_L=-m_{\bar{L}}$. 

\bigskip
\noindent Suppose $\Psi$ is real; then for all $L$, $m_L=-m_{\bar{L}}$ so $\Psi(\theta)$ is the exponential of a trace-zero matrix, so it has determinant 1.
$\Box$

\bigskip
\noindent \rem Any map which is trivial with respect to a unitary decomposition is unitary: let $\mathcal{L}$ be a unitary decomposition, let $\Phi$ be trivial with respect to $\mathcal{L}$ and let $L,L'\in \mathcal{L}$. 
Then for all $u\in \mathcal{L},v\in \mathcal{L}'$,

\begin{equation}\nonumber\langle \Phi(\theta) u,\Phi(\theta) v\rangle 
= \langle e^{2i\pi \langle m_L,\theta \rangle} u,e^{2i\pi \langle m_{L'},\theta \rangle} v\rangle 
=\langle  u, v\rangle 
\end{equation}

\begin{lem}\label{phisp}Let $\mathcal{L}$ be a real and symplectic decomposition and $\{m_L,L\in\mathcal{L}\}$ be a family of elements of $\frac{1}{2}\mathbb{Z}^d$. 
Let $\Psi=\sum_{L\in \mathcal{L}}e^{2i\pi \langle m_L,.\rangle}P^{\mathcal{L}}_L$.
Then $\Psi $ takes its values in $Sp(n,\mathbb{R})$ if and only if
\begin{itemize}
\item for all $L$, $m_L=-m_{\bar{L}}$
\item and if $\langle L,JL'\rangle \neq 0$, 
then $m_L=m_{L'}$. 
\end{itemize}

\end{lem}

\noindent \dem By Lemma \ref{réel}, $\Psi$ is real if and only if for all $L$, $m_L=-m_{\bar{L}}$. Assume now $\Psi$ is real. 

\bigskip
\noindent We show first that if for all $L,L'\in \mathcal{L}$, $\langle L,JL'\rangle \neq 0\Rightarrow m_L=m_{L'}$, then $\Psi$ takes its values in 
$Sp(n,\mathbb{R})$. 
Let $u,v\in \mathbb{R}^{n}$. Then

\begin{equation}\nonumber\langle u, \Psi(\theta)^*J\Psi(\theta)v\rangle 
=\langle  \Psi(\theta)u,J\Psi(\theta)v\rangle 
=\sum_{L} e^{2i\pi \langle m_L-m_{M(L)},\theta\rangle} \langle P^{\mathcal{L}}_L u,J P^{\mathcal{L}}_{M(L)}v\rangle
\end{equation}

\noindent where $M(L)$ stands for the unique subspace such that $\langle L,JM(L)\rangle \neq 0$. 
Assume that if $\langle L,JL'\rangle \neq 0$, 
then $m_L=m_{L'}$. 
This implies that

\begin{equation}\nonumber\langle u, \Psi(\theta)^*J\Psi(\theta)v\rangle 
=\sum_{L} \langle P^{\mathcal{L}}_L u,J P^{\mathcal{L}}_{M(L)}v\rangle
=\langle u,Jv\rangle
\end{equation}

\noindent so $\Psi(\theta)\in Sp(n,\mathbb{R})$. 

\bigskip
\noindent Now we will show that if $\Psi(\theta)\in Sp(n,\mathbb{R})$ and if $\langle L,JL'\rangle \neq 0$, 
then $m_L=m_{L'}$.  Suppose $\Psi(\theta)\in Sp(n,\mathbb{R})$. For any two vectors $u,v$,

\begin{equation}\nonumber\langle u,Jv\rangle= \langle u, \Psi(\theta)^*J\Psi(\theta)v\rangle 
=\langle  \Psi(\theta)u,J\Psi(\theta)v\rangle 
\end{equation}

\noindent If $u\in L$ and $v\in m(L)$ satisfy $\langle u,Jv\rangle \neq 0$, then

\begin{equation}\nonumber\langle u,Jv\rangle
=\langle  \Psi(\theta)u,J\Psi(\theta)v\rangle 
=e^{2i\pi \langle m_L-m_{M(L)},\theta\rangle} \langle  u,J v\rangle
\end{equation}

\noindent so $m_L=m_{M(L)}$.
 $\Box$

\bigskip
\noindent We will now define the periodicity properties. 

\bigskip
\noindent \Def Let $\mathcal{L}$ be a decomposition of $\mathbb{C}^n$. 
We say that $F\in C^0(2\mathbb{T}^d,
gl(n,\mathbb{R}))$ \textit{has nice periodicity properties with respect to $\mathcal{L}$} if there exists a map $\Phi$ which is trivial with respect to $\mathcal{L}$ and such that 
$\Phi^{-1}F\Phi$ is continuous on $\mathbb{T}^d$.

\noindent To make the family $(m_L)$ explicit, we say that $F$ \textit{has nice periodicity properties with respect to $\mathcal{L}$ and $(m_L)$}.

\bigskip
\noindent \rem 
\begin{itemize}
\item If $F\in C^0(2\mathbb{T}^d,gl(n,\mathbb{R}))$ has nice periodicity properties with respect to a decomposition $\mathcal{L}$ and $\Phi$ is 
trivial with respect to $\mathcal{L}$, 
then $\Phi F \Phi^{-1}$ has nice periodicity properties with respect to $\mathcal{L}$.

\item 
If $\mathcal{L}'$ is a decomposition of $\mathbb{C}^n$ which is finer than $\mathcal{L}$ and $F$ has nice periodicity properties with respect to $\mathcal{L}$,then 
$F$ has nice periodicity properties with respect to $\mathcal{L}'$. 

\item Let $\mathcal{L}$ be a decomposition of $\mathbb{C}^n$ and $(m_L)_{L\in \mathcal{L}}$ be a family of elements of $\frac{1}{2}\mathbb{Z}^d$. 
If $F_1,F_2\in C^0(2\mathbb{T}^d, gl(n,\mathbb{R}))$ have nice periodicity properties with respect to $\mathcal{L}$ and $(m_L)$, then the product $F_1F_2$ has nice periodicity properties with respect to $\mathcal{L}$ and $(m_L)$.

\end{itemize}

\subsection{Removing the resonances}\label{elimres}

\noindent In the following we will have to solve a homological equation and estimate the solution on a neighbourhood of the torus; in order to have a sufficient estimate, one will assume that the coefficients of the equation satisfy some diophantine conditions:

\bigskip
Let $A\in gl(n,\mathbb{R})$ and $0<\kappa'<1$. Let $N\in\mathbb{N}$. 

\bigskip
\Def Let $z\in\mathbb{C}, \nu\in\{1,2\}$. We say that $z$ is \textit{diophantine modulo $\nu$ with respect to $\omega$, with constant $\kappa'$, exponent $\tau$ and 
order $N$} if for every $m\in \frac{1}{\nu}\mathbb{Z}^d$ such that $0<|m|\leq N$, 

\begin{equation}\label{mnonnul}|z-2i\pi\langle m,\omega\rangle|\geq \frac{\kappa'}{|m|^\tau}
\end{equation}

\noindent This property will be denoted by

\begin{equation}z\in DC_{{\omega},{\nu}}^N(\kappa',\tau)\end{equation}

\noindent Note that

\begin{equation}DC_{\omega,2}^N(\kappa',\tau)\subset DC_{\omega,1}^N(\kappa',\tau)
\end{equation}

\noindent and that every real number $z$ is in $DC_{\omega,2}^N(\frac{\kappa}{2^\tau},\tau)$ since for all $m\in\frac{1}{2}\mathbb{Z}^d$,

\begin{equation}|z-2i\pi \langle m,\omega\rangle|=\left(|z|^2+(2\pi |\langle m,\omega\rangle|)^2\right)^\frac{1}{2}\geq \frac{\pi \kappa}{|2m|^\tau}\geq \frac{\kappa}{|2m|^\tau}
\end{equation} 

\rem In the definition above, the condition is required only for non vanishing $m$, so \eqref{mnonnul} has a meaning.

\bigskip
\Def $A$ is said to have \textit{$DC^N_\omega(\kappa',\tau)$ spectrum} if 

\begin{equation}\left\{\begin{array}{c}
\forall \alpha,\beta
\in \sigma(A),\ \alpha-\beta\in DC_{\omega,1}^N(\kappa',\tau)\\
\forall \alpha,\beta
\in \sigma(A),\ \alpha\neq \bar{\beta}\Rightarrow \alpha-\beta\in DC_{\omega,2}^N(\kappa',\tau)
\end{array}\right.\end{equation}

\bigskip
\noindent Let $N\in\mathbb{N}$. 
Let $A$ in a Lie algebra $\mathcal{G}$. 
The aim is to show that there exists $\kappa'>0$, $\tilde{A}\in \mathcal{G}$ 
such that $\tilde{A}$ has $DC_\omega^N(\kappa',\tau)$ spectrum and $A$ and $\tilde{A}$ are conjugate (in the acception of cocycles, following the definition given in the introduction). To achieve this, one has to find a family $(m_1,\dots ,m_n)$ satisfying

\begin{equation}\left\{\begin{array}{c}
 \forall \ \alpha_j,\alpha_k\in\sigma(A), 
\ \alpha_j-\alpha_k+2i\pi\langle
m_j-m_k,\omega\rangle\in DC_{\omega,1}^N(\kappa',\tau)\\
\forall \ \alpha_j,\alpha_k\in\sigma(A), 
\ \alpha_j\neq \bar{\alpha}_k\Rightarrow
\alpha_j-\alpha_k+2i\pi\langle
m_j-m_k,\omega\rangle\in DC_{\omega,2}^N(\kappa',\tau)
\end{array}\right.\end{equation}

\noindent We shall construct the so-called map of reduction of the eigenvalues $\Phi$ conjugating (in the sense of cocycles) $A$ to the matrix obtained from $A$ by substituting an eigenvalue $\alpha_j$ by $\alpha_j+2i\pi\langle m_j,\omega\rangle$, 
then we will prove that $\Phi$ is
$G$-valued.

\subsubsection{Diophantine conditions}

\begin{lem}\label{m_j} 
Let $\{\alpha_1,\dots ,\alpha_n\}\subset \mathbb{C}$. 
Let $\tilde{N}\in\mathbb{N}$ and $\kappa'\leq\frac{\kappa}{n(8\tilde{N})^\tau}$. 
There exists 
$m_1,\dots ,m_n\in\frac{1}{2}\mathbb{Z}^d$ such that $\sup_j |m_j|\leq \tilde{N}$, and such that letting for all $j$, $\tilde{\alpha}_j=\alpha_j-2i\pi\langle m_j,\omega\rangle$, then 

\begin{equation}\label{conjugué}  \{\alpha_1,\dots ,\alpha_n\}=\overline{\{\alpha_1,\dots ,\alpha_n\}}\Rightarrow \forall j,k,\ \alpha_j=\bar{\alpha}_k\Rightarrow m_j=-m_k\end{equation}

\begin{equation}\label{sl2C} n=2, \alpha_2=-\alpha_1 \Rightarrow m_1=-m_2
\end{equation}

\begin{equation}\label{oppconjugué}\forall j,k,\  \alpha_j=-\bar{\alpha}_k\Rightarrow m_j=m_k\end{equation}

\begin{equation}\label{prox} \forall j,k,\  |\alpha_j-\alpha_k|\leq\kappa'\Rightarrow 
m_j=m_k\end{equation}

\begin{equation}\label{décr}\forall j,\  |Im \tilde{\alpha}_j|\leq |Im \alpha_j|\end{equation}

\begin{equation}\label{mj1}
\forall j,k,\   \alpha_j= \bar{\alpha}_k\Rightarrow 
\tilde{\alpha}_j-\tilde{\alpha}_k\in DC_{\omega,1}^{\tilde{N}}
(\kappa',\tau)\end{equation}

\noindent and

\begin{equation}\label{mj2} 
\forall j,k,\ \alpha_j\neq \bar{\alpha}_k\Rightarrow 
\tilde{\alpha}_j-\tilde{\alpha}_k\in DC_{\omega,2}^{\tilde{N}}({\kappa'},\tau)
\end{equation}

\noindent and such that if not all $m_j$ vanish, then there exist $j,k$ such that

\begin{equation}\label{décr'}| {\alpha}_j- {\alpha}_k|\geq \kappa',\ | \tilde{\alpha}_j- \tilde{\alpha}_k|
<\kappa'\end{equation}

\noindent Moreover, there exist $m_1,\dots m_n\in \mathbb{Z}^d$, with $|m_j|\leq \tilde{N}$ for all $j$, fulfilling conditions \eqref{oppconjugué}, \eqref{prox}, \eqref{décr}, such that 

\begin{equation}\label{mj1C}
\forall j,k,\ \tilde{\alpha}_j-\tilde{\alpha}_k\in DC_{\omega,1}^{\tilde{N}}
(\kappa',\tau)\end{equation}

\noindent and such that if not all $m_j$ vanish, then there exist $j,k$ such that \eqref{décr'} holds.

\end{lem}

\dem We shall proceed in two steps. The first step consists in removing resonances which might occur between two eigenvalues whose imaginary parts are nearly opposite to each other. Once this first lot of resonances is removed, the second step consists in removing the resonances which might occur between two eigenvalues whose imaginary parts are far from opposite.

\bigskip
\noindent $\bullet$ Let $1\leq j\leq n$. Suppose that there is an $m\in \mathbb{Z}^d,0<\mid m\mid \leq \tilde{N}$ such that

\begin{equation}\nonumber\mid 2Im \alpha_j-2\pi \langle m,\omega\rangle \mid <\frac{\kappa'}{\mid m\mid ^\tau}
\end{equation}

then let $\alpha_j'=\alpha_j-2i\pi \langle \frac{m}{2},\omega\rangle$. Otherwise, let $\alpha_j'=\alpha_j$. Note that if $ |\alpha_j-\alpha_k|\leq\kappa'$ and if there exist $m_j\neq m_k$ such that

\begin{equation}\nonumber\mid 2Im \alpha_j-2\pi \langle m_j,\omega\rangle \mid <\frac{\kappa'}{\mid m_j\mid ^\tau};\ \mid 2Im \alpha_k-2\pi \langle m_k,\omega\rangle \mid <\frac{\kappa'}{\mid m_k\mid ^\tau}
\end{equation}

\noindent then

\begin{equation}\nonumber\begin{split}\mid 2i\pi \langle m_j-m_k,\omega\rangle \mid &
\leq \frac{\kappa}{\mid m_j-m_k\mid ^\tau}
\end{split}\end{equation}

\bigskip 
\noindent which is impossible since $\omega$ is diophantine. Therefore conditions \eqref{conjugué} to \eqref{mj1} hold with $\alpha_j'=\tilde{\alpha}_j$ and $m_j$ such that $ \alpha_j-{\alpha}'_j=2i\pi \langle m_j,\omega \rangle $.


\bigskip
\noindent $\bullet$ 
Let $I_{-r},\dots ,I_r$ be the finest partition of $\{1,\dots ,n\}$ such that 

\begin{equation}\nonumber\mid Im(\alpha'_j-\alpha'_k)\mid \leq \kappa' \Rightarrow \exists -r\leq r'\leq r \mid j,k\in I_{r'}
\end{equation}

\noindent and choose the indices in such a way that

\begin{equation}\nonumber r' < r'' \Rightarrow \forall {j\in I_{r'}}, \forall {k\in I_{r''}} , Im \alpha'_j \leq  Im \alpha'_k
\end{equation}

\noindent Note that $I_0$ might be empty. 
We will proceed by induction on $r'$ to prove the following property $\mathcal{P}(r')$: 
\begin{blank}
There are $m'_1,m'_{-1},\dots ,m'_{r'},m'_{-r'}\in \mathbb{Z}^d$ with $\sup_{\mid j\mid \leq r'} |m'_j|\leq \tilde{N}$ such that properties \eqref{conjugué} to \eqref{mj2} hold for all $-r'\leq r_1,r_2\leq r', j\in I_{r_1},k\in I_{r_2}$ with $m'_j$ instead of $m_j$ and 
$\alpha_j'$ instead of $\alpha_j$.

\end{blank}

\bigskip
$\bullet$ Case $r'=0$: if $I_0$ is empty, then $\mathcal{P}(0)$ trivially holds. Assume $I_0$ is non empty. Then for all $j,k\in I_0$ and all $m\in \frac{1}{2}\mathbb{Z}^d$ such that $0<\mid m\mid \leq \tilde{N}$, 

\begin{equation}\nonumber\mid \alpha'_j-\alpha'_k-2i\pi \langle m,\omega\rangle \mid 
\geq \mid Im(\alpha'_j-\alpha'_k)-2\pi \langle m,\omega\rangle \mid 
\geq \frac{\kappa}{\mid m\mid ^\tau}-n\kappa'\geq \kappa'
\end{equation}

\noindent so $\alpha'_j-\alpha'_k \in DC^{\tilde{N}}_{\omega,2}(\kappa',\tau)$ and $\mathcal{P}(0)$ holds true.

\bigskip
$\bullet$ Let $r'\leq r-1$. Assume $\mathcal{P}(r')$ holds. Consider 
$I_{r'+1}$ and $I_{-r'-1}$. There are two possible cases.

\begin{itemize}

\item There exist $-r'\leq r''\leq r',j\in I_{r''},k\in I_{r'+1}$ and $m\in \mathbb{Z}^d$ such that $\mid m\mid \leq \tilde{N}$ and 

\begin{equation}\nonumber\mid \alpha'_j-\alpha'_k-2i\pi \langle m_{r''}+m,\omega\rangle \mid < \frac{\kappa'}{\mid m \mid ^\tau}
\end{equation}

\item The case above does not hold.

\end{itemize}

\bigskip
\noindent 
In the first case, let $m'_{r'+1}=m=-m'_{-r'-1}$. In the second case, let 
$m'_{r'+1}=m'_{-r'-1}=0$. 

\bigskip
\noindent 
Now $m'_{r'+1}$ and $m'_{-r'-1}$ are independent from $j,k$. To see this, suppose there are $j_1,j_2\in I_{r_1}, 
k_1,k_2\in I_{r_2}, m_1\neq m_2\in \mathbb{Z}^d$ such that for $l=1,2$,

\begin{equation}\nonumber\mid \alpha'_{j_l}-\alpha'_{k_l}-2i\pi \langle m_l,\omega\rangle \mid < \frac{\kappa'}{\mid m_l \mid ^\tau}
\end{equation}

\noindent 
Then

\begin{equation}\nonumber\begin{split}\mid 2\pi\langle m_1-m_2,\omega\rangle \mid 
\leq \frac{\kappa}{\mid m_1-m_2\mid ^\tau}
\end{split}\end{equation}

\noindent which is impossible. Therefore $\mathcal{P}(r'+1)$ holds true.

\bigskip
\noindent $\bullet$ 
Once $m'_1,\dots, m'_r,m'_{-1},\dots,m'_{-r}\in \mathbb{Z}^d$ are defined, conditions \eqref{conjugué} to \eqref{mj2} hold with, for all $j\in I_{r'}, \tilde{\alpha}_j=\alpha'_j-2i\pi \langle m'_{r'},\omega\rangle$ and $m_j$ such that $ \alpha_j-\tilde{\alpha}_j=2i\pi \langle m_j,\omega \rangle $. Condition \eqref{décr'} is obvious by construction. 
 
 \bigskip
 \noindent $\bullet$ By proceeding only with the second step, one gets $m_1,\dots m_n\in \mathbb{Z}^d$, with $|m_j|\leq \tilde{N}$ for all $j$, satisfying conditions \eqref{oppconjugué}, \eqref{prox}, \eqref{décr}, such that

\begin{equation}\nonumber\label{mj1C'}
\forall j,k,\ \tilde{\alpha}_j-\tilde{\alpha}_k\in DC_{\omega,1}^{\tilde{N}}
(\kappa',\tau)\end{equation}

\noindent and such that if not all $m_j$ vanish, then there are $j,k$ such that \eqref{décr'} holds true.
 $\Box$

\begin{lem}\label{res} 
Let $\{\alpha_1,\dots ,\alpha_n\}\subset \mathbb{C}$. 
For every $R,N\in\mathbb{N},N\geq 2,R\geq 1$, there exists $\bar{N}\in [N, R^{\frac{1}{2}n(n-1)}N]$ and $m_1,\dots ,m_n\in\frac{1}{2}\mathbb{Z}^d$ with

\begin{equation}\sup_j |m_j|\leq 2\bar{N}\end{equation}

\noindent 
such that letting $\tilde{\alpha}_j=\alpha_j-2i\pi\langle m_j,\omega\rangle$ and

\begin{equation}\kappa''= \frac{\kappa}{n(8R^{\frac{1}{2}n(n-1)+1}{N})^{\tau}}
\end{equation}

\noindent conditions \eqref{conjugué} to \eqref{décr} of Lemma \ref{m_j} hold for $\kappa'=\kappa''$, and such that 

\begin{equation}\label{res2a}
\forall j,k,\ \tilde{\alpha}_j-\tilde{\alpha}_k \in DC_{\omega,1}^{R\bar{N}}
(\kappa'',\tau)\end{equation}

\noindent and

\begin{equation}\label{res2b}\forall j,k,\ 
\alpha_j\neq \bar{\alpha}_k\Rightarrow 
\tilde{\alpha}_j-\tilde{\alpha}_k \in DC_{\omega,2}^{R\bar{N}}(\kappa'',\tau)
\end{equation}

\noindent Moreover, there exist $m_1,\dots m_n\in \mathbb{Z}^d$ with $|m_j|\leq \bar{N}$ for all $j$ such that conditions \eqref{oppconjugué}, \eqref{prox}, \eqref{décr} and \eqref{res2a} hold true.

\end{lem}

\dem If $\alpha_j$ satisfy for all $j,k$

\begin{equation}\label{condjk}\left\{\begin{array}{c}
\alpha_j= \bar{\alpha}_k\Rightarrow  \alpha_j-\alpha_k
\in DC_{\omega,1}^{RN}
(\kappa'',\tau)\\
\alpha_j\neq \bar{\alpha}_k\Rightarrow 
\alpha_j-\alpha_k
 \in DC_{\omega,2}^{RN}(\kappa'',\tau)
\end{array}\right.
\end{equation}

\noindent then we are done with $\bar{N}=N$ and $m_1=\dots =m_n=0$. 

Suppose \eqref{condjk} does not hold. Then apply Lemma \ref{m_j} with $\tilde{N}=RN, \kappa'=\kappa''$ to get $m_1^1,\dots ,m^1_n$ such that

\begin{equation}\left\{\begin{array}{c}
\forall j,k, \ \alpha_j=\bar{\alpha}_k\Rightarrow m^1_j=-m^1_k\\
\forall j,k, \ \alpha_j=-\bar{\alpha}_k\Rightarrow m^1_j=m^1_k\\
\forall j,k, \ |\alpha_j-\alpha_k|\leq\kappa''\Rightarrow 
m^1_j=m^1_k\\
\forall j,\ |Im \alpha_{j}-2i\pi \langle m^1_{j}, \omega\rangle |\leq |Im \alpha_j|
\end{array}\right.\end{equation}

\noindent and

\begin{equation}\label{condjk1}\left\{\begin{array}{c}
\alpha_j= \bar{\alpha}_k\Rightarrow  \alpha_j-\alpha_k-2i\pi\langle m^1_j-m^1_k,\omega\rangle 
\in DC_{\omega,1}^{RN}
(\kappa'',\tau)\\
\alpha_j\neq \bar{\alpha}_k\Rightarrow 
\alpha_j-\alpha_k-2i\pi\langle m^1_j-m^1_k,\omega\rangle
 \in DC_{\omega,2}^{RN}(\kappa'',\tau)
\end{array}\right.
\end{equation}

\noindent and such that there exist $j_1,k_1$ satisfying $\mid Im (\alpha_{j_1}-\alpha_{k_1})-2i\pi \langle m^1_{j_1}-m^1_{k_1}, \omega\rangle \mid <\kappa''$.

\bigskip
\noindent Assume there are $m_1^r,\dots ,m^r_n$ such that 
$\sup |m_j^{r}|\leq (R+R^2+\dots + R^{r})N$ and that for all $j,k$,

\begin{equation}\left\{\begin{array}{c}
\forall j,k, \ \alpha_j=\bar{\alpha}_k\Rightarrow m^r_j=-m^r_k\\
\forall j,k, \ \alpha_j=-\bar{\alpha}_k\Rightarrow m^r_j=m^r_k\\
\forall j,k, \ |\alpha_j-\alpha_k|\leq\kappa''\Rightarrow 
m^r_j=m^r_k\\
\forall j,\ |Im \alpha_{j}-2i\pi \langle m^r_{j}, \omega\rangle |\leq |Im \alpha_j|
\end{array}\right.\end{equation}

\noindent and

\begin{equation}\label{condjkr}\left\{\begin{array}{c}
\alpha_j= \bar{\alpha}_k\Rightarrow  \alpha_j-\alpha_k-2i\pi\langle m^r_j-m^r_k,\omega\rangle 
\in DC_{\omega,1}^{R^rN}
(\kappa'',\tau)\\
\alpha_j\neq \bar{\alpha}_k\Rightarrow 
\alpha_j-\alpha_k-2i\pi\langle m^r_j-m^r_k,\omega\rangle
 \in DC_{\omega,2}^{R^rN}(\kappa'',\tau)
\end{array}\right.
\end{equation}

\noindent and suppose there exist distinct $(j_1,k_1),\dots ,(j_r,k_r)$ such that for all $l\leq r$, 

\begin{equation}\mid Im \alpha_{j_l}-Im \alpha_{k_l}-2i\pi \langle m^r_{j_l}-m^r_{k_l}, \omega\rangle \mid < 
\kappa''
\end{equation}

\noindent 
If moreover one has for all $j,k$

\begin{equation}\label{stopr}\left\{\begin{array}{c}
\alpha_j= \bar{\alpha}_k\Rightarrow  \alpha_j-\alpha_k-2i\pi\langle m^r_j-m^r_k,\omega\rangle 
\in DC_{\omega,1}^{R^{r+1}N}
(\kappa'',\tau)\\
\alpha_j\neq \bar{\alpha}_k\Rightarrow 
\alpha_j-\alpha_k-2i\pi\langle m^r_j-m^r_k,\omega\rangle
 \in DC_{\omega,2}^{R^{r+1}N}(\kappa'',\tau)
\end{array}\right.
\end{equation}

\noindent 
then the process ends with $\bar{N}=R^rN$ and $m_j=m^r_j$ since it is true that

\begin{equation}\mid m^r_j\mid \leq  (R+R^2+\dots +R^r)N\leq R^r N\frac{1-\frac{1}{R^r}}{1-\frac{1}{R}}\leq 2R^rN
\end{equation}

\noindent 
Otherwise, iterate once more Lemma \ref{m_j} with $\tilde{N}=R^{r+1}N$ and $\alpha_j-2i\pi \langle m_j^r,\omega\rangle $ in place of $\alpha_j$ to get $m_1^{r+1},\dots ,m_n^{r+1}$ such that 
$\sup |m_j^{r+1}|\leq (R+R^2+\dots +R^{r+1})N$ and for all $j,k$,

\begin{equation}\left\{\begin{array}{c}
\forall j,k, \ \alpha_j=\bar{\alpha}_k\Rightarrow m^{r+1}_j=-m^{r+1}_k\\
\forall j,k, \ \alpha_j=-\bar{\alpha}_k\Rightarrow m^{r+1}_j=m^{r+1}_k\\
\forall j,k, \ |\alpha_j-\alpha_k|\leq\kappa''\Rightarrow 
m^{r+1}_j=m^{r+1}_k\\
\forall j,\ |Im \alpha_{j}-2i\pi \langle m^{r+1}_{j}, \omega\rangle |\leq |Im \alpha_j|
\end{array}\right.\end{equation}

\noindent and

\begin{equation}\left\{\begin{array}{c} 
\alpha_j= \bar{\alpha}_k\Rightarrow 
\alpha_j-\alpha_k-2i\pi\langle m^{r+1}_j-m^{r+1}_k,\omega\rangle \in DC_{\omega,1}^{R^{r+1}N}
(\kappa'',\tau)\\
\alpha_j\neq \bar{\alpha}_k\Rightarrow 
\alpha_j-\alpha_k-2i\pi\langle m^{r+1}_j-m^{r+1}_k,\omega\rangle \in DC_{\omega,2}^{R^{r+1}N}({\kappa''}
,\tau)
\end{array}\right.
\end{equation}

\noindent and that there exist distinct $(j_1,k_1),\dots ,(j_{r+1},k_{r+1})$ such that for all $l\leq r+1$,

\begin{equation}\mid Im \alpha_{j_l}- Im \alpha_{k_l}-2i\pi \langle m^{r+1}_{j_l}-m^{r+1}_{k_l}, \omega\rangle \mid <\kappa''
\end{equation}

\noindent Therefore, for all $1\leq l\leq r+1$,

\begin{equation}|\alpha_{j_l}-\alpha_{k_l}-2i\pi\langle m^{r+1}_{j_l}-m^{r+1}_{k_l},\omega\rangle |
\ < \kappa''\end{equation}

\noindent This implies that for all $m\in \frac{1}{2}\mathbb{Z}^d$ such that $0<|m|\leq R\bar{N}$ and for all $l,1\leq l\leq r+1$,

\begin{equation}|\alpha_{j_l}-\alpha_{k_l}-2i\pi\langle m_{j_l}^{r+1}-m_{k_l}^{r+1},\omega\rangle -2i\pi\langle m,\omega\rangle|
\geq \frac{\kappa}{2^{\tau+1}(R\bar{N})^\tau}-\kappa''\geq \kappa''
\end{equation}

\noindent so for all $l\leq r+1$,

\begin{equation}\alpha_{j_l}-\alpha_{k_l}-2i\pi\langle m^{r+1}_{j_l}-m^{r+1}_{k_l},\omega\rangle 
\in 
DC^{R\bar{N}}_{\omega,2}
(\kappa'',
\tau)\end{equation}

\noindent 
Therefore, after $\bar{r}\leq \frac{n(n-1)}{2}$ steps, one gets conditions \eqref{res2a} and \eqref{res2b} with $m_j=m^{\bar{r}}_j$ and $\tilde{\alpha}_j=\alpha_j
-2i\pi \langle m_j,\omega\rangle$ and $|Im \alpha_j-2i\pi\langle m_j,\omega\rangle|\leq |Im \alpha_j|$. It is true that $\mid m^{\bar{r}}_j\mid \leq 2\bar{N}$ and conditions \eqref{conjugué} to \eqref{décr} of Lemma 
\ref{m_j} are also satisfied.

\bigskip
\noindent Lemma \ref{m_j} implies that if conditions \eqref{conjugué} and \eqref{res2b} are not required, 
then one can get $m_1,\dots m_n\in\mathbb{Z}^d$. 
$\Box$

\subsubsection{Reduction of the eigenvalue}\label{renorm} 

\noindent Now the preceding lemmas will be used to define the map of reduction of the eigenvalues $\Phi$ which will conjugate 
$A$ to a matrix with $DC^{RN}_\omega(\kappa'',\tau)$ spectrum for some $\kappa''$, with $R,N$ 
arbitrarily great and $\Phi$ bounded independently of $R$. 

\bigskip
\noindent In all that follows, $G$ will be a Lie group among 

$$GL(n,\mathbb{C}), GL(n,\mathbb{R}),Sp(n,\mathbb{R}) ,SL(2,\mathbb{C}), SL(n,\mathbb{R}), O(n),U(n)$$ 

\noindent and $\mathcal{G}$ will be the Lie algebra associated to $G$.

\begin{prop}\label{diophan}Let $A\in \mathcal{G}$, $R\geq 1$ and $N\in\mathbb{N}$. 
There exists $\bar{N}\in [N,R^{\frac{1}{2}n(n-1)}N]$ such that if

\begin{equation}\kappa''= \frac{\kappa}{n(8R^{\frac{1}{2}n(n-1)+1}{N})^{\tau}}\end{equation}

\noindent then there exists a map $\Phi$ which is trivial with respect to $\mathcal{L}_{A,\kappa''}$ and $G$-valued and such that

\begin{enumerate}

\item for all $r'\geq 0$,

\begin{equation}\label{normephi}|\Phi|_{r'}\leq nC_0\left(\frac{1+||A_{\mathcal{N}}||}{\kappa''}\right)^{n(n+1)}e^{4\pi \bar{N}r'},\ |\Phi^{-1}|_{r'}\leq 
nC_0\left(\frac{1+||A_{\mathcal{N}}||}{\kappa''}\right)^{n(n+1)}e^{4\pi \bar{N}r'}
\end{equation}

\item Let $\tilde{A}$ be such that 

\begin{equation}\forall \theta\in 2\mathbb{T}^d,\ \partial_\omega \Phi(\theta)=A\Phi(\theta)-\Phi(\theta)\tilde{A}
\end{equation}

\noindent then

\begin{equation}\label{normeA} ||\tilde{A}-A||\leq 4\pi \bar{N}\end{equation}

\noindent and $\tilde{A}$ has $DC^{R\bar{N}}_\omega (\kappa'',\tau)$ spectrum.

\item If $\mathcal{G}=gl(n,\mathbb{C})$ or $u(n)$, $\Phi$ is defined on $\mathbb{T}^d$.

\item If $\mathcal{G}=o(n)$ or $u(n)$, then 

\begin{equation}\label{normephiorth}|\Phi|_{r'}\leq ne^{4\pi \bar{N}r'},\ |\Phi^{-1}|_{r'}\leq 
ne^{4\pi \bar{N}r'}
\end{equation}

\item If $\mathcal{G}=sl(2,\mathbb{C})$ or $sl(2,\mathbb{R})$, then either $\Phi$ is the identity or $\mid \mid \tilde{A}\mid \mid\leq \kappa''$.
\end{enumerate}

\end{prop}

\dem Let $\{\alpha_1,\dots ,\alpha_n\}=\sigma(A)$. Two cases must be considered: 
\begin{itemize}
\item If $\mathcal{G}=gl(n,\mathbb{C})$ or $u(n)$, Lemma \ref{res} gives 
$\bar{N}$ and $m_j\in \mathbb{Z}^d$ for $j=1,\dots ,n$ such that

\begin{equation}\nonumber N\leq \bar{N}\leq R^{\frac{1}{2}n(n-1)}N;\ \sup_j |m_j|\leq 2\bar{N}\end{equation}

\noindent and such that conditions \eqref{oppconjugué} to \eqref{décr} of Lemma \ref{m_j} hold with $\kappa'=\kappa''$, as well as conditions \eqref{res2a}.

\item If $\mathcal{G}=gl(n,\mathbb{R}), sp(n,\mathbb{R}),sl(n,\mathbb{R}),sl(2,\mathbb{C})$ or $o(n)$, Lemma \ref{res} gives 
$\bar{N}$ and $m_j\in \frac{1}{2}\mathbb{Z}^d$ for $j=1,\dots ,n$ such that

\begin{equation}\nonumber N\leq \bar{N}\leq R^{\frac{1}{2}n(n-1)}N;\ \sup_j |m_j|\leq 2\bar{N}\end{equation}

\noindent and such that conditions \eqref{conjugué} to \eqref{décr} of Lemma \ref{m_j} hold with $\kappa'=\kappa''$, as well as conditions \eqref{res2a} and \eqref{res2b}.

\end{itemize}

\bigskip
\noindent For all $j$ 
there is a unique $L\in \mathcal{L}_{A,\kappa''}$ such that 
$\alpha_j\in \sigma(A_{|L})$. Let $m_L=m_j$. Then $m_L$ is independent of $j$ thanks to property \eqref{prox}. 

\noindent For all $\theta\in 2\mathbb{T}^d$, let

\begin{equation}\nonumber\Phi(\theta)=\sum_{L\in\mathcal{L}_{A,\kappa''}}e^{2i\pi\langle m_L,\theta\rangle}P^{\mathcal{L}_{A,\kappa''}}_L\end{equation}

\noindent By construction of the $(m_L)$, $\Phi$ is defined on $\mathbb{T}^d$ if $\mathcal{G}=gl(n,\mathbb{C})$ or $u(n)$. 
Let us prove that $\Phi$ is $G$-valued. 
\begin{itemize}
\item if $\mathcal{G}=gl(n,\mathbb{C})$, this is trivial;
\item if $\mathcal{G}=sl(2,\mathbb{C})$ or $sl(2,\mathbb{R})$, this comes from condition \eqref{sl2C};
\item if $\mathcal{G}=u(n)$, $\Phi$ has unitary values. 
\item if $\mathcal{G}=gl(n,\mathbb{R})$, this comes from Lemma \ref{réel}, since $\mathcal{L}_{A,\kappa''}$ is a real decomposition and according to Lemma 
\ref{res}, for all $L\in \mathcal{L}_{A,\kappa''}$, $m_L=-m_{\bar{L}}$. 
\item if $\mathcal{G}=o(n)$, the map $\Phi$ has values in real unitary matrices, i.e orthogonal matrices. 
\item if $\mathcal{G}= sp(n,\mathbb{R})$, $\mathcal{L}_{A,\kappa''}$ is a symplectic decomposition. 
Lemma \ref{res} ensures that for all $L\in \mathcal{L}_{A,\kappa''}$, $m_L=-m_{\bar{L}}$ and

\begin{equation}\nonumber\forall L,L'\in \mathcal{L}_{A,\kappa''},\ 
\langle L,JL' \rangle \neq 0\Rightarrow m_L=m_{L'}
\end{equation}

\noindent 
Therefore Lemma \ref{phisp} 
implies that for all $\theta$ the matrix $\Phi(\theta)$ is in $Sp(n,\mathbb{R})$. 
\end{itemize}

\bigskip
\noindent Properties \eqref{res2a} and \eqref{res2b} ensure that 
$\tilde{A}$ has $DC^{R\bar{N}}_\omega (\kappa'',\tau)$ spectrum.

\noindent Moreover, for all $L\in\bar{\mathcal{L}}$, $|m_L|\leq 2\bar{N}$. The estimate of each 
$P^{\bar{\mathcal{L}}}_L$ recalled in Lemma \ref{c0} implies that 
$\Phi$ satisfies the estimate

\begin{equation}\nonumber|\Phi|_{r'}\leq nC_0\left(\frac{1+||A_{\mathcal{N}}||}{\kappa''}\right)^{n(n+1)}
e^{4\pi \bar{N}r'}
\end{equation}

\noindent and so does $\Phi^{-1}$ since

\begin{equation}\nonumber\Phi^{-1}=\sum_{L\in\mathcal{L}_{A,\kappa''}}e^{-2i\pi\langle m_L,.\rangle}P^{\mathcal{L}_{A,\kappa''}}_L\end{equation}

\bigskip
\noindent Now if $\mathcal{G}$ is $o(n)$ or $u(n)$, then every projection $P^{\mathcal{L}_{A,\kappa''}}_L$ has norm 1 and therefore $\Phi$ and $\Phi^{-1}$ satisfy \eqref{normephiorth}. 
By definition of $\tilde{A}$, 

\begin{equation}\nonumber\forall L\in\mathcal{L}',\ \sigma(\tilde{A}_{|L})=\sigma(A_{|L})-2i\pi\langle
m_L,\omega\rangle\end{equation}

\noindent and by property \eqref{décr}, 

\begin{equation}\nonumber\forall \alpha\in \sigma(A_{|L}),\ |\alpha -2i\pi\langle
m_L,\omega\rangle| \leq |\alpha|\end{equation} 

\noindent Let $P$ be such that $PAP^{-1}$ is in Jordan normal form, let $\alpha_j$ be the eigenvalues of $A$ and $p_j$ the columns of $P$, then for all $j$,

\begin{equation}\nonumber ||(\tilde{A}-A)p_j||=||2i\pi\langle m_j,\omega\rangle
p_j||\leq 4\pi \bar{N}  ||p_j||
\end{equation}

\noindent So $||\tilde{A}-A||\leq 4\pi \bar{N}$, whence property \eqref{normeA}. Finally, if $\mathcal{G}=sl(2,\mathbb{C})$ or $sl(2,\mathbb{R})$, then either $\tilde{A}=A$, or $A$ is diagonalizable, and then $\tilde{A}$ is also diagonalizable, so their norms are the modulus of their eigenvalues and by condition \eqref{décr'}, $\mid \mid \tilde{A}\mid \mid \leq 
\kappa''$. $\Box$

\bigskip
\noindent \Def A map $\Phi$ satisfying the conclusion of Proposition \ref{diophan} will be called a \textit{map of reduction of the eigenvalues of $A$ 
at order $R,\bar{N}$}.

\bigskip
\noindent In dimension 2, 
the map of reduction of the eigenvalues $\Phi$ 
satisfies the following property: for every function $H$ continuous on $\mathbb{T}^d$ and with values in 
$gl(2,\mathbb{C})$, $\Phi H\Phi^{-1}$ and $\Phi^{-1} H\Phi$ are continuous on $\mathbb{T}^d$.

\noindent Dimension 2 has, indeed, the particularity that every decomposition $\mathcal{L}$ of $\mathbb{R}^2$ at most two subpaces $L_1,L_2$, 
in which case $m_{L_1}+m_{L_2}
\in\mathbb{Z}^d$ (if the decomposition is trivial, $m_L=0$). In any case, $\sum_{L\in \mathcal{L}}
m_L\in \mathbb{Z}^d$.

\subsection{Homological equation}

\noindent Solving the homological equation is a first step towards reducing the perturbation. 

\bigskip
\noindent 
\textbf{Notation:} For every function $F\in L^2(2\mathbb{T}^d)$ and every $N\in\mathbb{N}$, we will denote by $F^N$ and call 
truncation of $F$ at order $N$ the function that one obtains by truncating the Fourier series of $F$:

$$F^N(\theta)=\sum_{|m|\leq N} \hat{F}(m)e^{2i\pi\langle m,\theta\rangle}$$

\bigskip
\noindent The following lemma will be useful in the solving of the homological equation.

\begin{lem}\label{polyn} Let $f,g$ be trigonometric polynomial with $g$ real on $\mathbb{R}^d$. Let $r>0, r'\in ]0,r[$ and suppose that there exists $C$ such that $|f|_{r'}\leq C|g|_r$. Then for all $m\in \frac{1}{2}\mathbb{Z}^d$, 

\begin{equation}|fe^{2i\pi\langle m,.\rangle}|_{r'}\leq C|ge^{2i\pi\langle m,.\rangle}|_r \end{equation}

\end{lem}

\dem Since $g$ is real, 

\begin{equation}\label{ghatg}\forall m\in \mathbb{Z}^d, \ \overline{\hat{g}(-m)}=\hat{g}(m)
\end{equation}

\noindent 
so for all $x$ and all $y\in [-r,r]^d$,

\begin{equation}\nonumber g(x-iy)=\sum_m \hat{g}(m)e^{2i\pi \langle m,x-iy\rangle}
=\sum_m \overline{\hat{g}(-m)}e^{2i\pi \langle -m,-x+iy\rangle}=\overline{\sum_m \hat{g}(-m)e^{2i\pi \langle -m,x+iy\rangle}}
=\overline{g(x+iy)}
\end{equation}

\noindent which implies that for all $x,y$,

\begin{equation}\label{g}\mid g(x-iy)\mid=\mid g(x+iy)\mid 
\end{equation}

\noindent Let us show that for every $m\in \mathbb{Z}^d$,

\begin{equation}\label{mzd}|g|_re^{2\pi |m|r}
=|ge^{2i\pi\langle m,.\rangle}|_{r}\end{equation}

\noindent  By the maximum principle,

\begin{equation}\nonumber |g|_{r}=\sup_{x;|y_j|\leq r,1\leq j\leq d}|g(x+iy)|=\sup_{x;|y_j|=r, 1\leq j\leq d}|g(x+iy)|
\end{equation}

\noindent Let $y_0$ such that 

\begin{equation}\nonumber \mid g\mid _r=\sup_x\mid g(x+iy_0)\mid
\end{equation}

\noindent then, for $m$ having only one non-zero component $m_j$, either

\begin{equation}\nonumber 
|g|_re^{2\pi |m|r}=\sup_x \mid g(x+iy_0)\mid \ \mid e^{2i\pi \langle m, x+iy_0\rangle}\mid
=|ge^{2i\pi\langle m,.\rangle}|_{r}\end{equation}
 
\noindent if $m_j$ et $(y_0)_j$ have opposite signs, or

\begin{equation}\nonumber 
|g|_re^{2\pi |m|r}=\sup_x \mid g(x-iy_0)\mid \ \mid e^{2i\pi \langle m, x-iy_0\rangle}
=|ge^{2i\pi\langle m,.\rangle}|_{r}\end{equation}

\noindent if $m_j$ et $(y_0)_j$ have the same sign, 
whence \eqref{mzd} if $m$ has only one non-zero component. For $1\leq l\leq d$, let $\bar{m}_l=(m_1,\dots, m_l,0,\dots, 0)$. 
Assume that

\begin{equation}\nonumber \mid g\mid_re^{2\pi \mid m\mid r}=\mid g e^{2i\pi \langle  \bar{m}_{j-1},.\rangle}\mid _re^{2\pi (\mid m_j\mid+\dots +\mid m_d\mid) r}
\end{equation}

\noindent and that $\mid g e^{2i\pi \langle \bar{m}_{j-1},.\rangle}\mid _r$ is reached at $\bar{y}$. 
Let $\delta_j\in \{-1,1\}$ be such that $m_j$ and $\delta_j\bar{y}_j$ have opposite signs. 
Then

\begin{equation}\nonumber \begin{split}\mid g\mid_re^{2\pi \mid m\mid r}
&=\mid g e^{2i\pi \langle  \bar{m}_{j-1},.\rangle}\mid _re^{2\pi (\mid m_j\mid+\dots +\mid m_d\mid) r}\\
&=\sup_{x,y_k,k\neq j} \mid g(x+i(y_1,\dots,\bar{y}_j,\dots , y_d)) e^{2i\pi \langle  \bar{m}_{j-1},x+i(y_1,\dots,\bar{y}_j,\dots , y_d) \rangle}\mid e^{2\pi (\mid m_j\mid+\dots +\mid m_d\mid) r}\\
&=\sup_{x,y_k,k\neq j} \mid g(x+i(y_1,\dots,\delta_j\bar{y}_j,\dots , y_d)) e^{2i\pi \langle  \bar{m}_{j-1},x+i(y_1,\dots,\delta_j\bar{y}_j,\dots , y_d)\rangle}e^{2i\pi m_j(x_j+i\delta_j\bar{y}_j)}\mid \\
&.e^{2\pi (\mid m_{j+1}\mid+\dots +\mid m_d\mid) r}\\
&=\sup_{x,y_k,k\neq j} \mid g(x+i(y_1,\dots,\delta_j\bar{y}_j,\dots , y_d)) e^{2i\pi \langle \bar{m}_{j},x+i(y_1,\dots,\delta_j\bar{y}_j,\dots , y_d)\rangle}\mid
e^{2\pi (\mid m_{j+1}\mid+\dots +\mid m_d\mid) r}\\
&=\mid g e^{2i\pi \langle \bar{m}_{j},.\rangle}\mid _re^{2\pi (\mid m_{j+1}\mid+\dots +\mid m_d\mid) r}
\end{split}\end{equation}

\noindent and \eqref{mzd} is obtained through a simple iteration. Thus

\begin{equation}\nonumber |fe^{2i\pi\langle m,.\rangle}|_{r'}\leq |f|_{r'}e^{2\pi |m|r'}
\leq C|g|_re^{2\pi |m|r}
= C|ge^{2i\pi\langle m,.\rangle}|_r\ \Box\end{equation}

\bigskip
\noindent \rem If $f,g$ are matrix-valued trigonometric polynomials, $f=(f_{j,k}),g=(g_{j,k})$, 
and $g$ has real coefficients on $\mathbb{R}^d$, a similar statement holds. For if

\begin{equation}\nonumber |f|_{r'}=\sup_{x,\mid y_j\mid \leq r'} \mid \mid f(x+iy)\mid \mid \leq C|g|_r=C\sup_{x,\mid y_j\mid \leq r} \mid \mid g(x+iy)\mid \mid\end{equation}

\noindent as the norm of the greatest coefficient is equivalent to the operator norm, one has


\begin{equation}\nonumber \sup_{j,k}|f_{j,k}|_{r'}
\leq CC'\sup_{j,k}|g_{j,k}|_r \end{equation}

\noindent  for some $C'$ only depending on the dimension of the matrices.
So from Lemma \ref{polyn}, since there exists $j_0,k_0$ such that 

\begin{equation}\nonumber \forall j,k, \ \mid f_{j,k}\mid_{r'} \leq CC' \mid g_{j_0,k_0}\mid _r
\end{equation}

\noindent then

\begin{equation}\nonumber \sup_{j,k}|f_{j,k}e^{2i\pi\langle m,.\rangle}|_{r'}\leq CC'\sup_{j,k}|g_{j,k}e^{2i\pi\langle m,.\rangle}|_r \end{equation}

\noindent and as the norms are equivalent, the statement also holds in operator norm:

\begin{equation}\nonumber |fe^{2i\pi\langle m,.\rangle}|_{r'}\leq CC''|ge^{2i\pi\langle m,.\rangle}|_r \end{equation}

\noindent for some $C''$ depending only on the dimension of the matrices.


\begin{prop}\label{homol} 
Let
\begin{itemize}
\item ${N}\in\mathbb{N}$,
\item $\kappa'\in ]0, \kappa]$, 
\item $\gamma \geq n(n+1)$,
\item
 $0<r'<r$. 
 \end{itemize}
Let $\tilde{A}\in \mathcal{G}$ have $DC^{{N}}_\omega(\kappa',\tau)$ spectrum. 
Let $\tilde{F}\in C^\omega_r(2\mathbb{T}^d,  \mathcal{G})$ with nice periodicity properties with respect to an $(\tilde{A},\kappa',\gamma)$-decomposition $\mathcal{L}$. 
Then equation 

\begin{equation}\label{hom}\forall \theta\in 2\mathbb{T}^d, \ \partial_\omega \tilde{X}(\theta)=[\tilde{A},\tilde{X}(\theta)]+\tilde{F}^{{N}}(\theta)
-\hat{\tilde{F}}(0);\ \hat{\tilde{X}}(0)=0
\end{equation}

\noindent has a solution $\tilde{X}\in C^\omega_{r'}(2\mathbb{T}^d,  \mathcal{G})$
such that 
\begin{itemize}
\item if $\tilde{F}$ has nice periodicity properties with respect to
$\mathcal{L}$ and $(m_L)$, 
then $\tilde{X}$ has nice periodicity properties with respect to $\mathcal{L}$ and $(m_L)$; in particular, if $\tilde{F}$ is defined on $\mathbb{T}^d$, then so is $\tilde{X}$,

\item  
if $\Phi$ is trivial with respect to $\mathcal{L}$, 
then there exist $C',D$ depending only on 
$n,d,\tau$ such that





\begin{equation}\label{X2}|\Phi^{-1}\tilde{X}\Phi|_{r'}
\leq C'\left(\frac{1+||\tilde{A}_{\mathcal{N}}||}{(r-r')\kappa'}\right)^{2n^2\gamma+D}
|\Phi^{-1}\tilde{F}\Phi
|_r\end{equation}

\end{itemize}

\bigskip
Moreover, the truncation of 
$\tilde{X}$ at order $N$ is unique. 

\end{prop}

\dem $\bullet$ Let $C\in GL(n,\mathbb{C})$ be such that $C^{-1}\tilde{A}C$ is in Jordan normal form. 
Conjugating equation (\ref{hom}) by $C^{-1}$, decomposing into coefficients $x_{j,k}$ of $C^{-1}\tilde{X}C$ and 
developing into Fourier series, one gets for all $m\in\frac{1}{\nu}\mathbb{Z}^d$, with $\nu=1$ or 2 according to the periodicity of 
$(C^{-1}({\tilde{F}}^{{N}}-\hat{\tilde{F}}(0))C)_{j,k}$,

\begin{equation}\label{fou}i\langle m,\omega\rangle \hat{x}_{j,k}(m)
=(\tilde{\alpha}_j-\tilde{\alpha}_k)\hat{x}_{j,k}(m)
+\delta_1\hat{x}_{j,k+1}(m)+\delta_2\hat{x}_{j-1,k}(m)
+\hat{f}(m)\end{equation}

\noindent where $\delta_1,\delta_2$ are $0$ or $1$ and $\hat{f}(m)$ stands for the $m$-th Fourier coefficient of the function 
$(C^{-1}({\tilde{F}}^{{N}}-\hat{\tilde{F}}(0))C)_{j,k}$.

\noindent The diophantine conditions given by Proposition \ref{diophan} allow the existence of an analytic 
solution to the set of equations $(\ref{fou})$, 
therefore (\ref{hom}) has a solution $\tilde{X}$.

\bigskip
\noindent $\bullet$ Now we shall see that $\tilde{X}^{{N}}$ is unique. Suppose that $\tilde{X}$ and $\tilde{Y}$ 
are both solutions of (\ref{hom}). Then

\begin{equation}\label{uni}\partial_\omega(\tilde{X}-\tilde{Y})=[\tilde{A},\tilde{X}-\tilde{Y}];\ 
\hat{\tilde{X}}(0)-\hat{\tilde{Y}}(0)=0\end{equation}

\noindent The diophantine conditions on $\tilde{A}$ imply that the truncation at order ${N}$ of any solution of 
(\ref{uni}) is constant, 
and condition 
$\hat{\tilde{X}}(0)-\hat{\tilde{Y}}(0)=0$ implies that it vanishes, so $\tilde{X}^{{N}}=\tilde{Y}^{{N}}$.

\bigskip
\noindent $\bullet$ To check that $\tilde{X}$ is $\mathcal{G}$-valued, it is enough to show it for 
$\tilde{X}^{{N}}$, since one can assume that $\tilde{X}=\tilde{X}^{{N}}$.

\begin{itemize}

\item if $\mathcal{G}=gl(n,\mathbb{C})$, this is trivial.

\item if $\mathcal{G}=gl(n,\mathbb{R})$, this comes from the unicity of the solution up to order ${N}$, 
since $\tilde{X}$ and its complex conjugate are solutions of the same equation.

\item if $\mathcal{G}=sp(n,\mathbb{C})$, then



%


\begin{equation}\nonumber \begin{split}\forall \theta\in 2\mathbb{T}^d, \ \partial_\omega J(\tilde{X}(\theta)^*J+J\tilde{X}(\theta))&=-J(\tilde{X}(\theta)^*J
+J\tilde{X}(\theta))\tilde{A}-
J\tilde{A}^*(\tilde{X}(\theta)^*J+J\tilde{X}(\theta))\\
&=[\tilde{A},J(\tilde{X}(\theta)^*J+J\tilde{X}(\theta))]
\end{split}\end{equation}

\noindent Diophantine conditions on $\tilde{A}$ imply that $\tilde{X}^*J+J\tilde{X}$ is constant. 
Condition $\hat{\tilde{X}}(0)=0$ implies
that for every $\theta\in 2\mathbb{T}^d$, $\tilde{X}(\theta)^*J+J\tilde{X}(\theta)=0$, so $\tilde{X}$ takes its values in $sp(n,\mathbb{C})$.

\item if $\mathcal{G}=u(n)$, proceed as in the $sp(n,\mathbb{C})$ case, showing this time that 
$\tilde{X}^*+\tilde{X}$ is constant and thus is zero.

\item if $\mathcal{G}=sp(n,\mathbb{R})$ or $o(n)$, use the previous cases and the fact that 
$sp(n,\mathbb{R})=sp(n,\mathbb{C})\cap gl(n,\mathbb{R})$ and $o(n)=u(n)\cap gl(n,\mathbb{R})$. 

\item if $\mathcal{G}=sl(n,\mathbb{R})$ or $sl(2,\mathbb{C})$, note that the trace of $\tilde{X}$ is solution of

\begin{equation}\nonumber \forall \theta \in 2\mathbb{T}^d,\ \partial_\omega ( Tr \tilde{X}(\theta))=
Tr[\tilde{A},\tilde{X}(\theta)]=Tr (\tilde{A}\tilde{X}(\theta))-Tr (\tilde{X}(\theta)\tilde{A})=0
\end{equation}

\noindent so it is a constant, and as $Tr \hat{\tilde{X}}(0)=0$, it is identical to zero.

\end{itemize}

\bigskip
\noindent $\bullet$ As for periodicity properties, equation 
\eqref{hom} decomposes into blocks according to $\mathcal{L}$, then into Fourier coefficients: for $0<|m|\leq {N}$, 
%


\begin{equation}\label{fourier}2i\pi\langle m,\omega\rangle (P^{\mathcal{L}}_L \hat{\tilde{X}}(m)
P^{\mathcal{L}}_{L'})
=P^{\mathcal{L}}_L\tilde{A}P^{\mathcal{L}}_L\hat{\tilde{X}}(m)P^{\mathcal{L}}_{L'}
-P^{\mathcal{L}}_L\hat{\tilde{X}}(m)P^{\mathcal{L}}_{L'}\tilde{A}P^{\mathcal{L}}_{L'}
+P^{\mathcal{L}}_L\hat{\tilde{F}}(m)P^{\mathcal{L}}_{L'}
\end{equation}

\noindent Let $(m_L)$ be a family such that $\tilde{F}$ has nice periodicity properties with respect to $\mathcal{L}$ and $(m_L)$. If 
$m$ is not in $\mathbb{Z}^d+m_L-m_{L'}$, then $P^{\mathcal{L}}_L\hat{\tilde{F}}(m)
P^{\mathcal{L}}_{L'}
=0$ 
and since $\tilde{X}^{{N}}$ is unique, $P^{\mathcal{L}}_L\hat{\tilde{X}}(m)
P^{\mathcal{L}}_{L'}=0$. 
For $|m|>{N}$ one can assume $\hat{\tilde{X}}(m)=0$.
Therefore $\tilde{X}$ also has nice periodicity properties with respect to 
$\mathcal{L}$ et $(m_L)$. 

\bigskip
\noindent $\bullet$ Finally let us prove the estimate \eqref{X2}. 
Let $m\in\frac{1}{2}\mathbb{Z}^d, |m|\leq {N}$. First we shall prove that for all $L,L'\in \mathcal{L}$,

\begin{equation}\label{estimationre}||P^{\mathcal{L}}_L\hat{\tilde{X}}(m)
P^{\mathcal{L}}_{L'}||
\leq C'\frac{(1+||\tilde{A}_{\mathcal{N}}||)^{n^2-1}|m|^{(n^2-1)\tau}}{\kappa'^{(n^2-1)}}
||P^{\mathcal{L}}_L\hat{\tilde{F}}(m)
P^{\mathcal{L}}_{L'}||  (||P^{\mathcal{L}}_L|| \ 
||P^{\mathcal{L}}_{L'}||)^{n^2-1}
\end{equation}

\noindent where $C'$ only depends on $n$. The proof will be inspired by \cite{E1}, Lemma 2. Let $\mathcal{A}_{L,L'}$ be the linear operator from $gl(n,\mathbb{C})$ into itself such that for all $M\in gl(n,\mathbb{C})$,

\begin{equation}\nonumber \mathcal{A}_{L,L'}M=
\tilde{A} P^{\mathcal{L}}_L
M-M P^{\mathcal{L}}_{L'}
 \tilde{A} 
\end{equation}

\noindent Decomposing \eqref{hom} into blocks, then into Fourier series, one obtains for all $L,L'\in \mathcal{L}$ and all $m\in \frac{1}{2}\mathbb{Z}^d$ such that $0<\mid m\mid \leq N$,

%

%



%

%

\begin{equation}\label{invR} (P^{\mathcal{L}}_L \hat{\tilde{X}}(m)
P^{\mathcal{L}}_{L'})
=(2i\pi\langle m,\omega\rangle-\mathcal{A}_{L,L'})^{-1}P^{\mathcal{L}}_L\hat{\tilde{F}}(m)P^{\mathcal{L}}_{L'}
\end{equation}

\noindent 
Write $\mathcal{A}_{L,L'}$ as an $n^2$-dimensional matrix. 
Let $A_D\in gl(n^2,\mathbb{C})$ be a diagonal matrix and $A_N\in gl(n^2,\mathbb{C})$ a nilpotent matrix such that

\begin{equation}\nonumber (2i\pi\langle m,\omega\rangle-\mathcal{A}_{L,L'})=A_D-A_N
\end{equation}

\noindent 
Then $A_N$ coincides with the operator

\begin{equation}\nonumber A_N:B\mapsto (\tilde{A} P^{\mathcal{L}}_L)_{\mathcal{N}}
B-B (P^{\mathcal{L}}_{L'}
\tilde{A})_{\mathcal{N}}\end{equation}

\noindent Moreover,

\begin{equation}\nonumber (2i\pi\langle m,\omega\rangle-\mathcal{A}_{L,L'})^{-1}=A_D^{-1}(I+A_NA_D^{-1}+\dots + (A_NA_D^{-1})^{n^2-1})
\end{equation}

\noindent We will estimate $(2i\pi\langle m,\omega\rangle-\mathcal{A}_{L,L'})^{-1}$, for $m\in \mathbb{Z}^d$ if $L=\bar{L}'$ and 
$m\in \frac{1}{2}\mathbb{Z}^d$ if $L\neq \bar{L}'$.
Each coefficient of $A_D^{-1}(A_NA_D^{-1})^{j-1}$ has the form $\frac{p}{q}$ with $\mid p\mid \leq \mid \mid A_N\mid \mid ^{j-1}$ and $q=\beta_1\dots \beta_j$ where 
 $\beta_i$ are eigenvalues of $ 2i\pi\langle m,\omega\rangle-\mathcal{A}_{L,L'}$. Now

\begin{equation}\nonumber \sigma(\mathcal{A}_{L,L'})=\{\alpha-\alpha'\ \mid \ \alpha\in \sigma(\tilde{A}_{\mid L}),\alpha'\in \sigma(\tilde{A}_{\mid L'}) \}
\end{equation}

\noindent 
and for all $\alpha\in \sigma(\tilde{A}_{\mid L}),\alpha'\in \sigma(\tilde{A}_{\mid L'})$,

\begin{equation}\nonumber \mid \alpha-\alpha'- 2i\pi\langle m,\omega\rangle\mid \geq \frac{\kappa'}{\mid m\mid ^\tau}
\end{equation}

for all $m\in \mathbb{Z}^d$ if $L=\bar{L}'$ and all $m\in \frac{1}{2}\mathbb{Z}^d$ if $L\neq \bar{L}'$
Thus

\begin{equation}\nonumber \begin{split} \mid \mid (2i\pi\langle m,\omega\rangle-\mathcal{A}_{L,L'})^{-1}\mid \mid
&\leq n^22^{n^2} (1+\mid \mid      \tilde{A} _{\mathcal{N}}   \mid \mid \ (\mid \mid  P^{\mathcal{L}}_L\mid \mid \ 
+\mid \mid P^{\mathcal{L}}_{L'}  \mid \mid ))^{n^2-1}\left(\frac{\mid m\mid^\tau}{\kappa'}\right)^{n^2-1}
\end{split}\end{equation}

and \eqref{invR} implies \eqref{estimationre}.

%

%


%

%

\bigskip
\noindent $\bullet$ The estimate  \eqref{estimationre}
implies that

\begin{equation}\label{xr}\begin{split}|P^{\mathcal{L}}_L\tilde{X}P^{\mathcal{L}}_{L'}|_{r'}
&\leq C'\frac{(1+||\tilde{A}_{\mathcal{N}}||)^{n^2-1}}{\kappa'^{(n^2-1)}}\sum_m|m|^{(n^2-1)\tau}
|P^{\mathcal{L}}_L\tilde{F}
P^{\mathcal{L}}_{L'}|_re^{-2\pi\mid m\mid r} e^{2\pi \mid m\mid r'}
 (||P^{\mathcal{L}}_L|| \ 
||P^{\mathcal{L}}_{L'}||)^{n^2-1}\\
\end{split}\end{equation}

\noindent where $C'$ only depends on $n$. Now

\begin{equation}\nonumber \begin{split}\sum_m|m|^{(n^2-1)\tau}e^{-2\pi \mid m\mid(r- r')}
&\leq C_d \sum_{M\geq 1}M^{(n^2-1)\tau+d}e^{-2\pi M(r- r')}\\
&\leq C_d \int_0^\infty t^{(n^2-1)\tau+d}e^{-2\pi t(r- r')}dt
\leq \frac{C_d}{(2\pi (r-r'))^{(n^2-1)\tau+d+1}}
\end{split}\end{equation}

\noindent where $C_d$ only depends on $d$, so

\begin{equation}\label{xr'}\begin{split}|P^{\mathcal{L}}_L\tilde{X}P^{\mathcal{L}}_{L'}|_{r'}
&\leq \frac{C''}{(r-r')^{(n^2-1)\tau+d+1}}\frac{(1+||\tilde{A}_{\mathcal{N}}||)^{n^2-1}}{\kappa'^{(n^2-1)}}|P^{\mathcal{L}}_L\tilde{F}
P^{\mathcal{L}}_{L'}|_r
 (||P^{\mathcal{L}}_L|| \ 
||P^{\mathcal{L}}_{L'}||)^{n^2-1}\end{split}\end{equation}

\noindent where $C''$ only depends on $n,d,\tau$.

%

\noindent Let $(m'_L)_{L\in \mathcal{L}}$ a family of elements of $\frac{1}{2}\mathbb{Z}^d$ and $\Phi$ defined by

\begin{equation}\nonumber \Phi=\sum_{L\in\mathcal{L}}P^{\mathcal{L}}_Le^{2i\pi \langle m'_L,.\rangle}\end{equation}

\noindent then

\begin{equation}\nonumber \begin{split}|\Phi^{-1}\tilde{X}\Phi|_{r'}
& = |\sum_{L,L'\in \mathcal{L}}P^{\mathcal{L}}_L\tilde{X}e^{2i\pi \langle m'_L-m'_{L'},.\rangle}P^{\mathcal{L}}_{L'}|_{r'}\\
\end{split}\end{equation}

\noindent and since $\mathcal{L}$ is an $(\tilde{A},\kappa',\gamma)$-decomposition, then Lemma \ref{polyn} applied to \eqref{xr'} gives

\begin{equation}\nonumber \begin{split}
|\Phi^{-1}\tilde{X}\Phi|_{r'}
&\leq \frac{C_3}{(r-r')^{(n^2-1)\tau+d+1}}\left(\frac{1+||\tilde{A}_{\mathcal{N}}||}{\kappa'}\right)^{n^2(2\gamma+1)}
\sum_{L,L'}|P^{\mathcal{L}}_L
\Phi^{-1}\tilde{F}\Phi
P^{\mathcal{L}}_{L'}|_r\ \end{split}\end{equation}

\bigskip
\noindent where $C_3$ only depends on $n,d,\tau$, whence \eqref{X2}. $\Box$
%

%

\bigskip
\rem The loss of analyticity $r-r'$ is needed in order to have good estimates of the solution. Note that when $\mathcal{G}=o(n)$ or $u(n)$, then $\tilde{A}_\mathcal{N}$ is zero, thus 
the estimate does not depend on $\tilde{A}$.

\subsection{Inductive lemma without reduction of the eigenvalues}

\subsubsection{Auxiliary lemmas}

\noindent The first lemma will be used to iterate the inductive lemma without having to 
perform reduction of the eigenvalues at each step, which will greatly improve the final estimates.

\begin{lem}\label{NRdurable} Let
\begin{itemize}
\item $\kappa'\in ]0,1[$, $C>0$,
\item $\tilde{F}\in \mathcal{G}$,  
\item $\tilde{\epsilon}=||\tilde{F}||$,
\item $\tilde{N}\in \mathbb{N}$,
\item $\tilde{A}\in \mathcal{G}$ with $DC^{\tilde{N}}_\omega(\kappa',\tau)$ spectrum.

\end{itemize}

\noindent There exists a constant $c$ only depending on $n\tau $ 
such that if $\tilde{\epsilon}$ satisfies

\begin{equation}\label{cond1}\tilde{\epsilon}\leq c\left(\frac{C^\tau \kappa'}{1+||\tilde{A}||}\right) ^{2n}
\end{equation} 

\noindent and

\begin{equation}\label{cond2} \tilde{N}
\leq  \frac{|\log \tilde{\epsilon}|}{C}
\end{equation}

\noindent then
$\tilde{A}+{\tilde{F}}$ has $DC^{\tilde{N}}_\omega(\frac{3\kappa'}{4},\tau)$ spectrum.
\end{lem}

\dem If $\tilde{\alpha}\in\sigma(\tilde{A}+\tilde{F})$, 
by Lemma \ref{appendice} given as an appendix, there exists
$\alpha\in \sigma(\tilde{A})$ such that $|\alpha-\tilde{\alpha}|\leq 2n(||\tilde{A}||+1)\tilde{\epsilon}^{\frac{1}{n}}$.

\bigskip
\noindent By assumption $\tilde{A}$ has $DC^{\tilde{N}}_\omega(\kappa',\tau)$ spectrum. 
Thus for all $\alpha,\beta\in \sigma(
\tilde{A}+{\tilde{F}})$ and all $m\in\mathbb{Z}^d, 0<|m|\leq \tilde{N}$,

\begin{equation}\label{diopert}|\alpha-\beta-	2i\pi
\langle m,\omega\rangle|\geq \frac{\kappa'}{|m|^\tau}-4n(||\tilde{A}||+1)\tilde{\epsilon}^\frac{1}{n}\end{equation} 

\noindent and if $\alpha\neq \bar{\beta}$, \eqref{diopert} holds for every $m\in\frac{1}{2}\mathbb{Z}^d$, $0<|m|\leq \tilde{N}$. 
Therefore it is enough to show that

\begin{equation}\nonumber 4n\tilde{N}^\tau(||\tilde{A}||+1)\tilde{\epsilon}^\frac{1}{n}\leq \frac{\kappa'}{4}
\end{equation}

\noindent Now there is a constant $c\leq 1$ which only depends on $n\tau $ such that if $\tilde{\epsilon}\leq c$, then

\begin{equation}\nonumber \tilde{\epsilon} \left({|\log \tilde{\epsilon}|}\right)^{n\tau }\leq  \tilde{\epsilon}^\frac{1}{2}\end{equation} 

\noindent so if

\begin{equation}\nonumber \tilde{\epsilon}\leq c\left(\frac{C^\tau\kappa'}{16n(||\tilde{A}||+1)}\right)^{2n} \end{equation} 

%

\noindent by asumption \eqref{cond2}, then

\begin{equation}\nonumber 4n(||\tilde{A}||+1)\tilde{\epsilon}^\frac{1}{n} \tilde{N}^{\tau}
\leq 4n(||\tilde{A}||+1)\tilde{\epsilon}^{\frac{1}{2n}}C^{-\tau}\leq 
\frac{\kappa'}{4}\end{equation} 

\noindent which proves the Lemma. $\Box$

\bigskip
\noindent If $G$ is compact, then by lemma \ref{appendice2}, the same is true under a smallness condition which does not depend on $\tilde{A}$.

\begin{lem}\label{NRdurableorth} Let
\begin{itemize}
\item $\kappa'\in ]0,1[$, $C>0$,
\item $\tilde{F}\in \mathcal{G}$,  
\item $\tilde{\epsilon}=||\tilde{F}||$,
\item $\tilde{N}\in \mathbb{N}$,
\item $\tilde{A}\in \mathcal{G}$ with $DC^{\tilde{N}}_\omega(\kappa',\tau)$ spectrum.

\end{itemize}

\noindent There exists a constant $c$ only depending on $\tau $ 
such that if $\tilde{\epsilon}$ satisfies

\begin{equation}\label{cond1orth}\tilde{\epsilon}\leq c(C^\tau \kappa')^2
\end{equation} 

\noindent and

\begin{equation}\label{cond2orth} \tilde{N}
\leq  \frac{|\log \tilde{\epsilon}|}{C}
\end{equation}

\noindent then
$\tilde{A}+{\tilde{F}}$ has $DC^{\tilde{N}}_\omega(\frac{3\kappa'}{4},\tau)$ spectrum.
\end{lem}

\dem If $\tilde{\alpha}\in\sigma(\tilde{A}+\tilde{F})$, 
by Lemma \ref{appendice2}, there exists
$\alpha\in \sigma(\tilde{A})$ such that $|\alpha-\tilde{\alpha}|\leq \tilde{\epsilon}$.
Since $\tilde{A}$ has $DC^{\tilde{N}}_\omega(\kappa',\tau)$ spectrum, 
then for all $\alpha,\beta\in \sigma(
\tilde{A}+{\tilde{F}})$ and all $m\in\mathbb{Z}^d, 0<|m|\leq \tilde{N}$,

\begin{equation}\label{diopert'}|\alpha-\beta-	2i\pi
\langle m,\omega\rangle|\geq \frac{\kappa'}{|m|^\tau}-2\tilde{\epsilon}\end{equation} 

\noindent and if $\alpha\neq \bar{\beta}$, \eqref{diopert'} holds for every $m\in\frac{1}{2}\mathbb{Z}^d$, $0<|m|\leq \tilde{N}$. 
There is a constant $c\leq 1$ which only depends on $\tau $ such that if $\tilde{\epsilon}\leq c$, then

\begin{equation}\nonumber \tilde{\epsilon} \left({|\log \tilde{\epsilon}|}\right)^{\tau }\leq  \tilde{\epsilon}^\frac{1}{2}\end{equation} 

\noindent so it is enough that

\begin{equation}\nonumber \tilde{\epsilon}\leq c\left(\frac{C^\tau\kappa'}{8}\right)^{2} \ \Box\end{equation} 

%

\bigskip
\noindent The following lemma will be used to avoid doubling the period more than once.

\begin{lem}\label{bppG} Let $A,A'\in gl(n,\mathbb{R})$ and $H:2\mathbb{T}^d\rightarrow gl(n,\mathbb{R})$. 
Assume that $H$ has nice periodicity properties with respect to an $A$-decomposition $\mathcal{L}$ and assume

\begin{equation}\label{parité}\forall L,L'\in\mathcal{L}, P_L^{\mathcal{L}}(A'-A)P^{\mathcal{L}}_{L'}\neq 0
\Rightarrow P_L^{\mathcal{L}}HP^{\mathcal{L}}_{L'}\in C^0(\mathbb{T}^d,gl(n,\mathbb{R}))
\end{equation}

\noindent Then $H$ has nice periodicity properties with respect to an $A'$-decomposition which is less fine than $\mathcal{L}$.

\end{lem}

\dem 
Define a decomposition ${\mathcal{L}}'$ of $\mathbb{C}^n$  as follows: for all $L,L'\in \mathcal{L} $,

\begin{equation}\nonumber  (\exists L_0\in {\mathcal{L}'}\ \mid \ L\subset L_0, L'\subset L_0 )
\Leftrightarrow P_L^{\mathcal{L}}HP^{\mathcal{L}}_{L'}\in C^0(\mathbb{T}^d,gl(n,\mathbb{R}))
\end{equation}

Let $(m_L)$ be a family such that $H$ has nice periodicity properties with respect to $\mathcal{L}$ and $(m_L)$. 
For all $L'\in {\mathcal{L}'}$, let $L$ be a subspace of $\mathcal{L}$ contained in $L'$ and let $\bar{m}_{L'}=m_{L}$; the class of $\bar{m}_{L'}$ in the equivalence relation 

\begin{equation}\nonumber m\sim m' \Leftrightarrow m-m'\in \mathbb{Z}^d
\end{equation}

\noindent does not depend on a particular choice of $L$. 
Then for all $L'\in {\mathcal{L}'}$,

\begin{equation}\nonumber e^{2i\pi\langle \bar{m}_{L'},.\rangle}
P^{{\mathcal{L}'}}_{L'}=\sum_{L\in\mathcal{L},L\subset L'}e^{2i\pi\langle \bar{m}_{L'},.\rangle}
P^{\mathcal{L}}_{L}\end{equation}

\noindent so for all $L_1,L_2\in \mathcal{L}'$,

\begin{equation}\nonumber P_{L_1}^{{\mathcal{L}'}}H P_{L_2}^{{\mathcal{L}'}}e^{2i\pi\langle \bar{m}_{L_1}
-\bar{m}_{L_2},.\rangle}=\sum_{L_1'\subset L_1,L_2'\subset L_2}P_{L'_1}^{{\mathcal{L}}}H
P_{L'_2}^{{\mathcal{L}}}e^{2i\pi\langle {m}_{L'_1}
-{m}_{L'_2},.\rangle}e^{2i\pi\langle \bar{m}_{L_1}-m_{L'_1}-(\bar{m}_{L_2}-m_{L'_2}),.\rangle}
\end{equation}

\noindent which is continuous on $\mathbb{T}^d$. 
Moreover, let $L_0\in {\mathcal{L}'}$, then

\begin{equation}\nonumber P^{{\mathcal{L}'}}_{L_0}HP^{{\mathcal{L}'}}_{\bar{L}_0}
=\sum_{L,L'\in\mathcal{L},L\subset L_0, L'\subset \bar{L}_0} 
P^{{\mathcal{L}}}_{L}HP^{{\mathcal{L}}}_{{L'}}\end{equation}

\noindent which is continuous on $\mathbb{T}^d$. Thus $H$ has nice periodicity properties with respect to ${\mathcal{L}'}$.

\bigskip
\noindent By definition, ${\mathcal{L}'}$ is $A$-invariant. 
Moreover, assumption \eqref{parité} implies 

$$A'-A=\sum_{L'\in {\mathcal{L}'}}P^{{\mathcal{L}'}}_{L'}(A'-A)P^{{\mathcal{L}'}}_{L'}$$ 

\noindent so it also implies that ${\mathcal{L}'}$ is $A'-A$-invariant. Thus, ${\mathcal{L}'}$ is $A'$-invariant and so it is an $A'$-decomposition. $\Box$

\bigskip
\noindent Here is a standard lemma on the estimate of the rest of the Fourier series for an analytic function.

\begin{lem}\label{tronc}Let $H\in C^\omega_r(2\mathbb{T}^d,gl(n,\mathbb{C}))$. Soit $N\in\mathbb{N}$ 
and $H^N$ the truncation of $H$ at order $N$. Then for all $r'<r$, 

\begin{equation}|H-H^N|_{r'}\leq \frac{CN^d}{(r-r')^{d+1}}|H|_re^{-2\pi N(r-r')}
\end{equation}

\noindent where $C$ only depends on $d$. 

\end{lem}

\dem It is a simple computation. Since

\begin{equation}\nonumber H-H^N=\sum_{|m|>N}\hat{H}(m)e^{2i\pi \langle m,.\rangle}
\end{equation}

\noindent then

\begin{equation}\nonumber \begin{split}|H-H^N|_{r'}&\leq \sum_{|m|>N}||\hat{H}(m)||e^{2\pi |m| r'}
\leq |H|_r \sum_{|m|>N}e^{-2\pi |m|(r- r')}\\
&\leq C |H|_r \sum_{M>N} M^d e^{-2\pi M (r-r')}\leq C |H|_r \frac{N^d}{(r-r')^{d+1}}e^{-2\pi N(r-r')}\ \Box
\end{split}\end{equation}

\subsubsection{Inductive lemma}

\begin{prop}\label{iter} 
Let
\begin{itemize}
\item $\tilde{\epsilon}>0,\tilde{r}\leq 1$, $\tilde{r}'\in [\frac{\tilde{r}}{2},\tilde{r}[, \kappa'>0,\tilde{N}\in\mathbb{N},\gamma\geq n(n+1),C>0$; 

\item ${\tilde{F}}\in C^\omega_{\tilde{r}}(2\mathbb{T}^d,\mathcal{G}), \tilde{A}\in \mathcal{G}$,
\item $\mathcal{L}$ an $(\tilde{A},\kappa',\gamma)$-decomposition. 
\end{itemize}
There exists a constant $C''>0$ depending only on $\tau,n$ such that if

\begin{enumerate}


\item $\tilde{A}$ has $DC^{\tilde{N}}_\omega(\kappa',\tau)$ spectrum;

\item 

\begin{equation}\label{moypt}||\hat{\tilde{F}}(0)||\leq  \tilde{\epsilon}\leq C''\left(\frac{C^\tau\kappa'}{1+||\tilde{A}||}\right)^{2n}\end{equation}

\noindent and

\begin{equation}\label{cond2'} \tilde{N}
\leq  \frac{|\log \tilde{\epsilon}|}{C}
\end{equation} 

\item $\tilde{F}$ has nice periodicity properties with respect to $\mathcal{L}$

\end{enumerate}
 
\noindent then there exist
\begin{itemize}
\item $C'\in\mathbb{R}$ depending only on $n,d,\kappa,\tau$, 
\item $D\in\mathbb{N}$ depending only on $n,d,\tau$,
\item ${X}\in C^\omega_{\tilde{r}'}(2\mathbb{T}^d, \mathcal{G})$, 
\item ${A}'\in \mathcal{G}$ 
\item an $(A',\frac{3\kappa'}{4 },\gamma) $-decomposition $\mathcal{L}'$
\end{itemize}
satisfying the following properties: 

\begin{enumerate}

\item \label{nr} $A'$ has $DC^{\tilde{N}}_\omega(\frac{3\kappa'}{4},\tau)$ spectrum,

\item \label{0-} $||{A}'-\tilde{A}||\leq  \tilde{\epsilon}$;

\item \label{bpper-} the map
$F'\in C^\omega_{\tilde{r}'}(2\mathbb{T}^d,\mathcal{G})$ defined by

\begin{equation}\label{4-}
\forall \theta\in 2\mathbb{T}^d,\ \partial_\omega e^{{X}(\theta)}=(\tilde{A}+\tilde{F}(\theta))e^{{X}(\theta)}
-e^{{X}(\theta)}({A}'+{F}'(\theta))\end{equation}

\noindent has nice periodicity properties with respect to $\mathcal{L}'$

\item \label{9-} If $\Phi$ is trivial with respect to $\mathcal{L}$,


\noindent 
then

\begin{equation}\label{pxp}|\Phi^{-1}X\Phi|_{\tilde{r}'}
\leq C'\left(\frac{1+||\tilde{A}_{\mathcal{N}}||}{\kappa'(\tilde{r}-\tilde{r}')}\right)^{D\gamma}
|\Phi^{-1}\tilde{F}\Phi|_{\tilde{r}}\end{equation}

\item \label{10-} and if $\Phi$ is trivial with respect to $\mathcal{L}$, 


\begin{equation}\begin{split}|\Phi^{-1}F'\Phi|_{\tilde{r}'}&\leq C'\left(\frac{ 1+||\tilde{A}_{\mathcal{N}}||}{\kappa'
(\tilde{r}-\tilde{r}')}\right)^{D\gamma}
e^{|\Phi^{-1}X\Phi|_{\tilde{r}'}}|\Phi^{-1}\tilde{F}\Phi|_{\tilde{r}}\\
&(|\Phi|^2_{\tilde{r}}|\Phi^{-1}|^2_{\tilde{r}}\tilde{N}^de^{-2\pi \tilde{N}(\tilde{r}-\tilde{r}')}+|\Phi^{-1}\tilde{F}\Phi|_{\tilde{r}'}
(1+e^{|\Phi^{-1}X\Phi|_{\tilde{r}'}}))\end{split}\end{equation}

\end{enumerate}

\noindent Moreover, if $\tilde{F}$ is continuous on $\mathbb{T}^d$, then so are $X$ and $F'$. If $\mathcal{G}=o(n)$ or $u(n)$, then the same holds replacing condition \eqref{moypt} by

\begin{equation}\label{moypt'}||\hat{\tilde{F}}(0)||\leq  \tilde{\epsilon}\leq C''(C^\tau\kappa')^2\end{equation}

\end{prop}

\dem By assumption,   
$\tilde{F}$ has nice periodicity properties with respect to $\mathcal{L}$ and some family $(m_L)$ and $\tilde{A}$ has $DC^{\tilde{N}}_\omega(\kappa',\tau)$ spectrum, so one can apply Proposition \ref{homol}. Let ${X}\in C^\omega_{r'}(2\mathbb{T}^d,\mathcal{G})$ be a solution 
of

\begin{equation}\nonumber \forall \theta\in 2\mathbb{T}^d,\ 
\partial_\omega {X}(\theta)=[\tilde{A},{X}(\theta)]+\tilde{F}^{\tilde{N}}(\theta)
-\hat{\tilde{F}}(0)
\end{equation}

\noindent satisfying the conclusion of Proposition \ref{homol}.


\bigskip
\noindent Let ${A}'=\tilde{A}+\hat{\tilde{F}}(0)$. Then 
${A}'\in \mathcal{G}$ and $||\tilde{A}-A'||=||\hat{\tilde{F}}(0)||$, so property \ref{0-} holds. 




\noindent 
Moreover, let $c$ be the constant given by Lemma \ref{NRdurable}, 
and assume $C''\leq c $. 
Assumptions \eqref{moypt} and \eqref{cond2'} make it possible to apply Lemma \ref{NRdurable} and infer that $A'$ has $DC^{\tilde{N}}_\omega
(\frac{3\kappa'}{4},\tau)$ spectrum, thus property \ref{nr} holds. If $\mathcal{G}=o(n)$ or $u(n)$, one can apply lemma \ref{NRdurableorth} instead of lemma \ref{NRdurable} to get 
the same result with the weaker smallness condition \eqref{moypt'}.

\noindent Let $F'\in C^\omega_{r'}(2\mathbb{T}^d,\mathcal{G})$ the map defined in \eqref{4-}. Then

\begin{equation}\label{F'bar-}\begin{split}&{F}'=e^{-{X}}
(\tilde{F} -\tilde{F}^{\tilde{N}} )
+e^{-{X} }
\tilde{F} 
(e^{{X} }-Id)+(e^{-{X} }-Id)\hat{\tilde{F}}(0)
-e^{-{X} }\sum_{k\geq 2}\frac{1}{k!}\sum_{l=0}^{k-1}{X} ^l
(\tilde{F}^{\tilde{N}} -\hat{\tilde{F}}(0))
{X} ^{k-1-l}\end{split}\end{equation}

\noindent We shall appply Lemma \ref{bppG} with $A=\tilde{A}$ and $G=F'$, in order to get property \ref{bpper-}. 
The map ${F}'$ 
has nice periodicity properties with respect to $\mathcal{L}$ 
and some family $(m_L)$ 
since ${X}$ and 
$\tilde{F}$ have them. Moreover, as ${\tilde{F}}$ has nice periodicity properties with respect to
$\mathcal{L}$, 

\begin{equation}\nonumber P^{\mathcal{L}}_L\hat{\tilde{F}}(0)P^{\mathcal{L}}_{L'}\neq 0\Rightarrow 
P^{\mathcal{L}}_L\tilde{F}P^{\mathcal{L}}_{L'}\in C^0(\mathbb{T}^d)
\end{equation}

\noindent and since

\begin{equation}\nonumber 
P^{\mathcal{L}}_L\tilde{F}P^{\mathcal{L}}_{L'}\in C^0(\mathbb{T}^d)
\Rightarrow m_L-m_{L'} \in \mathbb{Z}^d
\Rightarrow P^{\mathcal{L}}_L{F}'P^{\mathcal{L}}_{L'}\in C^0(\mathbb{T}^d)
\end{equation}

%


\noindent then assumption \eqref{parité} of Lemma \ref{bppG} is fulfilled. By Lemma \ref{bppG}, ${F}'$ has therefore nice periodicity properties with respect to an $A'$-decomposition $\mathcal{L}'$ which is less fine than $\mathcal{L}$, so $\mathcal{L}'$ is an $(\tilde{A},
\kappa',\gamma)$-decomposition. As it is an $(\tilde{A},
\kappa',\gamma)$-decomposition, and by property \ref{0-}, each subspace $L\in\mathcal{L}'$ satisfies
%


\begin{equation}\nonumber  \mid \mid P^{\mathcal{L}'}_L\mid \mid 
\leq C_0 \left(\frac{1+\mid \mid A' _{\mathcal{N}}\mid \mid +2\tilde{\epsilon}}{\kappa'}\right)^{\gamma} 
\leq C_0 \left(\frac{1+\mid \mid {A}' _{\mathcal{N}}\mid \mid }{\frac{3\kappa'}{4}}\right)^{\gamma} 
\end{equation}

\noindent and so $\mathcal{L}'$ is an $(A',\frac{3\kappa'}{4},\gamma)$-decomposition, thus property \ref{bpper-} is satisfied.

\bigskip
\noindent Property \ref{9-} is given by Proposition \ref{homol}.

\bigskip
$\bullet$ By Lemma \ref{tronc},

\begin{equation}\nonumber \begin{split}|\tilde{F}-\tilde{F}^{\tilde{N}}|_{\tilde{r}'}
&\leq C_1\tilde{N}^d|\tilde{F}|_{\tilde{r}}\frac{e^{-2\pi \tilde{N}(\tilde{r}-\tilde{r}')}}{(\tilde{r}-\tilde{r}')^{d+1}}
\end{split}\end{equation}

\noindent where $C_1$ only depends on $d$. By \eqref{F'bar-}, \eqref{X2} and Lemma \ref{tronc}, it is true that

%



\begin{equation}\nonumber \begin{split}|\Phi^{-1}{F}'\Phi|_{\tilde{r}'}
&\leq C'\left(\frac{1+||\tilde{A}_{\mathcal{N}}||}{\kappa'(\tilde{r}-\tilde{r}')}\right)^{D\gamma}e^{|\Phi^{-1} X\Phi|_{\tilde{r}'}}
|\Phi^{-1}\tilde{F}\Phi|_{\tilde{r}}(|\Phi|^2_{\tilde{r}}|\Phi^{-1}|^2_{\tilde{r}}\tilde{N}^de^{-2\pi \tilde{N}(\tilde{r}-\tilde{r}')}\\
&+|\Phi^{-1}\tilde{F}\Phi|_{\tilde{r}'}(1+e^{|\Phi^{-1} X\Phi|_{\tilde{r}'}}))
\end{split}\end{equation}

\noindent where $C'$ only depends on $n,d,\kappa,\tau$ and $D$ only depends on $n,d,\tau$, whence property \ref{10-}. $\Box$

\subsection{Inductive step}

\noindent Now we are able to state the whole inductive step. In the following we will denote 

\begin{equation}\label{kappa''}\left\{\begin{array}{c}
N(r,\epsilon)=\frac{1}{2\pi r}|\log \epsilon|\\
R(r,r')= \frac{1}{(r-r')^8}80^{4}(\frac{1}{2}n(n-1)+1)^2\\
\kappa''(r,r',\epsilon)=\frac{\kappa}{n(8R(r,r')^{\frac{1}{2}n(n-1)+1}{N(r,\epsilon)})^{\tau}}
\end{array}\right.\end{equation}
%

%


%




\begin{prop}\label{iter3}
Let
\begin{itemize}
\item $A\in \mathcal{G}$,
\item $r\leq \frac{1}{2},r''\in [\frac{95}{96}r,r[$,$\gamma\geq n(n+1)$, 
\item $\bar{A},\bar{F}\in C^\omega_{r}(2\mathbb{T}^d,\mathcal{G})$ and $\Psi\in C^\omega_{r}(2\mathbb{T}^d,G)$, 
\item $\epsilon=|\bar{F}|_{r}$,
\end{itemize}

\noindent There exists $\tilde{C}'>0$ depending only on $n,d,\kappa,\tau,\gamma$ and there exists $D_3\in\mathbb{N}$ depending only on 
$n,d,\tau$ 
such that if

\begin{enumerate}

\item $\bar{A}$ is reducible to $A$ 
by $\Psi$, 

\item \label{bpperinput} $\Psi^{-1}\bar{F}\Psi$ has nice periodicity properties with respect to an $(A,\kappa''(r,r'',\epsilon),\gamma)$-decomposition $\mathcal{L}$

\item

\begin{equation}\label{epsA+} \epsilon 
\leq \frac{\tilde{C}'}{(||A||+1)^{D_3\gamma}}(r-r'')^{D_3\gamma}
\end{equation}

\item $|\Psi|_{r}\leq (\frac{1}{\epsilon})^{-\frac{1}{2}(r-r'')}$ et 
$|\Psi^{-1}|_{r}\leq (\frac{1}{\epsilon})
^{-\frac{1}{2}(r-r'')}$, 


\end{enumerate}
 
\noindent then there exist
\begin{itemize}
\item $\epsilon'\in [\epsilon^{R(r,r'')^{n^2}}, \epsilon^{100}]$;
\item $Z'\in C^\omega_{r''}(2\mathbb{T}^d,G)$, 
\item $\bar{A}',\bar{F}'\in C^\omega_{r''}(2\mathbb{T}^d,\mathcal{G})$, 
\item $\Psi'\in C^\omega_{r}(2\mathbb{T}^d, G)$, 
\item ${A}'\in \mathcal{G}$ 

\end{itemize}
\noindent satisfying the following properties:

\begin{enumerate}

\item \label{5} $\bar{A}'$ is reducible by $\Psi'$ to ${A}'$,

\item \label{bpper+} the map $\Psi'^{-1}\bar{F}'\Psi'$ has nice periodicity properties with respect to an 
 $(A',{\kappa''(r'',r''-\frac{r-r''}{2},\epsilon')},2\gamma))$-decomposition 
$\mathcal{L}'$

\item \label{8+} $|\bar{F}'|_{r''} \leq  \epsilon '$,

\item \label{ren+} $|\Psi'|_{r''}\leq (\frac{1}{\epsilon'})^{\frac{1}{4}(r-r'')}$ and
$|\Psi'^{-1}|_{r''}\leq 
(\frac{1}{\epsilon'})^{\frac{1}{4}(r-r'')}$,

\item \label{0} $||A'||\leq ||A||+\mid \log \epsilon\mid\left(\frac{1}{r-r'}\right) ^{D_3}$;

\item \label{4+}

\begin{equation}\partial_\omega Z'=(\bar{A}+\bar{F})Z'
-Z'(\bar{A}'+\bar{F}')\end{equation}

\item \label{6} 

\begin{equation}|Z'-Id|_{r''}\leq \frac{1}{\tilde{C}'}\left(\frac{(1+||A||)|\log\epsilon|}
{r-r''}\right)^{D_3\gamma}
\epsilon^{1-4(r-r'')}
\end{equation}

\noindent and so does $(Z')^{-1}-Id$.
%


\end{enumerate}

\noindent Moreover, 
\begin{itemize}
\item in dimension 2, if $\bar{A},\bar{F}$ are continuous on $\mathbb{T}^d$, 
and if assumption \ref{bpperinput} is replaced by

\bigskip
\ref{bpperinput}'. $\Psi$ is such that for all function $H$ continuous on 
$\mathbb{T}^d$, $\Psi H \Psi^{-1}$ is continuous on $\mathbb{T}^d$, 

\bigskip
\noindent then 
$Z',\bar{A}',\bar{F}'$ are continuous on $\mathbb{T}^d$ and property \ref{bpper+} is replaced by

\bigskip
\ref{bpper+}'. $\Psi'$ is such that for every function $H$ continuous on 
$\mathbb{T}^d$, $\Psi' H \Psi'^{-1}$ is continuous on $\mathbb{T}^d$.

\item If $\mathcal{G}=gl(n,\mathbb{C})$ or $u(n)$ and if $\bar{A},\bar{F},\Psi$ are continuous on $\mathbb{T}^d$, then 
$Z',\bar{A}',\bar{F}',\Psi'$ are continuous on $\mathbb{T}^d$. 

\item if $\mathcal{G}=o(n)$ or $u(n)$, the same holds with the weaker condition

\begin{equation}\label{epsorth+} \epsilon 
\leq \tilde{C}'(r-r'')^{D_3\gamma}
\end{equation}

\noindent instead of \eqref{epsA+};

\item if $\mathcal{G}=sl(2,\mathbb{C})$ or $sl(2,\mathbb{R})$, then either $\Psi'^{-1}\Psi$ is the identity or $\mid \mid A'\mid \mid\leq \kappa''(r,r'',\epsilon) +\epsilon^{\frac{1}{2}}$.

\end{itemize}

\end{prop}

\bigskip
\noindent  The proof will be made in two steps: the first step is to 
reduce the perturbation when there are resonances. 
The second step is to iterate Proposition \ref{iter} as many times as possible using the fact that resonances, once removed, do not reappear immediately. 

\paragraph{First step: removing the resonances}

\noindent 
Let $r'=\frac{r+r''}{2}$. 
Let $R=R(r,r');N=N(r,\epsilon);\kappa''=\kappa''(r,r',\epsilon)$. Let $\bar{N}$ be given by Proposition \ref{diophan} and $\Phi$ a map of reduction of the eigenvalues of $A$ at order $R,\bar{N}$. Let $\Psi'=\Psi \Phi$ and $\tilde{F}=(\Psi')^{-1}\bar{F}\Psi'$. 

\noindent We shall apply Proposition \ref{iter} with

\begin{equation}\nonumber \begin{split}&
\tilde{\epsilon}=\epsilon^{1-{2(r-r')-\frac{1}{48}}},\ 
\tilde{r}=r,\ 
\tilde{r}'=r',\ 
\kappa'=\frac{\kappa''}{C_0 }
, \ \tilde{N}=R\bar{N}, C=\frac{2\pi r}{R^{\frac{1}{2}n(n-1)+1 } } 
\end{split}\end{equation}

and $\tilde{A}\in \mathcal{G}$ 
such that

\begin{equation}\nonumber \forall \theta\in 2\mathbb{T}^d,\ \partial_\omega \Phi(\theta)=A\Phi(\theta)-\Phi(\theta)\tilde{A}
\end{equation}

\noindent Let $C''$ be given by Proposition \ref{iter} (depending only on $n$ and $\tau$).

\bigskip
\noindent The matrix $\tilde{A}$ has $DC_\omega^{R\bar{N}}(\kappa'',\tau)$ spectrum. 
By assumption, $\Psi^{-1}\bar{F}\Psi$ has nice periodicity properties with respect to an $(A,\kappa'',\gamma)$-decomposition $\mathcal{L}$ and some family $(m_L)$. 
Moreover $\Phi$ is trivial with respect to $\mathcal{L}_{A,\kappa''}$. 
Since $\mathcal{L} $ is an $A$-decomposition, there is a Jordan decomposition which is finer than $\mathcal{L} $; and since $\mathcal{L}_{A,\kappa''}$ is less fine than any Jordan decomposition, one can define an $A$-decomposition $\bar{\mathcal{L}}$ in the following way: 

\begin{equation}\nonumber L\in \bar{\mathcal{L}} \Leftrightarrow \exists L_1\in \mathcal{L},L_2\in \mathcal{L}_{A,\kappa''}\ \mid \ 
L=L_1\cap L_2
\end{equation}

\noindent 
$\bar{\mathcal{L}}$ is an $(A,\frac{\kappa''}{C_0},2\gamma)$-decomposition since 
$\mathcal{L} $ and $\mathcal{L}_{A,\kappa''}$ are $(A,\kappa'',\gamma)$-decompositions and 
$\tilde{F}$ has nice periodicity properties with respect to $\bar{\mathcal{L}}$. Since $\bar{\mathcal{L}}$ is an $(A,\frac{\kappa''}{C_0},2\gamma)$-decomposition, it is also an 
$(\tilde{A},\frac{\kappa''}{C_0},2\gamma)$-decomposition (because the nilpotent parts of $A$ and $\tilde{A}$ coincide, and because any Jordan decomposition for $A$ is a Jordan decomposition for $\tilde{A}$).

\noindent Moreover,

\begin{equation}\nonumber ||\hat{\tilde{F}}(0)||\leq |\tilde{F}|_{0}\leq |\Phi|_0|\Phi^{-1}|_0
|\Psi|_0|\Psi^{-1}|_0|\bar{F}|_0 \end{equation}

\noindent Now by \eqref{normephi}, for all $s'\geq 0$,

\begin{equation}\label{phis'}\mid \Phi \mid_{s'} \leq C_0\left(\frac{1+||A_{\mathcal{N}}||}{\kappa''}\right)^{n(n+1)}e^{4\pi \bar{N}s'}
\end{equation}

\noindent and so does $\Phi^{-1}$.
Thus

%


\begin{equation}\nonumber ||\hat{\tilde{F}}(0)|| \leq \epsilon^{1-{2(r-r')}}C_0^2\left(\frac{1+||A_{\mathcal{N}}||}{\kappa''}\right)^{2n(n+1)}\end{equation}

\noindent therefore, if 
$\tilde{C}'\leq C_0^{96}$ and $D_3\gamma \geq 96n(n+1)$, then

\begin{equation}\nonumber ||\hat{\tilde{F}}(0)|| \leq \epsilon^{1-2(r-r')-\frac{1}{48}}\end{equation}



\noindent
Assumption \eqref{epsA+}, which implies \eqref{moypt} with 

\begin{equation}\nonumber \tilde{C}'\leq C''^4\left(\frac{C}{(r-r')^{4n(n-1)+9}}\right)^{8n\tau},\ 
D_3\gamma \geq 64n(n(n-1)+2)\tau\end{equation}

\noindent (note that 
$\frac{C}{(r-r')^{4n(n-1)+9}}$ has a lower bound which is independent of $r-r'$), together with the choice of $\tilde{N}$ 
which implies \eqref{cond2'},  
%
make it possible to apply Proposition \ref{iter} to obtain $C'>0$ depending only on $n,d,\kappa,\tau$, $D\in \mathbb{N}$ depending only on $n,d,\tau$ and 
functions 
$X\in C^\omega_{r'}(2\mathbb{T}^d,\mathcal{G})$, 
${F}_1\in C^\omega_{r'}(2\mathbb{T}^d,\mathcal{G})$, and a matrix ${A}_1\in \mathcal{G}$ 
such that

\begin{itemize}

\item $A_1$ has $DC^{R\bar{N}}_\omega(\frac{3}{4}\left(\frac{\kappa''}{C_0 }\right),\tau)$ spectrum 

\item $||A_1-\tilde{A}||\leq  \epsilon^{\frac{23}{24}}$, which implies

\begin{equation}\label{1/96}\mid \mid A_1-A\mid \mid \leq ||A_1-\tilde{A}||+||A-\tilde{A}||\leq  \epsilon^{\frac{23}{24}}+4\pi \bar{N}\end{equation}

\noindent 
If $\mathcal{G}=sl(2,\mathbb{C})$ or $sl(2,\mathbb{R})$, then 

\begin{equation}\nonumber ||A_1||\leq \mid \mid \tilde{A}  \mid \mid +\epsilon^{\frac{23}{24}}\leq \kappa'' +\epsilon^{\frac{23}{24}}
\end{equation}

\item $\partial_\omega e^X=(\tilde{A}+\tilde{F})e^X
-e^X({A}_1+{F}_1)
$,

\item ${F}_1$ has nice periodicity properties with respect to an $(A_1,\frac{3\kappa''}{4C_0 },2\gamma)$-decomposition $\mathcal{L}'$

\item

\noindent and since $\Phi$ is trivial with respect to $\bar{\mathcal{L}}$, 



\begin{equation}\nonumber 
|\Phi X\Phi^{-1}|_{r'}\leq C'\left(\frac{C_0(1+||A_{\mathcal{N}}||)}{\kappa''(r-r')}\right)^{D\gamma}
|\Phi\tilde{F}\Phi^{-1}|_r\tag{\ref{pxp}}
\end{equation}

\noindent and

\begin{equation}\label{pfp}\begin{split}|\Phi{F}_1\Phi^{-1}|_{r'}
&\leq C'\left(\frac{C_0(1+||A_{\mathcal{N}}||)}{\kappa''(r-r')}\right)^{D\gamma}e^{|\Phi X\Phi^{-1}|_{r'}}
|\Phi\tilde{F}\Phi^{-1}|_{r}
(|\Phi|^2_r|\Phi^{-1}|^2_r(R\bar{N})^de^{-2\pi R\bar{N}(r-r')}\\
&+|\Phi\tilde{F}\Phi^{-1}|_{r'}(1+e^{|\Phi X\Phi^{-1}|_{r'}}))\end{split}\end{equation}

\end{itemize}

\noindent Now

\begin{equation}\nonumber |\Phi\tilde{F}\Phi^{-1}|_{r}\leq |\Psi|_r|\Psi^{-1}|_r|\bar{F}|_r
\leq \epsilon^{1-2(r-r')}\end{equation}


%

\noindent so, by 
\eqref{epsA+}, if $D_3$ is great enough as a function of $n,\gamma,D$, then









\begin{equation}\nonumber |\Phi{F}_1\Phi^{-1}|_{r'}\leq \epsilon^{-\frac{1}{96}} \epsilon^{1-2(r-r')}((R\bar{N})^d\epsilon^{100}
+\epsilon^{1-2(r-r')})
\end{equation}

\noindent There exists a constant $c_d$ which only depends on $D,\gamma,\tau$ such that if $\epsilon \leq c_d$, then
%


\begin{equation}\nonumber \epsilon^{\frac{1}{2}}\mid \log \epsilon\mid^{D\gamma\tau} \leq 1
\end{equation}

%

%

\noindent thus if $\tilde{C}'$ is small enough and $D_3$ big enough (as a function of $n,d,\gamma,\tau$),

\begin{equation}\nonumber |\Phi{F}_1\Phi^{-1}|_{r'}
\leq  \epsilon^{2-4(r-r')-\frac{1}{96}}
\end{equation}




%


%







%

\bigskip
\noindent The estimate \eqref{phis'}, the assumption \eqref{epsA+} and the fact that $\mid \mid A_{\mathcal{N}}\mid \mid \leq \mid \mid A\mid \mid$, 
imply that



\begin{equation}\nonumber |\Psi\Phi|_{r}
\leq |\Psi|_r |\Phi|_r\leq \epsilon^{-(r-r')-\frac{1}{96}} e^{4\pi r\bar{N}}
\end{equation}


%

%




\bigskip
\noindent We shall estimate $|\Psi\Phi e^X(\Psi\Phi) ^{-1}-Id|_{r'}$. The estimate \eqref{pxp} implies 

\begin{equation}\nonumber \begin{split}|\Phi e^X \Phi^{-1} -Id|_{r'}
&\leq C''\left(\frac{(1+||A_{\mathcal{N}}||)R^{\frac{1}{2}(n(n-1)+1)\tau }{N}^{\tau }}{r-r'}\right)^{D\gamma}|\Phi\tilde{F}\Phi^{-1}|_r
\end{split}\end{equation}


\noindent for some $C''$ only depending on $n,d,\kappa,\tau$, so
%


\begin{equation}\nonumber \begin{split}|\Psi\Phi e^X(\Psi\Phi) ^{-1}-Id|_{r'}&\leq C_3\left(\frac{(1+||A_{\mathcal{N}}||)|\log\epsilon|}
{r-r'}\right)^{D'_1\gamma}|\bar{F}|_r
(\frac{1}{\epsilon})^{4(r-r')}
\end{split}\end{equation}

\noindent for some $C_3$ depending only on $n,d,\kappa,\tau$ and $D'_1$ depending only on $n,d,\tau$. 
The same estimate holds for $|\Psi\Phi e^{-X}(\Psi\Phi) ^{-1}-Id|_{r'}$.

%



%


\bigskip
\noindent Let $\bar{F}_1=\Psi\Phi{F}_1(\Psi\Phi)^{-1}$ 
and let $\bar{A}_1\in 
C^\omega_r(2\mathbb{T}^d,\mathcal{G})$ such that

\begin{equation}\nonumber \partial_\omega \Psi\Phi=\bar{A}_1\Psi\Phi-\Psi\Phi {A}_1
\end{equation}

Thus 
we have obtained
\begin{itemize}
\item $\bar{N}\in [N,R^{\frac{1}{2}n(n-1)}N]$,
\item $Z_1,\Psi'\in C^\omega_{r'}(2\mathbb{T}^d,G)$,
\item $A_1\in \mathcal{G}$ 
\item $\bar{A}_1\in C^\omega_{r'}(2\mathbb{T}^d,\mathcal{G})$ 

\item and ${F}_1=(\Psi')^{-1}\bar{F}_1\Psi'$ 
\end{itemize}
\noindent such that

\begin{enumerate}

\item $\bar{A}_1$ is reducible to $A_1$ by $\Psi'$ 

\item $F_1$ has nice periodicity properties with respect to an $(A_1,\frac{3\kappa''}{4C_0},2\gamma)$-decomposition $\mathcal{L}_1$

\item $\mid \Psi'\mid_{r'}\leq \epsilon^{-(r-r')-\frac{1}{96}}e^{4\pi r\bar{N}}$ and $\mid \Psi'^{-1}\mid_{r'}\leq \epsilon^{-(r-r') -\frac{1}{96}}e^{4\pi r\bar{N}}$

\item $A_1$ has $DC^{R\bar{N}}_\omega (\frac{3}{4}\kappa'',\tau)$ spectrum, 

\item $\partial_\omega Z_1=(\bar{A}+\bar{F})Z_1
-Z_1(\bar{A}_1+\bar{F}_1)$,

\item $||A_1||\leq ||A||+\epsilon ^{\frac{23}{24}}+4\pi \bar{N}$, and, if $\mathcal{G}=sl(2,\mathbb{C})$ or $sl(2,\mathbb{R})$ and $\Psi'^{-1}\Psi$ is not the identity, $||A_1||\leq 
\kappa''(r,r'',\epsilon)+\epsilon ^{\frac{23}{24}}$;

\item

\begin{equation}\label{Z1}|Z_1-Id|_{r'}\leq \frac{1}{\tilde{C}'}\left(\frac{(1+||A_{\mathcal{N}}||)|\log\epsilon|}
{r-r'}\right)^{D_1\gamma}\epsilon^{1-4(r-r')}
\end{equation}

\noindent and so does $|Z_1^{-1}-Id|_{r'}$, 

\item 

\begin{equation}\label{oulala}
|\Psi^{-1}\bar{F}_1\Psi|_{r'}\leq \epsilon^\frac{3}{2}\end{equation}

\item $\Psi'^{-1}\Psi$ is trivial with respect to $\mathcal{L}_{A,\kappa''}$,

\item and for every $s'\geq 0$,

\begin{equation}\label{psps'}\mid \Psi'^{-1}\Psi\mid _{s'} \leq C_n\left(
\frac{1+\mid \mid A_{\mathcal{N}}\mid \mid }{\kappa''}
\right)^{n(n+1)}e^{4\pi \bar{N}s'}
\end{equation}

\noindent and so does $\mid \Psi^{-1}\Psi'\mid _{s'}$, where $C_n$ only depends on $n$.

\end{enumerate}

\paragraph{Second step: iteration far from resonances}
Let $l$ such that

\begin{equation}\nonumber \epsilon^{(\frac{4}{3})^{l+1}}\leq  e^{-2\pi (r-r'')\sqrt[4]{R}\bar{N}}\leq \epsilon^{(\frac{4}{3})^{l}}\end{equation}

\noindent Let $\epsilon'=e^{-2\pi (r-r'')\sqrt[4]{R}\bar{N}}$. Define the sequence
$\epsilon_j=\epsilon^{(\frac{3}{2})^{j}-\frac{1}{48}}$. We shall iterate $l-1$ times Proposition \ref{iter}, starting with $j=2$, with 
\begin{itemize}
\item $\tilde{\epsilon}=\epsilon_{j-1}$,
\item $C=\left(\frac{r-r''}{160(\frac{1}{2}n(n-1)+1)}\right)^{8(\frac{1}{2}n(n-1)+1)}$
\item $\tilde{r}=r_{j-2}=\frac{r+r''}{2}-(j-2)\frac{r-r''}{2l}$,
\item $\tilde{r}'=r_{j-1}=\frac{r+r''}{2}-(j-1)\frac{r-r''}{2l}$, 
\item $\kappa'=(\frac{3}{4})^{j-1}\frac{\kappa''}{C_0 }$,
\item $\tilde{N}=R\bar{N}$,

\item $\tilde{F}=F_{j-1}$,
\item $\tilde{A}=A_{j-1}$,
\item $\Phi=\Psi^{-1}\Psi'$,
\item $\mathcal{L}=\mathcal{L}_1$,
\end{itemize}

\bigskip
\noindent Note that for every $j$, 

\begin{equation}\nonumber \epsilon_{j} \leq C''\left(\frac{ C^\tau(\frac{3}{4})^{j}\frac{\kappa''}{C_0 }}{1+||A_{1}||+\sum_{l=1}^{j-1}\epsilon_l}\right)^{2n}
\end{equation}

\noindent Estimates \eqref{oulala} and \eqref{psps'} imply

\begin{equation}\nonumber ||\hat{F}_1(0)||\leq \mid F_1\mid_0\leq \mid \Psi'^{-1} \Psi\mid _0 \  \mid \Psi^{-1} \Psi'\mid _0 \ |\Psi^{-1}\bar{F}_1\Psi|_0
\leq  C_n^2\left(\frac{1+||A_{\mathcal{N}}||}{\kappa''}\right)^{2n(n+1)}\epsilon^\frac{3}{2}
\leq \epsilon^{\frac{3}{2}-\frac{1}{48}}\end{equation}

\noindent Moreover, $A_1$ has $DC^{R\bar{N}}_\omega(\frac{3}{4}\kappa'',\tau)$ spectrum and $F_1$ has nice periodicity properties with respect to $\mathcal{L}$. Let $C''$ be the constant given by Proposition \ref{iter}. 
By assumption on  $\epsilon$, with $C'$ depending only on $n,d,\kappa,\tau$ and $D_3$ depending only on $n,\tau$, one has
%

\noindent 

\begin{equation}\nonumber \begin{split}\tilde{\epsilon}
& \leq \frac{C''}{(1+||A_1||)^{2n}} \left(\frac{3\kappa''}{4C_0 }\right)^{2n} C^{2n\tau} 
\end{split}\end{equation}

\noindent Moreover,

\begin{equation}\nonumber R\bar{N}\leq R^{{n_0}+1}N\leq \frac{1}{C}|\log\epsilon|
\end{equation}

\noindent  so the assumptions \eqref{moypt} and \eqref{cond2'} of Proposition \ref{iter} hold with 
$\tilde{F}=F_1, \kappa'=\kappa'', \tilde{N}=R\bar{N}$.


\bigskip
\noindent 
Fix $j$ and assume ${A}_{j-1}$ has $DC^{\tilde{N}}_\omega(\kappa',\tau)$ spectrum, $F_{j-1}$ has nice periodicity properties with respect to
an $(A_{j-1},(\frac{3}{4})^{j-1}\frac{\kappa''}{C_0 },2\gamma)$-decomposition, 

\begin{equation}\nonumber 
\mid \mid \hat{F}_{j-1}(0)\mid \mid \leq \epsilon_{j-1}\end{equation}

\noindent and

\begin{equation}\nonumber CR\bar{N} \leq \mid \log \epsilon_{j-1}\mid 
\end{equation}

\noindent One obtains functions $F_j,X_j$ and a matrix $A_j$ 
such that
\begin{enumerate}

\item $A_j$ has $DC^{R\bar{N}}((\frac{3}{4})^j\frac{\kappa''}{C_0},\tau)$ spectrum, 
\item $||A_j||\leq ||A_{j-1}||+\epsilon_{j-1}$, 

\item

\begin{equation}\nonumber \partial_\omega e^{X_j}=(A_{j-1}+F_{j-1})e^{X_j}-e^{X_j}(A_j+F_j)\end{equation}

\noindent and
$F_j$ has nice periodicity properties with respect to an $(A_j,(\frac{3}{4})^j\frac{\kappa''}{C_0},2\gamma)$-decomposition

\item

\begin{equation}\label{blabla}\mid \Psi^{-1}\Psi' X_j\Psi'^{-1}\Psi\mid_{r_{j-1}}\leq
C'\left(\frac{1+||(A_{j-1})_{\mathcal{N}}||}{\kappa''(r_{j-2}-r_{j-1})}\right)^{D\gamma}| \Psi^{-1}\Psi'F_j\Psi'^{-1}\Psi|_{r_{j-1}}
\end{equation}

\noindent for some $C'$ depending only on $n,d,\kappa,\tau$ and some
$D$ depending only on $n,d,\tau$,

\item and

\begin{equation}\label{bla}\begin{split}| \Psi^{-1}\Psi' F_j\Psi'^{-1}\Psi |_{r_{j-1}}&\leq C'\left(\frac{1+||(A_{j-1})_{\mathcal{N}}||}{\kappa''
(r_{j-2}-r_{j-1})}\right)^{D\gamma}
e^{|\Psi^{-1}\Psi' X_{j-1}\Psi'^{-1}\Psi|_{r_{j-2}}}|  \Psi^{-1}\Psi'F_{j-1}\Psi'^{-1}\Psi|_{r_{j-2}}\\
&(\mid \Psi'^{-1}\Psi\mid_{r_{j-2}}^4(R\bar{N})^de^{-2\pi R\bar{N}(r_{j-2}-r_{j-1})}\\
&+(1+2e^{|  \Psi^{-1}\Psi'X_{j-1}\Psi'^{-1}\Psi|_{r_{j-2}}})
| \Psi^{-1}\Psi' F_{j-1}\Psi'^{-1}\Psi|_{r_{j-2}})
\end{split}\end{equation}

\end{enumerate}

\noindent We shall bound $\mid \mid \hat{F}_j(0)\mid \mid$ to iterate Proposition \ref{iter}. Estimates \eqref{oulala} and \eqref{epsA+} imply
%


\begin{equation}\label{n0'}\begin{split}|\Psi^{-1}\Psi' F_j\Psi'^{-1}\Psi|_{r_{j-1}}&\leq 3C'\left(\frac{1+||(A_{j-1})_{\mathcal{N}}||}{\kappa''(r_{j-2}-r_{j-1})}
\right)^{D\gamma}
|\Psi^{-1}\Psi'F_{j-1}\Psi'^{-1}\Psi|_{r_{j-2}}\\
&(\mid \Psi'^{-1}\Psi\mid_{r_{j-2}}^4(R\bar{N})^de^{-2\pi R\bar{N}(r_{j-2}-r_{j-1})}+|\Psi^{-1}\Psi'F_{j-1}\Psi'^{-1}\Psi|_{r_{j-2}})\\
\end{split}\end{equation}

\noindent and since $r_{j-2}-r_{j-1}=\frac{r-r''}{2l}$,

\begin{equation}\label{n0}\begin{split}|  \Psi^{-1}\Psi'F_j\Psi'^{-1}\Psi|_{r_{j-1}}
&\leq |\Psi^{-1}\Psi'F_{j-1}\Psi'^{-1}\Psi|_{r_{j-2}}^{\frac{3}{4}}(\mid \Psi'^{-1}\Psi\mid_{r_{j-2}}^4(R\bar{N})^de^{-2\pi \frac{R\bar{N}(r-r'')}{2l}}+|\Psi^{-1}\Psi'F_{j-1}\Psi'^{-1}\Psi|_{r_{j-2}})\\
&\leq |\Psi^{-1}\Psi'F_{j-1}\Psi'^{-1}\Psi|_{r_{j-2}}^{\frac{3}{4}}(\mid \Psi'^{-1}\Psi\mid_{r_{j-2}}^4(R\bar{N})^d\epsilon'^{ \frac{R^{\frac{3}{4}}}{2l}}+|\Psi^{-1}\Psi'F_{j-1}\Psi'^{-1}\Psi|_{r_{j-2}})
\end{split}\end{equation}
%




\noindent Now $l$ is bounded by

\begin{equation}\label{l'}\begin{split}
l
&\leq 8 (\frac{1}{2}n(n-1)+1)\sqrt[4]{R}
\end{split}\end{equation}

\noindent Moreover

\begin{equation}\nonumber \mid \Psi'^{-1}\Psi\mid _{r_{j-2}} \leq C_n\left(
\frac{1+\mid \mid A\mid \mid }{\kappa''}
\right)^{n(n+1)}e^{4\pi \bar{N}r_{j-2}}
\leq C_n\left(
\frac{1+\mid \mid A\mid \mid }{\kappa''}
\right)^{n(n+1)}\epsilon'^{- \frac{2r_{j-2}}{(r-r'')\sqrt[4]{R}}}
\end{equation}

\noindent and so

\begin{equation}\label{n0''}\begin{split}| \Psi^{-1}\Psi'F_j\Psi'^{-1}\Psi|_{r_{j-1}}
&\leq |\Psi^{-1}\Psi'F_{j-1}\Psi'^{-1}\Psi|_{r_{j-2}}^{\frac{3}{4}}(\epsilon'+|\Psi^{-1}\Psi'F_{j-1}\Psi'^{-1}\Psi|_{r_{j-2}})\\
&\leq |\Psi^{-1}\Psi'F_{j-1}\Psi'^{-1}\Psi|_{r_{j-2}}^{\frac{3}{2}}
\end{split}\end{equation}

\noindent By a simple induction, for every $j$,

\begin{equation}\label{fj}
|\Psi^{-1}\Psi'  F_j\Psi'^{-1}\Psi|_{r_{j-1}}\leq 
|\Psi^{-1}\Psi'  F_{1}\Psi'^{-1}\Psi|_{r_0}^{(\frac{3}{2})^{j-1}}\leq \epsilon^{(\frac{3}{2})^j}\end{equation}

\noindent Finally

\begin{equation}\nonumber \mid \mid \hat{F}_j(0)\mid \mid \leq |\Psi^{-1}\Psi'F_{j}\Psi'^{-1}\Psi|_{r_{j-1}} \mid \mid \Psi^{-1}\Psi' \mid \mid_0 \ \mid \mid \Psi'^{-1}\Psi \mid \mid_0
\leq  \epsilon_j
\end{equation}

\noindent so it is possible to iterate Proposition \ref{iter}.

\paragraph{Conclusion}

\noindent 
After $l-1$ steps,

\begin{equation}\nonumber 
|\Psi^{-1}\Psi'  F_{l+1}\Psi'^{-1}\Psi|_{r_{l}}\leq \epsilon'^{\frac{17}{16}}\end{equation}

\noindent Let
$Z=e^{X_2}\dots e^{X_{{l+1}}}\in C^\omega_{r'}(2\mathbb{T}^d,G)),{A}'=A_{l+1},{F}'=F_{l+1}$. Then

\begin{equation}\nonumber \partial_\omega Z=({A}_1+{F}_1)Z
-Z({A}'+{F}')\end{equation}

\noindent and

\begin{equation}\nonumber ||A'||\leq ||A_1||+\sum_{j=1}^{l}||\hat{F}_j(0)||+4\pi \bar{N}
\leq ||A||+\mid \log \epsilon\mid\left(\frac{1}{r-r'}\right) ^{D_4}
\end{equation}

\bigskip
\noindent for $D_4$ great enough depending only on $n$, whence property \ref{0}. 
If $\mathcal{G}=sl(2,\mathbb{C})$ or $sl(2,\mathbb{R})$ and $\Psi'^{-1}\Psi$ is not the identity, then

\begin{equation}\nonumber ||A'||\leq ||A_1||+\sum_{j=1}^{l}||\hat{F}_j(0)||
\leq \kappa''(r,r'',\epsilon)+\epsilon^{\frac{1}{2}}
\end{equation}

\bigskip To prove that $\mathcal{L}_{l+1}$ is indeed an 
$(A_{l+1},\kappa''(r'',r-\frac{r-r''}{2},\epsilon'),2\gamma)$-decomposition, it is enough to show that

\begin{equation}\nonumber \kappa''(r'',r''-\frac{r-r''}{2},\epsilon') \leq (\frac{3}{4})^{l+1}\frac{\kappa''}{C_0 }
\end{equation}

\noindent which comes from the definition of the function $\kappa''$.

%





%

%




%


%


\bigskip
\noindent Let us prove property \ref{ren+}. 
It is true that

%


\begin{equation}\label{(r-r')^2}|\Psi'|_{r''}\leq \epsilon^{-\frac{1}{2}(r-r'')}
\epsilon^{-\frac{1}{96}}e^{4\pi r\bar{N}}
\leq \epsilon^{-\frac{1}{2}(r-r'')}\epsilon^{-\frac{1}{96}}
\epsilon'^{-\frac{r-r''}{200}}
\end{equation}

\noindent and property \ref{ren+} comes from it, since

\begin{equation}\nonumber \epsilon=\epsilon'^{\frac{|\log \epsilon|}{2\pi \sqrt[4]{R}\bar{N}(r-r'')}}\end{equation}
%

%


%



\bigskip
\noindent Moreover,

\begin{equation}\nonumber |\Psi'{F}'\Psi'^{-1}|_{r''}\leq |\Psi|_r|\Psi^{-1}|_r|\Psi^{-1}\Psi' {F}'\Psi'^{-1}\Psi|_{r''}\leq  \epsilon'\end{equation}

\noindent whence \ref{8+}. Let $Z'=Z_1 \Psi' Z \Psi'^{-1}$, $\bar{F}'=\Psi'F'\Psi^{-1}$ (which satisfies property \ref{bpper+}) and $\bar{A}$ such that

\begin{equation}\nonumber \partial_\omega \Psi'=\bar{A}'\Psi'-\Psi' A'
\end{equation}

\noindent
Then \ref{4+} and \ref{5} hold, and by \eqref{Z1},

%


\begin{equation}\nonumber \begin{split}|Z'-Id|_{r''}&\leq |Z_1-Id|_{r_1}+ |\Psi|_r|\Psi^{-1}|_r\sum_j|\Psi^{-1}\Psi'  X_j\Psi'^{-1}\Psi|_{r_j}\\
&\leq \frac{1}{\tilde{C}'}\left(\frac{l(1+||A_{\mathcal{N}}||)|\log\epsilon|}
{r-r''}\right)^{D_1\gamma}
(\frac{1}{\epsilon})^{4(r-r'')}(\epsilon+\sum_j|\Psi^{-1}\Psi'  F_j\Psi'^{-1}\Psi|_{r_j})
\end{split}\end{equation}

\noindent and by \eqref{oulala} and \eqref{fj},

\begin{equation}\nonumber 
|Z'-Id|_{r''}\leq  \frac{2}{\tilde{C}'}\left(\frac{l(1+||A_{\mathcal{N}}||)|\log\epsilon|}
{r-r''}\right)^{D_1\gamma}
(\frac{1}{\epsilon})^{4(r-r'')}\epsilon
\end{equation}

\noindent whence property \ref{6} with $D_3\gamma\geq 2D_1\gamma$ if $C'\leq\frac{\tilde{C}'}{2(l(r-r''))^{D_1\gamma}}$, since $l(r-r'')$ has a bound which is independent of $r-r''$. $\Box$

\bigskip
\noindent This proposition is the inductive step which can be iterated as a whole. 
It is necessary to obtain an $\epsilon'$ which is much smaller than $\epsilon$ so as to control $|\Psi'|_{r'}$ as a function of $\epsilon'$ and make sure that the output be similar 
to the input.

\subsection{Main theorem}

\noindent First let us give a lemma which will enable us to iterate Proposition \ref{iter3}.

\begin{lem}\label{eps'1}

\noindent Let $C'\leq 1,b_0>0,r\leq \frac{1}{2}$ and $r'\in [\frac{95}{96}r,r[$. Let $D_5,\gamma_0\in\mathbb{N}$.There exists $C$ depending only on $C',D_5,\gamma_0$
 such that for all $\epsilon\leq C\left(\frac{r-r'}{b_0+1}\right)^{2\gamma_0D_5}$, 
choosing a sequence $(\epsilon_k)$ such that for all $k$,

\begin{equation}\nonumber \epsilon_k\leq  \epsilon_{k-1}^{100}<1\end{equation}

\noindent and letting for all $k$

\begin{equation}\nonumber \left\{\begin{array}{c}
\gamma_k=2^k\gamma_0\\
r_k=r'+\frac{r-r'}{2^k}\\
b_k=b_{k-1}+\mid \log \epsilon_{k-1}\mid\left(\frac{2^k}{r-r'}\right) ^{D_5}\\
\end{array}\right.\end{equation}

\noindent  
then for every $k\in\mathbb{N}$,



\begin{equation}\label{epsk0}\mid \log \epsilon_k\mid ^{2D_5\gamma_k} \leq \epsilon_k^{-\frac{1}{4}}
\end{equation} 

\noindent and

\begin{equation}\label{epsk1}
\left(\frac{b_k+1 }{r_k-r_{k+1}}\right)^{D_5\gamma_k}\epsilon_k\leq C'
\end{equation}


\end{lem}

\dem Let us first prove \eqref{epsk0}. It is equivalent to
%


\begin{equation}\nonumber 2^{k+3}D_5\gamma_0\leq \frac{ \mid \log \epsilon_k\mid}{\log\mid \log \epsilon_k\mid}
\end{equation}

\noindent The function $t\mapsto \frac{ \mid \log t\mid}{\log\mid \log t\mid}$ is decreasing for $t\in ]0, e^{-\frac{1}{e}}]$ so it is enough to show that

\begin{equation}\nonumber 2^{k+3}D_5\gamma_0\leq \frac{100^k \mid \log \epsilon\mid}{k\log 100+\log\mid \log \epsilon\mid}
\end{equation}

\noindent which is true if we choose $C$ as a function of $D_5,\gamma_0$.

\noindent 
$\bullet$ Let $a_k=\left(\frac{b_k+1 }{r_k-r_{k+1}}\right)^{D_5\gamma_k}\epsilon_k$. For all $k$, 

\begin{equation}\nonumber \begin{split}
a_{k+1}
&=\left(\frac{(b_{k+1}+1)2^{k+2}}{r-r'}\right)^{D_5\gamma_{k+1}}\epsilon_{k+1}\\
&\leq \left(\frac{(b_0+(k+1)\mid \log \epsilon_{k}\mid)2^{k+2}}{r-r'}\right)^{2D_5\gamma_{k+1}}\frac{\epsilon_{k+1}}{\epsilon_k}a_k\\
\end{split}\end{equation}

\noindent so by \eqref{epsk0}, 

\begin{equation}\nonumber \begin{split}
a_{k+1}
& \leq \left(\frac{(b_0+1)}{r-r'}\right)^{\gamma_016^{k+1}D_5}\epsilon^{100^k.98}a_k
\end{split}\end{equation}

\noindent thus, if $\epsilon$ is also smaller than $(\frac{r-r'}{b_0+1})^{16\gamma_0D_5}$, then $a_{k+1}\leq a_k$. If $\epsilon$ 
is also small enough to satisfy

\begin{equation}\nonumber a_0=\left(\frac{b_0+1 }{r-r'}\right)^{D_5\gamma_0}\epsilon\leq C'
\end{equation}

\noindent for instance

\begin{equation}\nonumber \left(\frac{b_0+1 }{r-r'}\right)^{D_5\gamma_0} \epsilon^{\frac{3}{4}}\leq C'
\end{equation}

\noindent then \eqref{epsk1} is true for all $k$. $\Box$




\bigskip
\noindent 
Lemma \ref{eps'1} implies that assumption \eqref{epsA+} of Proposition \ref{iter3} holds for all $k$ with $\epsilon\leq \epsilon_k$, $||A||=b_k$, $r=r_k$ and $r''=r_{k+1}$.

\bigskip
\noindent As a consequence, one gets the main result, of which we will give various formulations.

\begin{thm}\label{PR} Let $r\leq \frac{1}{2},A\in \mathcal{G}$ and $F\in C^\omega_r(2\mathbb{T}^d, \mathcal{G})$ with nice periodicity properties with respect to 
$\mathcal{L}_A$. 
Let

$$r'\in [\frac{95}{96}r,r[$$

\noindent 
There exists $D_7$ depending only on $n,d,\tau,\kappa,A$
such that if

\begin{equation}\label{petitesse}|F|_r\leq\epsilon_0'(r,r')= \left(\frac{r-r'}{\mid \mid A\mid \mid +1}\right)^{D_7}
\end{equation}

\noindent then for any $\epsilon\leq \epsilon_0'$, there exists 

\begin{itemize}
\item $Z_\epsilon,\Psi_\epsilon\in C^\omega_{r'}(2\mathbb{T}^d,G)$,
\item $A_\epsilon\in \mathcal{G}$, 
\item $\bar{A}_\epsilon,\bar{F}_\epsilon\in C^\omega_{r'}
(2\mathbb{T}^d,\mathcal{G})$, 
\end{itemize}
 
\noindent such that 
\begin{enumerate}
\item \label{Ared} $\bar{A}_\epsilon$ is reducible to ${A}_\epsilon$ by $\Psi_\epsilon$,

\item \label{3'} $|\bar{F}_\epsilon|_{r'}\leq \epsilon$

\item \label{2'} for every $\theta\in 2\mathbb{T}^d$,

$$\partial_\omega Z_\epsilon(\theta)=(A+F(\theta))Z_\epsilon(\theta)-Z_\epsilon(\theta)(\bar{A}_\epsilon(\theta)
+\bar{F}_\epsilon(\theta))
$$

\item \label{1'} 

$$|Z_\epsilon-Id|_{r'}\leq 2^{D_7}\epsilon_0^{\frac{1}{4}-4(r-r')}$$ 

\noindent and so does $Z_\epsilon^{-1}-Id$,


\item \label{4''} $Z_\epsilon, \partial_\omega Z_\epsilon$ are bounded in $C^\omega_{r'}(2\mathbb{T}^d,
gl(n,\mathbb{C}))$ uniformly in $\epsilon$,

\item \label{borneconj}

\begin{equation}\nonumber \mid \Psi_\epsilon\mid_{r'} \leq  \epsilon^{-(\frac{1}{2})^{ c' \sqrt{\log \mid \log \epsilon \mid}}}
\end{equation}

\noindent where $c'$ only depends on $n,d,\kappa,\tau,A$.

\end{enumerate}

\noindent Moreover, 
\begin{itemize}
\item in dimension 2 or if $\mathcal{G}=gl(n,\mathbb{C})$ or $u(n)$, if $F$ is continuous on $\mathbb{T}^d$, then 
$\bar{A}_\epsilon, \bar{F}_\epsilon$ and $Z_\epsilon$ are continuous on $\mathbb{T}^d$. 
\item If $\mathcal{G}$ is $o(n)$ or $u(n)$, then $D_7$ does not depend on $A$ and the same holds replacing \eqref{petitesse} by $\mid F\mid _r\leq (r-r')^{D_7}$.
\item if $\mathcal{G}=sl(2,\mathbb{C})$ or $sl(2,\mathbb{R})$ and $A+F$ is not reducible, then there exists a sequence $\epsilon_k\rightarrow 0$ such that $\mid \mid A_{\epsilon_k}\mid \mid \ \mid \log \epsilon_k\mid^\tau $ is bounded.
\end{itemize} 

\end{thm}

\dem 
The proof will be made by induction as follows. Let $r''=\frac{r+r'}{2}$.  
Let $R(r,r''),N(r,\epsilon)$, $\kappa''(r,r'',\epsilon)$ be as in \eqref{kappa''}. There exists $\gamma_0\in \mathbb{N}$ depending only on $n,d,\tau,\kappa,A$, such that $\mathcal{L}_A$ is an $(A,\kappa,\gamma_0)$-decomposition (one can assume $\gamma_0\geq n(n+1)$). 
Let $C'$, $D_3$ be as in Proposition \ref{iter3}. Let $D_5=2D_3$.
Let $C$ be as in Lemma \ref{eps'1} and $D_7$ such that

\begin{equation}\nonumber \left(\frac{r-r''}{||A||+1}\right)^{D_7}\leq C \left(\frac{r-r''}{||A||+1}\right)^{4\gamma_0D_5}
\end{equation} 

\noindent Let

\begin{equation}
 \epsilon_0'=\left(\frac{r-r''}{||A||+1}\right)^{D_7}
 \end{equation}

\noindent Before carrying on with the proof, note that if $\mathcal{G}$ is $o(n)$ or $u(n)$, then $\mathcal{L}_A$ is a unitary decomposition, therefore it is an $(A,\kappa,0)$-decomposition, 
so one can take $\gamma_0=n(n+1)$ and then $\gamma_0, D_3,D_5$ and $D_7$ do not depend on $A$. 
For all $k\in\mathbb{N}$, let

\begin{equation}\nonumber \left\{\begin{array}{c}
r_k=r''+\frac{r-r''}{2^k},\\ 
b_0=||A||,\\
b_{k+1}=||A||+\sum_{j\leq k} \frac{\mid \log \epsilon_j\mid}{(r_{j-1}-r_j)^{D_5}}\end{array}\right.\end{equation}

\noindent where $(\epsilon_j)$ will be defined by induction in the following. 
Suppose that $|F|_r\leq \epsilon_0'$. Let $\bar{F}_1=F, \bar{A}_1=A_1=A$ and $\Psi_0=Id$.
Iterate Proposition \ref{iter3} using lemma \ref{eps'1} to find, for all $k\geq 1$, 
\begin{itemize}
\item ${Z}_{k+1}\in C^\omega_{r_{k+1}}(2\mathbb{T}^d, G)$, 
\item ${A}_{k+1}\in
\mathcal{G}$,
\item $\bar{A}_{k+1}\in C^\omega_{r}
(2\mathbb{T}^d, \mathcal{G})$,
\item $\Psi_k\in C^\omega_{r}
(2\mathbb{T}^d,G)$,
\item $\bar{F}_{k+1}\in C^\omega_{r_{k+1}}(2\mathbb{T}^d, \mathcal{G})$
\item $\epsilon_{k+1}\in [\epsilon_k^{R(r_k,r_{k+1})^{n^2}},\epsilon_k^{100}]$
\end{itemize}
such that 
 
\begin{enumerate}

\item $\bar{A}_{k+1}$ is reducible to $A_{k+1}$ by $\Psi_k$,

\item $\Psi_k^{-1}\bar{F}_{k+1}\Psi_k$ has nice periodicity properties with respect to an
 $(A_{k+1},\kappa''(r_{k+1},r_{k+2},\epsilon_{k+1}),2^{k+1}\gamma_0)$-decomposition,

\item $|\bar{F}_{k+1}|_{r_{k+1}}
\leq \epsilon_{k+1}$,

\item $|\Psi_k|_{r}\leq \epsilon_{k+1}^{-\frac{1}{2}(r_{k+1}-r_{k+2})}$ and 
$|\Psi_k^{-1}|_{r}\leq \epsilon_{k+1}^{-\frac{1}{2}(r_{k+1}-r_{k+2})}$,

\item $||{A}_{k+1}||\leq b_{k+1}$, and, if $\mathcal{G}=sl(2,\mathbb{C})$ or $sl(2,\mathbb{R})$ and $\Psi_{k}^{-1}\Psi_{k-1}$ is not the identity, $||{A}_{k+1}||\leq 
\kappa''(r_k,r_{k+1},\epsilon_k)+ \epsilon_k^{\frac{1}{2}}$;

\item 

$$\partial_\omega {Z}_{k+1}=(\bar{A}_k+\bar{F}_k)
{Z}_{k+1}
-{Z}_{k+1}(\bar{A}_{k+1}+\bar{F}_{k+1})$$

\item 

$$|{Z}_{k+1}-Id|_{r_{k+1}}\leq 
\frac{1}{C'}\left(\frac{(1+||A_k||)|\log\epsilon_k|}
{r_k-r_{k+1}}\right)^{2^kD_3\gamma_0}
\epsilon_k^{1-4(r_k-r_{k+1})}$$ 

\noindent which implies, using Lemma \ref{eps'1}, that
%



$$|{Z}_{k+1}-Id|_{r_{k+1}}\leq 
\frac{1}{C'}\epsilon_k^{\frac{1}{4}-4(r_k-r_{k+1})}$$ 

\noindent and so does ${Z}_{k+1}^{-1}-Id$.
%


\end{enumerate}

\bigskip
\noindent $\bullet$ Let $\epsilon\leq \epsilon_0'$ and $k_\epsilon\in\mathbb{N}$ such that 
$\epsilon_{k_\epsilon+1}\leq \epsilon\leq \epsilon_{k_\epsilon}$. 
Let

\begin{equation}\nonumber \left\{\begin{array}{c}
Z_\epsilon={Z}_1\dots {Z}_{k_\epsilon}\\ 
\bar{A}_\epsilon=\bar{A}_{k_\epsilon}\\
\bar{F}_\epsilon=\bar{F}_{k_\epsilon}\\
\end{array}\right.\end{equation}
 
\noindent then properties \ref{Ared} and \ref{3'} hold. 
Thus for all $\theta\in 2\mathbb{T}^d$,

\begin{equation}\label{Zeps}\partial_\omega Z_\epsilon(\theta)=(A+F(\theta))Z_\epsilon(\theta)-
Z_\epsilon(\theta)(\bar{A}_{\epsilon}(\theta)
+\bar{F}_{\epsilon}(\theta))
\end{equation}

\bigskip
\noindent whence property \ref{2'}. Moreover, let $a_k:=|{Z}_1\dots {Z}_k-Id|_{r''}$, then 

\begin{itemize}
\item 

\begin{equation}\nonumber a_1=|Z_1-Id|_{r''}\leq \frac{1}{C'}\epsilon_0^{\frac{1}{4}-4(r_0-r_1)}
\end{equation}
%

%

\noindent so

\begin{equation}\nonumber |Z_1|_{r''}\leq 1+\frac{1}{C'}
\epsilon_0^{\frac{1}{4}-4(r_0-r_1)}
\end{equation}

\item let $k\geq 2$ and assume that for all $j\leq k-1$,

\begin{equation}\nonumber |Z_1\dots Z_j|_{r''}\leq 1+ \frac{3}{C'}
\epsilon_0^{\frac{1}{4}-4(r_0-r_1)}\end{equation}

\noindent then

\begin{equation}\nonumber \begin{split}a_k&\leq |{Z}_k-Id|_{r''}|{Z}_1\dots 
{Z}_{k-1}|_{r''}+a_{k-1}\\
&\leq   a_1+ \frac{1}{C'}\sum_{j=1}^{k-1} |{Z}_1\dots 
{Z}_{j}|_{r''}\epsilon_j^{\frac{1}{4}-4(r_j-r_{j+1})} \leq  \frac{3}{C'}
\epsilon_0^{\frac{1}{4}-4(r_0-r_1)}
\end{split}\end{equation}

%

%


%

%



\noindent and

\begin{equation}\nonumber \begin{split}|{Z}_1\dots {Z}_k|_{r''}
& \leq 1+\frac{3}{C'}
\epsilon_0^{\frac{1}{4}-4(r_0-r_1)}\end{split}\end{equation}

\end{itemize}
 
\noindent whence property \ref{1'}.
This also implies that

\begin{equation}\nonumber |Z_\epsilon|_{r''}\leq 2+\frac{3}{C'}\epsilon_0^{\frac{1}{4}-4(r_0-r_1)}\end{equation} 

\noindent Moreover, by a Cauchy estimate,

\begin{equation}\nonumber |\partial_\omega Z_\epsilon|_{r'}\leq \frac{1}{r''-r'}|Z_\epsilon|_{r''}\end{equation} 

 
\noindent so \ref{4''} is true. Also note that 

\begin{equation}\nonumber \mid \Psi_\epsilon\mid_{r''}=\mid \Psi_{k_\epsilon-1}\mid_{r''} \leq \epsilon_{k_\epsilon}^{-\frac{1}{2}(r_{k_\epsilon}-r_{k_\epsilon+1})}
 \leq \epsilon^{-\frac{1}{2^{k_\epsilon+2}}}
\end{equation}

\noindent Since

\begin{equation}\nonumber k_\epsilon(k_\epsilon+1)\geq c\log \left(\frac{\mid \log \epsilon \mid}{\mid \log (\mid F\mid_r)\mid} \right)
\end{equation}

\noindent where $c$ only depends on $n,d,\kappa,\tau,A$, then

\begin{equation}\nonumber \mid \Psi_\epsilon\mid_{r''} 
\leq  \epsilon^{-(\frac{1}{2})^{ c' \sqrt{\log \mid \log \epsilon \mid}}}
\end{equation}

\noindent  where $c'$ only depends on $n,d,\kappa,\tau,A$; therefore property \ref{borneconj} holds.


\bigskip
\noindent If $\mathcal{G}$ is either $gl(n,\mathbb{C})$ or $u(n)$ or in dimension 2, if $F$ is continuous on $\mathbb{T}^d$, 
each step will give functions ${Z}_{k+1},{A}_{k+1},\bar{A}_{k+1},\bar{F}_{k+1}$ continuous on $\mathbb{T}^d$ 
so, at the end of the process, the functions $Z_\epsilon,\bar{A}_\epsilon$ et $\bar{F}_\epsilon$ are continuous on 
$\mathbb{T}^d$. The fact that $\mid \mid A_{\epsilon_k}\mid \mid\ \mid \log \epsilon_k\mid^\tau$ is bounded for some sequence $\epsilon_k$ if $\mathcal{G}=sl(2,\mathbb{C})$ or $sl(2,\mathbb{R})$ and $A+F$ is not reducible comes from property 5 in the iteration. 
$\Box$

\bigskip
\noindent This proves Theorem \ref{th1}.

%

\bigskip
\noindent In general, almost reducibility does not imply reducibility. Reducibility happens if there are a finite number of steps at which one has to reduce the eigenvalues, or if the sequence $(\Psi_k)$ given by Theorem \ref{PR} converges in 
$C^\omega_{r'}(2\mathbb{T}^d,G)$. In general, this sequence is not even bounded in $C^\omega_{0}(2\mathbb{T}^d,G)$.
However, if the method above has been used to conjugate the system $A+F$ to a system $\bar{A}_\epsilon+\bar{F}_\epsilon$ where $\bar{A}_\epsilon$ is reducible by $\Psi_\epsilon$ to a constant 
$A_\epsilon$, and where $\bar{F}_\epsilon$ is bounded by $\epsilon$, 
one can also bound $\Psi_\epsilon^{-1}\bar{F}_\epsilon\Psi_\epsilon$.

\begin{cor}\label{infty} Let $r\leq \frac{1}{2},A\in \mathcal{G}$ and $F\in C^\omega_r(2\mathbb{T}^d, \mathcal{G})$ with nice periodicity properties with respect to $\mathcal{L}_A$. Let $r'\in [\frac{95}{96}r,r[$. There exists $D_8$ depending only on $n,d,\kappa,\tau,A$
such that if 

$$|F|_r\leq (r-r')^{D_8}$$
%

\noindent then there exists 
\begin{itemize}
\item $Z\in C^\omega
_{r'}(2\mathbb{T}^d,G)$, 
\item a family $(A_l)$ of reducible functions in 
$C^\omega_{r'}(2\mathbb{T}^d,\mathcal{G})$ 
\item and $A_\infty\in C^\omega_{r'}(2\mathbb{T}^d,\mathcal{G})$ 
\end{itemize}
such that

\begin{equation}\label{lim}\partial_\omega Z(\theta)=(A+F(\theta))Z(\theta)-Z(\theta)A_\infty(\theta)
\end{equation}

\noindent and

\begin{equation}\nonumber \lim_{l\rightarrow \infty} |A_l-A_\infty|_{r'}=0
\end{equation}

\noindent Moreover, in dimension 2 or if $\mathcal{G}=gl(n,\mathbb{C})$ or $u(n)$, if $F$ is continuous on $\mathbb{T}^d$, then $Z$, $A_l$ and $A_\infty$ are continuous on $\mathbb{T}^d$. Finally, if $\mathcal{G}=o(n)$ or $u(n)$, then $D_8$ does not depend on $A$.

\end{cor}

\dem Let $D_7$ be as in Theorem \ref{PR} and $D_8$ such that 

\begin{equation}(r-r')^{D_8}\leq \left(\frac{r-r'}{1+\mid \mid A\mid \mid }\right)^{D_7}
\end{equation}

Let $Z_\epsilon\in C^\omega_{r'}(2\mathbb{T}^d,G), 
A_\epsilon\in C^\omega_{r'}(2\mathbb{T}^d,\mathcal{G})$ be as in Theorem \ref{PR}. 
Then $Z_\epsilon$ and $\partial_\omega Z_\epsilon$ remain bounded in $C^\omega_{r'}
(2\mathbb{T}^d,G)$ when $\epsilon\rightarrow 0$. Let $Z$ be the limit in $C^\omega_{r'}
(2\mathbb{T}^d,G)$ of a subsequence 
$(Z_\frac{1}{k_l})$ 
of $(Z_\frac{1}{k})_{k\in\mathbb{N}\setminus\{0\}}$ 
and 

$$A_\infty(\theta):=Z(\theta)^{-1}(A+F(\theta))Z(\theta)-Z(\theta)^{-1}\partial_\omega Z(\theta)$$ 

\noindent then 

$$A_\infty
\in C^\omega_{r'}(2\mathbb{T}^d,\mathcal{G}), 
\lim_{l\rightarrow \infty} |A_\frac{1}{k_l}-A_\infty|_{r'}=0$$
 
\noindent and so equation 
\eqref{lim} holds. 

\bigskip
\noindent In dimension 2 or if $\mathcal{G}=gl(n,\mathbb{C})$ or $u(n)$, if $F$ is continuous on $\mathbb{T}^d$, all functions that one has to consider are continuous on $\mathbb{T}^d$. $\Box$

\bigskip
\rem In Corollary \ref{infty}, the function $A_\infty$ is not reducible in general, it is only a limit of reducible functions.

\begin{cor}\label{dense}Let $0<r'<r\leq \frac{1}{2}, A\in \mathcal{G}$ and $F\in C^\omega_r(\mathbb{T}^d,\mathcal{G})$. There exists 
$\epsilon_0'$ depending only on $n,d,\tau,\kappa,A,r-r'$ such that if $|F-A|_r\leq \epsilon_0'$, 
then for all $\epsilon>0$ there exists $H\in C^\omega_{r'}(2\mathbb{T}^d,\mathcal{G})$ such that 
$|F-H|_{r'}\leq \epsilon$ and $H$ is reducible.
\end{cor}

\dem Let $D_7$ be as in Theorem \ref{PR}. 
Assume that

\begin{equation}\nonumber |F-A|_r\leq (r-r')^{D_7}=:\epsilon_0'\end{equation}

\noindent Let $\epsilon>0$. 
By Theorem \ref{PR}, there exist $Z_\epsilon\in C^\omega_{r'}(2\mathbb{T}^d,G)$, 
$\bar{A}_\epsilon,\bar{F}_\epsilon\in C^\omega_{r'}(2\mathbb{T}^d,\mathcal{G})$ and $A_\epsilon\in 
\mathcal{G}$ 
such that

\begin{itemize}

\item $\bar{A}_\epsilon$ is reducible to $A_\epsilon$,

\item $\partial_\omega Z_\epsilon=FZ_\epsilon-Z_\epsilon(\bar{A}_\epsilon+\bar{F}_\epsilon)$,

\item $|Z_\epsilon|_{r'}\leq 2$, $|Z_\epsilon^{-1}|_{r'}\leq 2$,

\item $|\bar{F}_\epsilon|_{r'}\leq \frac{\epsilon}{4}$.

\end{itemize}

\noindent Therefore 

\begin{equation}\nonumber \partial_\omega Z_\epsilon=HZ_\epsilon-Z_\epsilon\bar{A}_\epsilon
\end{equation}

\noindent where $H=F-Z_\epsilon\bar{F}_\epsilon Z_\epsilon^{-1}$ is reducible to $A_\epsilon$ and satisfies

\begin{equation}\nonumber |H-F|_{r'}\leq 4|\bar{F}_\epsilon|_{r'}\leq \epsilon \ \Box
\end{equation}

\begin{cor} \label{dense'}Let $0<r'<r\leq \frac{1}{2}, A\in sl(2,\mathbb{R})$ and $F\in C^\omega_r(\mathbb{T}^d,sl(2,\mathbb{R}))$. There exists 
$\epsilon_0'$ depending only on $n,d,\tau,\kappa,A,r-r'$ such that if $|F-A|_r\leq \epsilon_0'$, 
then for any $\epsilon>0$ there exist $H\in C^\omega_{r'}(\mathbb{T}^d,sl(2,\mathbb{R}))$ such that 
$|F-H|_{r'}\leq \epsilon$ and $H$ is reducible.
\end{cor}

\dem Do the same construction as in Corollary \ref{dense}. Theorem \ref{PR} gives functions 
$\bar{A}_\epsilon,\bar{F}_\epsilon,Z_\epsilon$ which are, in fact, continuous on $\mathbb{T}^d$. Thus $H$ is continuous on $\mathbb{T}^d$. $\Box$

\bigskip
\noindent Corollary \ref{dense'} also holds with $gl(n,\mathbb{C})$ or $u(n)$ instead of $sl(2,\mathbb{R})$. 
This proves Theorem \ref{th3}. Again, note that if $G$ is a compact group, then the smallness condition does not depend on $A$.

\section{Strong almost reducibility for quasi-periodic cocycles in a Gevrey class}\label{PRGev}

Let $\beta>1$ and $r>0$. Let $C^{G,\beta}_r(2\mathbb{T}^d,gl(n,\mathbb{C}))$ be the functions of class Gevrey $\beta$ with parameter $r$, i.e the functions 
$F\in C^\infty(2\mathbb{T}^d,gl(n,\mathbb{C}))$ satisfying

$$\sum_{\alpha\in \mathbb{N}^d}\frac{r^{\beta \mid \alpha\mid}}{\alpha!^\beta}\sup_{\theta}\mid \mid \partial^\alpha F(\theta) \mid \mid <+\infty$$

\noindent 
Denote by $\mid \mid .\mid \mid_{\beta,r}$ the norm

$$\mid \mid F\mid \mid_{\beta,r}=\sum_{\alpha\in \mathbb{N}^d} \frac{r^{\beta \mid \alpha\mid}}{\alpha!^\beta} \sup_\theta \mid \mid \partial^\alpha F(\theta) \mid \mid$$

The main theorem in this part is formulated analogously to Theorem \ref{th1}.

\begin{thm}\label{th1G'}Let $0<r'<r\leq \frac{1}{2}$, $A\in \mathcal{G}$, $F\in C^{G,\beta}_r(\mathbb{T}^d, \mathcal{G})$. 
There is $\epsilon_0<1$ 
depending only on $n,d,\kappa,\tau,A,r-r'$ such that if

$$||F||_{\beta,r}\leq \epsilon_0$$ 

\noindent 
then for all $\epsilon>0$, there exists $\bar{A}_\epsilon,\bar{F}_\epsilon
\in C^{G,\beta}_{r'}
(2\mathbb{T}^d,\mathcal{G})$, 
$\Psi_\epsilon,Z_\epsilon\in C^{G,\beta}_{r'}(2\mathbb{T}^d, G)$ and $A_\epsilon\in \mathcal{G}$ 
such that for all $\theta\in 2\mathbb{T}^d$,

$$\partial_\omega Z_\epsilon(\theta)=(A+F(\theta))Z_\epsilon(\theta)-Z_\epsilon(\theta)(\bar{A}_\epsilon(\theta)
+\bar{F}_\epsilon(\theta))
$$

\noindent with
\begin{itemize}
\item $\bar{A}_\epsilon$ reducible to $A_\epsilon$ by $\Psi_\epsilon$, 
\item $||\bar{F}_\epsilon||_{\beta,r'}\leq \epsilon$, 

\item $\mid\mid  \Psi_\epsilon \mid \mid _{\beta,r'} \leq \epsilon^{-\frac{1}{8}}$,

\item and $||Z_\epsilon
-Id||_{\beta,r'}\leq 2\epsilon_0^\frac{1}{2}$. 
\end{itemize}

\bigskip
\noindent Moreover, 
\begin{itemize}
\item in dimension 2 or if $G=GL(n,\mathbb{C})$ or $U(n)$, $Z_\epsilon,\bar{A}_\epsilon,\bar{F}_\epsilon$ 
are continuous on $\mathbb{T}^d$; 
\item If $\mathcal{G}$ is $o(n)$ or $u(n)$, then $\epsilon_0$ does not depend on $A$;
\item if $\mathcal{G}=sl(2,\mathbb{C})$ or $sl(2,\mathbb{R})$ and $A+F$ is not reducible, then there exists a sequence $\epsilon_k\rightarrow 0$ such that $\mid \mid A_{\epsilon_k}\mid \mid \ \mid \log \epsilon_k\mid^\tau $ is bounded.
\end{itemize}
\end{thm}

\subsection{Preliminaries on Gevrey class functions}\label{prelGev}

\rem \begin{itemize}
\item For all $0<r'<r$, one has the inclusion $C^{G,\beta}_r( 2\mathbb{T}^d,\mathcal{G})\subset C^{G,\beta}_{r'}( 2\mathbb{T}^d,\mathcal{G})$ and 

$$\mid \mid f \mid \mid_{\beta,r'}
\leq \mid \mid f \mid \mid_{\beta,r}$$

\item For $f,g\in C^{G,\beta}_r(2\mathbb{T}^d,\mathcal{G})$, one has $||fg||_{\beta,r}\leq ||f||_{\beta,r}||g||_{\beta,r}$ (see 
\cite{MS}, appendix).
\end{itemize}

\begin{lem}\label{gevreyexp} For all $m\in\mathbb{Z}^d$ and all $r'>0$, the map $\theta\mapsto e^{2i\pi \langle m,\theta\rangle}$ satisfies 

\begin{equation}\nonumber \mid \mid e^{2i\pi \langle m,.\rangle}\mid \mid _{\beta,r'}\leq  e^{\beta\pi r' d |m|^{\frac{1}{\beta}}}
\end{equation}

 \end{lem}

\dem
For all $\alpha\in \mathbb{N}^d$ and all $\theta\in \mathbb{T}^d$,

\begin{equation}\nonumber \begin{split}\frac{r'^{\beta\mid \alpha\mid }}{(\alpha!)^\beta} \mid \partial^\alpha (e^{2i \pi \langle m,\theta \rangle})\mid 
&\leq \frac{r'^{\beta\mid \alpha\mid }}{(\alpha!)^\beta}  \prod_j \mid 2\pi m_j\mid ^{\alpha_j} \\
&\leq \prod_j \frac{( r'^\beta\mid 2\pi m_j\mid )^{ \alpha_j }}{(\alpha_j!)^\beta}\\
&\leq \prod_j \left( \frac{( r' \mid 2\pi m_j\mid^{\frac{1}{\beta}} )^{ \alpha_j }}{\alpha_j!}\right)^\beta
\end{split}\end{equation}

\noindent thus

\begin{equation}\nonumber \begin{split}\sum_{\alpha}\frac{r'^{\beta\mid \alpha\mid }}{(\alpha!)^\beta} \mid \partial^\alpha (e^{2i \pi \langle m,\theta \rangle})\mid 
&\leq \prod_j \sum_{\alpha_j}\left( \frac{( r' \mid 2\pi m_j\mid^{\frac{1}{\beta}} )^{ \alpha_j }}{\alpha_j!}\right)^\beta\\
&\leq  \prod_j \left(\sum_{\alpha_j} \frac{( r' \mid 2\pi m_j\mid^{\frac{1}{\beta}} )^{ \alpha_j }}{\alpha_j!}\right)^\beta\\
&\leq \prod_j e^{\beta r'\mid 2\pi m_j\mid^{\frac{1}{\beta}}} \leq  e^{\beta\pi r' d \mid m\mid^{\frac{1}{\beta}}} \ \Box
\end{split}\end{equation}

\noindent \rem This implies that the functions which are analytic on an $r$-neighbourhood of the torus or the double torus are Gevrey $\beta$ with parameter $r$ for all $\beta>1$;

\begin{soulem}\label{fouGevrey}Let $f\in C^{G,\beta}_r(2\mathbb{T}^d,gl(n,\mathbb{C}))$. Then for all $m\in \frac{1}{2}\mathbb{Z}^d$, 

\begin{equation}\nonumber 
\mid \mid \hat{  f} (m)\mid \mid \leq \mid \mid f\mid \mid_{\beta,r} (1-\frac{1}{2^{\frac{\beta}{\beta-1}}})^{1-\beta} e^{-\sum_j(2\pi  \mid m_j\mid)^{\frac{1}{\beta}}r}
\end{equation}

%


\end{soulem}

\dem By definition of $\mid\mid f\mid\mid_{\beta,r}$, 

\begin{equation}\nonumber \sum_{\alpha\in \mathbb{N}^d}\mid  \mid\widehat{ \partial^\alpha f} (m)\mid  \mid \frac{r^{\beta\mid \alpha \mid }}{\alpha!^\beta}
\leq\sum_{\alpha\in \mathbb{N}^d} \sup_\theta \mid \mid \partial^\alpha f(\theta)\mid \mid \frac{r^{\beta\mid \alpha \mid }}{\alpha!^\beta}
= \mid \mid f\mid \mid_{\beta,r}
\end{equation}

\noindent Now 

\begin{equation}\nonumber \partial^\alpha f(\theta)\sim \sum_m \hat{f}(m) \partial^\alpha(e^{2i\pi \langle m,\theta\rangle})
\sim \sum_m \hat{f}(m) \prod_j( 2i\pi m_j )^{\alpha_j}.(e^{2i\pi \langle m,\theta\rangle})
\end{equation}

\noindent thus

\begin{equation}\nonumber  \widehat{ \partial^\alpha f} (m)
=\prod_j ( 2i\pi m_j )^{\alpha_j}\hat{f}(m)\end{equation}

\noindent and therefore

\begin{equation}\label{alphagen}\mid \mid \hat{  f} (m)\mid \mid \sum_{\alpha\in \mathbb{N}^d}\prod_j( 2\pi \mid m_j\mid )^{\alpha_j}\frac{r^{\beta \alpha_j  }}{\alpha_j!^\beta}
\leq \mid \mid f\mid \mid_{\beta,r}
\end{equation}

\noindent that is to say,

\begin{equation}\label{alphagen'}\mid \mid \hat{  f} (m)\mid \mid \prod_j\sum_{\alpha_j\in \mathbb{N}}( 2\pi \mid m_j\mid )^{\alpha_j}\frac{r^{\beta\mid \alpha_j \mid }}{\alpha_j!^\beta}
\leq \mid \mid f\mid \mid_{\beta,r}
\end{equation}

\noindent Now, by Lemma \ref{sauzin},

\begin{equation}\nonumber \sum_{\alpha_j\in \mathbb{N}}( 2\pi  \mid m_j\mid r^\beta )^{\alpha_j}\frac{1}{\alpha_j!^\beta}
\geq (1-\frac{1}{2^{\frac{\beta}{\beta-1}}})^{\beta-1} e^{(2\pi \mid m_j\mid)^{\frac{1}{\beta}}r}
\end{equation}

\noindent therefore

\begin{equation}\nonumber \mid \mid \hat{  f} (m)\mid \mid \leq \mid \mid f\mid \mid_{\beta,r} (1-\frac{1}{2^{\frac{\beta}{\beta-1}}})^{1-\beta} e^{-\sum_j(2\pi  \mid m_j\mid)^{\frac{1}{\beta}}r}\ \Box
\end{equation}

\begin{lem}\label{troncG}Let $0<r\leq 1$, $f\in C^{G,\beta}_{r}(2\mathbb{T}^d,gl(n,\mathbb{C}))$, $N\in\mathbb{N}$ and $f^N$ the truncation of $f$ at order 
$N$. Then for all $r'<r$,

\begin{equation}\nonumber ||f-f^N||_{\beta,r'}\leq C_{d,\beta} \mid \mid f\mid\mid_{\beta,r} N^{d+1}\frac{1}{(r-r')^{2(d+1)}}e^{-2 (r-r') N^{\frac{1}{\beta}}}
\end{equation}

\noindent where $C_{d,\beta}$ only depends on $d,\beta$.
\end{lem}

\dem By definition,

\begin{equation}\nonumber  ||f-f^N||_{\beta,r'}= \sum_\alpha \sup_{\theta} \frac{r'^{\beta\mid \alpha\mid}} {\alpha!^\beta} \mid \partial^\alpha(f-f^N)(\theta)\mid
\end{equation}

\noindent Now

\begin{equation}\nonumber 
 \mid \partial^\alpha(f-f^N)(\theta)\mid \leq \sum_{\mid m\mid >N} \mid \mid \hat{f}(m)\mid\mid\ \prod_{j=1}^d \mid m_j\mid^{\alpha_j}
\end{equation}

\noindent so by Sublemma \ref{fouGevrey},

\begin{equation}\nonumber 
 \mid \partial^\alpha(f-f^N)(\theta)\mid \leq(1-\frac{1}{2^{\frac{\beta}{\beta-1}}})^{1-\beta}  \mid \mid f\mid \mid_{\beta,r}\sum_{\mid m\mid >N} e^{-\sum_l(2\pi  \mid m_l\mid)^{\frac{1}{\beta}}r} \prod_{j=1}^d \mid m_j\mid^{\alpha_j}
\end{equation}

\noindent whence

\begin{equation} \nonumber \begin{split}||f-f^N||_{\beta,r'}&\leq (1-\frac{1}{2^{\frac{\beta}{\beta-1}}})^{1-\beta}  \mid \mid f\mid \mid_{\beta,r} \sum_{\alpha} \frac{r'^{\beta\mid \alpha\mid}} {\alpha!^\beta} \sum_{\mid m\mid >N} e^{-\sum_l(2\pi  \mid m_l\mid)^{\frac{1}{\beta}}r} \prod_{j=1}^d \mid m_j\mid^{\alpha_j}\\
& \leq(1-\frac{1}{2^{\frac{\beta}{\beta-1}}})^{1-\beta}  \mid \mid f\mid \mid_{\beta,r} \sum_{\mid m\mid >N} e^{-\sum_l(2\pi  \mid m_l\mid)^{\frac{1}{\beta}}r} \prod_{j=1}^d \sum_{\alpha_j} \frac{r'^{\beta \alpha_j}} {\alpha_j!^\beta} \mid m_j\mid^{\alpha_j}
\end{split}\end{equation}

\noindent thus, using Lemma \ref{sauzin}, 

\begin{equation}\nonumber
\begin{split} ||f-f^N||_{\beta,r'}& \leq (1-\frac{1}{2^{\frac{\beta}{\beta-1}}})^{1-\beta}  \mid \mid f\mid \mid_{\beta,r} \sum_{\mid m\mid >N} e^{-2(r-r')\sum_j \mid m_j\mid ^{\frac{1}{\beta}}}\\
&\leq (1-\frac{1}{2^{\frac{\beta}{\beta-1}}})^{1-\beta}  \mid \mid f\mid \mid_{\beta,r} \sum_{\mid m\mid >N} e^{-2(r-r') \mid m\mid ^{\frac{1}{\beta}}}
\end{split}\end{equation}

\noindent and finally

\begin{equation}\nonumber ||f-f^N||_{\beta,r'} \leq C_{d,\beta}  \mid \mid f\mid \mid_{\beta,r} \sum_{M >N} M^de^{-2(r-r') M^{\frac{1}{\beta}}}
\leq C'_{d,\beta} \mid \mid f\mid \mid_{\beta,r} \frac{N^{d+1}}{(r-r')^{2(d+1)}}e^{-2(r-r')N^{\frac{1}{\beta}}} \end{equation}

\noindent where $C_{d,\beta},C'_{d,\beta}$ only depend on $d,\beta$. $\Box$

\subsection{Reduction of the eigenvalues}\label{redeigenG}

\noindent 
The reduction of the eigenvalues of a matrix $A$ at order $R,\bar{N}$ satisfies a good estimate in the Gevrey norm, as shows the following proposition:

\begin{lem}\label{renormG} Let $R,N\in \mathbb{N}\setminus\{0\},A\in gl(n,\mathbb{C})$ and $\Phi$ a map of reduction of the eigenvalues of $A$ at order $R,\bar{N}$. 
Then $\Phi$ satisfies for all $r'$ the Gevrey norm estimate

\begin{equation}\label{normephiG}||\Phi||_{\beta,r'}\leq nC.C_0 \left(\frac{1+||A_{\mathcal{N}}||}{\kappa''}\right)^{n(n+1)} 
e^{2\beta\pi r'd{\bar{N}^{\frac{1}{\beta}}}}
\end{equation}

\noindent where $C$ only depends on $d$, and so does $\Phi^{-1}$. Moreover, if $\mathcal{G}=o(n)$ or $u(n)$, then

\begin{equation}\label{normephiG'}||\Phi||_{\beta,r'}\leq nC
e^{2\beta\pi r'd{\bar{N}^{\frac{1}{\beta}} }}
\end{equation}

\noindent and so does $\Phi^{-1}$.
\end{lem}

\noindent \dem For all $m\in \mathbb{Z}^d$ and all $r'>0$, by Lemma \ref{gevreyexp},

\begin{equation}\nonumber \begin{split}||e^{2i\pi\langle m,.
\rangle}||_{\beta,r'}
\leq Ce^{\beta\pi r'd |m|^{\frac{1}{\beta}}}
\end{split}\end{equation}

\noindent where $C$ only depends on $d$. Therefore

\begin{equation}\nonumber \begin{split}
||\Phi||_{r'}&\leq \sum_{L\in \mathcal{L}_{A,\kappa''}} \mid \mid P^{\mathcal{L}_{A,\kappa''}}_L \mid \mid \ \mid\mid e^{2i\pi \langle m_L,.\rangle}\mid\mid_{r'}\\
&\leq C \sum_{L\in \mathcal{L}_{A,\kappa''}} \mid \mid P^{\mathcal{L}_{A,\kappa''}}_L \mid \mid e^{\beta\pi r'd |m_L|^{\frac{1}{\beta}}}\\
&\leq C \sum_{L\in \mathcal{L}_{A,\kappa''}} \mid \mid P^{\mathcal{L}_{A,\kappa''}}_L \mid \mid e^{2\beta\pi r'd \bar{N}^{\frac{1}{\beta}}}
\end{split}\end{equation}

\noindent Now by Lemma \ref{c0}, 

\begin{equation}\nonumber \mid \mid P^{\mathcal{L}_{A,\kappa''}}_L \mid \mid \leq 
C_0\left(\frac{1+\mid \mid A_{\mathcal{N}}\mid \mid }{\kappa''}\right)^{n(n+1)}
\end{equation}

\noindent so \eqref{normephiG} holds. 
If $\mathcal{G}$ is either $o(n) $ or $u(n)$, then $\mathcal{L}_{A,\kappa''}$ is a unitary decomposition and thus $P^{\mathcal{L}_{A,\kappa''}}_L$ has norm 1, 
thus \eqref{normephiG'} holds. 
$\Box$

\subsection{Homological equation}\label{homoleqG}


\begin{lem}\label{gevrey}Let $0<r'<r\leq 1$, $f\in C^{G,\beta}_r(2\mathbb{T}^d,\mathcal{G})$ and 
$g\in C^{G,\beta}_{r'}(2\mathbb{T}^d,\mathcal{G})$. Let ${C>0,D\geq 0}$. Assume that for all $m\in \frac{1}{2}\mathbb{Z}^d$, 

\begin{equation}\nonumber ||\hat{g}(m)||\leq C |m|^{D}||\hat{f}(m)||\end{equation}

\noindent 
Then

\begin{equation}\nonumber ||g||_{\beta,r'}\leq C'C ||f||_{\beta,r}\left(\frac{1}{r-r'}\right)^{2\beta(D+2)d}
\end{equation}

\noindent where $C'$ only depends on $d,D,\beta$.
\end{lem}

\dem For all $\theta\in N\mathbb{T}^d$ and all $\alpha\in \mathbb{N}^d$,

\begin{equation}\nonumber \begin{split}||\partial^\alpha g(\theta)||&\leq \sum_{m\in \frac{1}{2}\mathbb{Z}^d} \mid\mid \hat{g}(m)\mid \mid \ \mid  \partial^\alpha e^{2i\pi \langle m,\theta\rangle}    \mid\\
&\leq \sum_{m\in \frac{1}{2}\mathbb{Z}^d}  \mid\mid \hat{g}(m) \mid \mid \ \prod_j \mid 2\pi m_j\mid^{\alpha_j}\\
\end{split}\end{equation}

\noindent Therefore, by assumption,

\begin{equation}\nonumber \begin{split}||\partial^\alpha g(\theta)||
&\leq C\sum_{m\in \frac{1}{2}\mathbb{Z}^d}  |m|^{D}||\hat{f}(m)||    \prod_j \mid 2\pi m_j\mid^{\alpha_j}\\
&\leq C'C \sum_{m\neq 0} \frac{1}{\mid 2\pi m\mid ^{2d}}
\prod_j
\mid 2\pi m_j\mid ^{\alpha_j+D+2} ||\hat{f}(m)||\\
\end{split}\end{equation}

\noindent so, letting $\bar{1}=(1,\dots,1)$,

\begin{equation}\nonumber \begin{split}||\partial^\alpha g(\theta)||
&\leq C'C \sum_{m\neq 0} 
\frac{1}{\mid 2\pi m\mid ^{2d}}
 ||\widehat{\partial^{\alpha+(D+2)\bar{1}}f}(m)||\\
&\leq C' C\sup_\theta ||\partial^{\alpha+(D+2)\bar{1}}f(\theta)||
\end{split}\end{equation}



\noindent where $C'$ only depends on $d,D$, thus

\begin{equation}\nonumber \begin{split}||g||_{\beta,r'}&=\sum_{\alpha\in \mathbb{N}^d} (r')^{\beta\mid \alpha\mid}\frac{1}{(\alpha!)^\beta}\sup_\theta ||\partial^\alpha g(\theta)||\\
&\leq C'C\sum_{\alpha} (r')^{\beta\mid \alpha\mid }\frac{1}{(\alpha!)^\beta} \sup_{\theta} ||\partial^{\alpha+(D+2)\bar{1}}f(\theta)||\\
& \leq C' C\sum_{\alpha} \frac{r^{\beta\mid \alpha+(D+2)\bar{1}\mid }}{((\alpha+(D+2)\bar{1})!)^\beta} 
\sup_\theta  ||\partial^{\alpha+(D+2)\bar{1}} f(\theta)||
\frac{r'^{\beta\mid \alpha\mid}}{r^{\beta\mid \alpha+(D+2)\bar{1}\mid }}\left(\frac{(\alpha+(D+2)\bar{1})!}{\alpha!}\right)^\beta\\
& \leq C'C  \mid \mid f\mid \mid _{\beta,r} 
\sum_{\alpha}\frac{r'^{\beta\mid \alpha\mid}}{r^{\beta(\mid \alpha\mid +(D+2)d) }}\prod_j\left(\frac{(\alpha_j+D+2)!}{\alpha_j!}\right)^\beta\\
& \leq C'C  \mid \mid f\mid \mid _{\beta,r} 
\sum_{\alpha}\frac{r'^{\beta\mid \alpha\mid}}{r^{\beta(\mid \alpha\mid +(D+2)d) }}\left(\frac{(\mid \alpha\mid +D+2)!}{\mid \alpha\mid !}\right)^{\beta d}\\
& \leq C' C \mid \mid f\mid \mid _{\beta,r} 
\sum_{\alpha}\frac{r'^{\beta\mid \alpha\mid}}{r^{\beta(\mid \alpha\mid +(D+2)d) }}\left( \mid \alpha\mid +D+2\right)^{\beta(D+2)d}\\
\end{split}\end{equation}

\noindent Now the function

$$\phi:[0,+\infty[\rightarrow [0,+\infty[
,\ t\mapsto\left( \frac{r'}{r}\right)^t t^{\beta(D+2)d}$$ 

\noindent has its maximum at $t=\frac{\beta(D+2)d}{\ln \frac{r}{r'}}$ where it takes the value 
$e^{-\beta(D+2)d}\left(\frac{\beta(D+2)d}{\ln \frac{r}{r'}}\right)^{\beta(D+2)d}$. Therefore

\begin{equation}\nonumber \begin{split}||g||_{\beta,r'}&\leq C'C ||f||_{\beta,r}e^{-\beta(D+2)d}\left(\frac{\beta(D+2)d}{{r'}\ln \frac{r}{r'}}\right)^{\beta(D+2)d}\\
&\leq C'C ||f||_{\beta,r}e^{-\beta(D+2)d}\frac{(\beta(D+2)d)^{\beta(D+2)d}}{(r-r')^{2\beta(D+2)d}}\ \Box
\end{split}\end{equation}

\begin{prop}\label{homolG} 
Let
\begin{itemize}
\item ${N}\in\mathbb{N}$,
\item $\kappa'\in ]0, \kappa]$, 
\item $\gamma \geq n(n+1)$,
\item
 $0<r'<r$. 
 \end{itemize}
Let $\tilde{A}\in \mathcal{G}$ with a $DC^{{N}}_\omega(\kappa',\tau)$ spectrum.
Let $\tilde{F}\in C^{G,\beta}_r(2\mathbb{T}^d,  \mathcal{G})$ with nice periodicity properties with respect to an $(\tilde{A},\kappa',\gamma)$-decomposition $\mathcal{L}$. 
Then there exists a solution $\tilde{X}\in C^{G,\beta}_{r'}(2\mathbb{T}^d,  \mathcal{G})$ of equation 

\begin{equation}\label{homG}\forall \theta\in 2\mathbb{T}^d, \ \partial_\omega \tilde{X}(\theta)=[\tilde{A},\tilde{X}(\theta)]+\tilde{F}^{{N}}(\theta)
-\hat{\tilde{F}}(0);\ \hat{\tilde{X}}(0)=0
\end{equation}

\noindent such that 
\begin{itemize}
\item if $\tilde{F}$ has nice periodicity properties with respect to 
$\mathcal{L}$ et $(m_L)$, 
then so does $\tilde{X}$; in particular, if $\tilde{F}$ is defined on $\mathbb{T}^d$, then so is $\tilde{X}$,

\item Let $\Phi$ be trivial with respect to $\mathcal{L}$ and a family $(m_L)$ such that for all $L$, $\mid m_L\mid \leq N'$. 
There exists $C',D'$ depending only on 
$n,d,\tau,\beta$ such that

\begin{equation}\label{X2G}\mid \mid \Phi\tilde{X}\Phi^{-1}\mid \mid _{\beta,r'}
\leq C'\left(\frac{(1+||\tilde{A}_{\mathcal{N}}||)N'}{(r-r')\kappa'}\right)^{D'\gamma}
\mid \mid \Phi\tilde{F}\Phi^{-1}\mid \mid _{\beta,r}\end{equation}

\end{itemize}

\bigskip
Moreover, the truncation of 
$\tilde{X}$ at order $N$ is unique. 

\end{prop}

%


%

\dem The existence of $\tilde{X}$, its unicity up to order $N$, the fact that it takes its values in $\mathcal{G}$ and its nice periodicity properties with respect to $\mathcal{L}$ are proved as in part \ref{C2}, proposition \ref{homol}.

\bigskip
\noindent To get the estimate \eqref{X2G}, one first shows that for all $m\in \frac{1}{2}\mathbb{Z}^d$ and all $L,L'\in \mathcal{L}'$,

\begin{equation}\label{estimationreG}||P^{\mathcal{L}}_L\hat{\tilde{X}}(m)
P^{\mathcal{L}}_{L'}||
\leq C'\frac{(1+||\tilde{A}_{\mathcal{N}}||)^{n^2-1}|m|^{(n^2-1)\tau}}{\kappa'^{(n^2-1)}}
||P^{\mathcal{L}}_L\hat{\tilde{F}}(m)
P^{\mathcal{L}}_{L'}||  (||P^{\mathcal{L}}_L|| \ 
||P^{\mathcal{L}}_{L'}||)^{n^2-1}
\end{equation}

\noindent It is done exactly as in proposition \ref{homol} to get \eqref{estimationre}. The estimate \eqref{estimationreG} and Lemma \ref{gevrey} imply

\begin{equation}\nonumber ||P^{\mathcal{L}'}_L\tilde{X}e^{2i\pi \langle m_L-m_{L'}\rangle} 
P^{\mathcal{L}'}_{L'}||_{\beta,r'}
\leq C'' 
\left( \frac{(1+||\tilde{A}_{\mathcal{N}}||)N'}{(r-r')\kappa'}\right)^{D\gamma}
||P^{\mathcal{L}'}_L
\tilde{F}e^{2i\pi \langle m_L-m_{L'}\rangle}
P^{\mathcal{L}'}_{L'}||_{\beta,r}
\end{equation}

\noindent where $C'',D$ only depend on $n,d,\tau,\beta$. Thus

\begin{equation}\nonumber \begin{split}
||\Phi\tilde{X}\Phi^{-1}||_{\beta,r'}&\leq \sum_{L,L'}||P^{\mathcal{L}'}_L
\tilde{X}e^{2i\pi \langle m_L-m_{L'}\rangle}
P^{\mathcal{L}'}_{L'}||_{\beta,r}
\leq C''
\left(\frac{(1+||\tilde{A}_{\mathcal{N}}||)N'}{(r-r')\kappa'}\right)^{D\gamma}
\sum_{L,L'}||P^{\mathcal{L}'}_L
\tilde{F} e^{2i\pi \langle m_L-m_{L'}\rangle}P^{\mathcal{L}'}_{L'}||_{\beta,r}\\
\end{split}\end{equation}

\noindent and therefore

\begin{equation}\nonumber ||\Phi\tilde{X}\Phi^{-1}||_{\beta,r'}\leq C_3
\left(\frac{(1+||\tilde{A}||)N'}{(r-r')\kappa'}\right)^{D'\gamma}
||\Phi\tilde{F}\Phi^{-1}||_{\beta,r}\end{equation}

\noindent where $D',C_3$ only depend on $n,d,\tau,\beta$. $\Box$

\subsection{Inductive lemmas}\label{indlemG}

In Gevrey regularity, we will need a Lemma which is almost identical to Lemma \ref{NRdurable}, apart from the presence of the parameter $\beta$, which is fixed and does not modify the proof:

\begin{lem}\label{NRdurableG} Let
\begin{itemize}
\item $\kappa'\in ]0,1[$, $C>0$, $\beta>1$;
\item $\tilde{F}\in \mathcal{G}$,  
\item $\tilde{\epsilon}=||\tilde{F}||$,
\item $\tilde{N}\in \mathbb{N}$,
\item $\tilde{A}\in \mathcal{G}$ with $DC^{\tilde{N}}_\omega(\kappa',\tau)$ spectrum.

\end{itemize}

\noindent There exists a constant $c$ only depending on $n\tau,\beta $ 
such that if $\tilde{\epsilon}$ satisfies

\begin{equation}\label{cond1G}\tilde{\epsilon}\leq c\left(\frac{C^\tau \kappa'}{1+||\tilde{A}||}\right) ^{2n}
\end{equation} 

\noindent and

\begin{equation}\label{cond2G} \tilde{N}
\leq  \frac{|\log \tilde{\epsilon}|^\beta}{C}
\end{equation}

\noindent then
$\tilde{A}+{\tilde{F}}$ has $DC^{\tilde{N}}_\omega(\frac{3\kappa'}{4},\tau)$ spectrum.
\end{lem}

Exactly as in the analytic case, one obtains the following inductive lemmas, where only the estimates slightly differ from their analytic analogues; since $\beta$ is fixed, the proofs go on in quite the same way:

\begin{prop}\label{iterG} 
Let
\begin{itemize}
\item $\tilde{\epsilon}>0,\tilde{r}\leq 1$, $\tilde{r}'\in [\frac{\tilde{r}}{2},\tilde{r}[, \kappa'>0,\tilde{N}\in\mathbb{N},\gamma\geq n(n+1),C>0$; 

\item ${\tilde{F}}\in C^{G,\beta}_{\tilde{r}}(2\mathbb{T}^d,\mathcal{G}), \tilde{A}\in \mathcal{G}$,
\item $\mathcal{L}$ an $(\tilde{A},\kappa',\gamma)$-decomposition. 
\end{itemize}
There exists a constant $C''>0$ depending only on $\tau,n,\beta$ such that if

\begin{enumerate}

\item $\tilde{A}$ has $DC^{\tilde{N}}_\omega(\kappa',\tau)$ spectrum;

\item 

\begin{equation}\label{moyptG}||\hat{\tilde{F}}(0)||\leq  \tilde{\epsilon}\leq C''\left(\frac{C^\tau\kappa'}{1+||\tilde{A}||}\right)^{2n}\end{equation}

\noindent and

\begin{equation}\label{cond2'G} \tilde{N}
\leq  \frac{|\log \tilde{\epsilon}|^\beta}{C}
\end{equation} 

\item $\tilde{F}$ has nice periodicity properties with respect to $\mathcal{L}$

\end{enumerate}
 
\noindent then there exist
\begin{itemize}
\item $C'\in\mathbb{R}$ depending only on $n,d,\kappa,\tau,\beta$, 
\item $D\in\mathbb{N}$ depending only on $n,d,\tau,\beta$,
\item ${X}\in C^{G,\beta}_{\tilde{r}'}(2\mathbb{T}^d, \mathcal{G})$, 
\item ${A}'\in \mathcal{G}$ 
\item an $(A',\frac{3\kappa'}{4 },\gamma) $-decomposition $\mathcal{L}'$
\end{itemize}
satisfying the following properties: 

\begin{enumerate}

\item \label{nrG} $A'$ has $DC^{\tilde{N}}_\omega(\frac{3\kappa'}{4},\tau)$ spectrum,

\item \label{0-G} $||{A}'-\tilde{A}||\leq  \tilde{\epsilon}$;

\item \label{bpper-G} the map
$F'\in C^{G,\beta}_{\tilde{r}'}(2\mathbb{T}^d,\mathcal{G})$ defined by

\begin{equation}\label{4-G}
\forall \theta\in 2\mathbb{T}^d,\ \partial_\omega e^{{X}(\theta)}=(\tilde{A}+\tilde{F}(\theta))e^{{X}(\theta)}
-e^{{X}(\theta)}({A}'+{F}'(\theta))\end{equation}

\noindent has nice periodicity properties with respect to $\mathcal{L}'$

\item \label{9-G} If $\Phi$ is trivial with respect to $\mathcal{L}$ and a family $(m_L)$ satisfying, for all $L$, $\mid m_L\mid \leq N'$, then



\begin{equation}\label{pxpG}||\Phi^{-1}X\Phi||_{\beta,\tilde{r}'}
\leq C'\left(\frac{(1+||\tilde{A}_{\mathcal{N}}||)N'}{\kappa'(\tilde{r}-\tilde{r}')}\right)^{D\gamma}
||\Phi^{-1}\tilde{F}\Phi||_{\beta,\tilde{r}}\end{equation}

\item \label{10-G} and if $\Phi$ is trivial with respect to $\mathcal{L}$ and a family $(m_L)$ satisfying, for all $L$, $\mid m_L\mid \leq N'$, then 


\begin{equation}\nonumber \begin{split}||\Phi^{-1}F'\Phi||_{\beta,\tilde{r}'}&\leq C'\left(\frac{ (1+||\tilde{A}_{\mathcal{N}}||)N'}{\kappa'
(\tilde{r}-\tilde{r}')}\right)^{D\gamma}
e^{||\Phi^{-1}X\Phi||_{\beta,\tilde{r}'}}||\Phi^{-1}\tilde{F}\Phi||_{\beta,\tilde{r}}\\
&(||\Phi||^2_{\beta,\tilde{r}}||\Phi^{-1}||^2_{\beta,\tilde{r}}\tilde{N}^de^{-2\pi \tilde{N}(\tilde{r}-\tilde{r}')}+||\Phi^{-1}\tilde{F}\Phi||_{\beta,\tilde{r}'}
(1+e^{||\Phi^{-1}X\Phi||_{\beta,\tilde{r}'}}))\end{split}\end{equation}

\end{enumerate}

\noindent Moreover, if $\tilde{F}$ is continuous on $\mathbb{T}^d$, then so are $X$ and $F'$. If $\mathcal{G}=o(n)$ or $u(n)$, then the same holds replacing condition \eqref{moyptG} by

\begin{equation}\label{moypt'G}||\hat{\tilde{F}}(0)||\leq  \tilde{\epsilon}\leq C''(C^\tau\kappa')^2\end{equation}

\end{prop}

\noindent 
In Proposition \ref{iterG}, the only difference with Proposition \ref{iter} is the introduction of the parameters $\beta$ and $N'$. 
The inductive step is formulated exactly as proposition \ref{iter3}, with the only difference that the parameters $N,R$ will be chosen as

\begin{equation}\label{kappa''G}\left\{\begin{array}{c}
N(r,\epsilon)=(\frac{1}{2\pi r}|\log \epsilon|)^\beta\\
R(r,r')=[ \frac{1}{(r-r')^8}80^{4}(\frac{1}{2}n(n-1)+1)^2]^\beta\\
\end{array}\right.\end{equation}

\noindent and not as in \eqref{kappa''}. Note that the parameter $N'$ in the estimates of properties \ref{9-G} and \ref{10-G} in proposition \ref{iterG} does not modify essentially the proof, once it is instantiated by $\bar{N}$. The statement and the proof of the main theorem are identical.

\section{Appendix}

\subsection{Spectrum of a one-parameter family of matrices}

\begin{lem}\label{appendice}Let $\mathcal{G}$ be a Lie algebra and $A,F\in\mathcal{G}$ with $||F||\leq 1$. Let $\alpha_1(\lambda),
\dots ,\alpha_n(\lambda)$ be a continuous choice of the eigenvalues of $A+\lambda F$ as $\lambda$ varies from 
0 to 1. 
Then for all $1\leq j\leq n$, there exists $1\leq j'\leq n$
such that

\begin{equation}\nonumber 
\mid \alpha_{j'}(\lambda)-\alpha_j(0)\mid \leq 
2n\lambda^{\frac{1}{n}}(\mid \mid A\mid \mid +1)
\end{equation}

\end{lem}

\dem Fix $j\leq n$. For every $\lambda$, let 

\begin{equation}\nonumber A(\lambda)=A+\lambda F
\end{equation}

\noindent and

\begin{equation}\nonumber f(\lambda)=det (\alpha_j(0)I-A(\lambda ))
\end{equation}

\noindent Then $f(0)=0$ and for every $\lambda$,

\begin{equation}\nonumber f(\lambda)=det (\alpha_j(0)I-A(\lambda))=\prod_{j'} (\alpha_j(0)-\alpha_{j'}(\lambda))
\end{equation}

\noindent so

\begin{equation}\nonumber \mid \prod_{j'} (\alpha_j(0)-\alpha_{j'}(\lambda))\mid
= \mid f(0)-f(\lambda)\mid \leq  \sup_{\lambda''} \mid f'(\lambda'')\mid \ \mid \lambda\mid 
\end{equation}

\noindent and since

\begin{equation}\nonumber \begin{split} \mid f'(\lambda'')\mid & =\mid \sum_\sigma \frac{d}{d\lambda''}\prod_k (\alpha_j(0)I-A-\lambda'' F)_{k,\sigma(k)} \mid \\
&\leq n n! [\mid \mid A\mid \mid + \mid \mid A(\lambda'')\mid \mid  ]^{n-1}\\
&\leq 2^{n-1}n n!  [\mid \mid A\mid \mid + 1  ]^{n-1}
\end{split}\end{equation}

\noindent 
then there exists $j'$ such that

\begin{equation}\nonumber \mid \alpha_j(0)-\alpha_{j'}(\lambda)\mid
\leq 2n\mid \lambda\mid ^{\frac{1}{n}} [\mid \mid A\mid \mid + 1] \ \Box  \end{equation}
%

%

\noindent In case $G$ is compact, we have the following lemma (see \cite{El97}, lemma A.5): 

\begin{lem}\label{appendice2} Let $\mathcal{G}=o(n)$ or $u(n)$ and $A,F\in\mathcal{G}$ with $||F||\leq 1$. Let $\alpha_1(\lambda),
\dots ,\alpha_n(\lambda)$ be an analytic choice of the eigenvalues of $A+\lambda F$ as $\lambda$ varies from 
0 to 1. 
Then for all $1\leq j\leq n$,

\begin{equation}\mid \alpha_j(\lambda)-\alpha_j(0)\mid \leq \lambda
\end{equation}

\end{lem}

\dem For each $\lambda$, let $p_1(\lambda), \dots, p_n(\lambda)$ be an orthonormal basis of eigenvectors of $A+\lambda F$ (take them analytic in $\lambda$). Then for each $1\leq j \leq n$, one can assume

\begin{equation}\nonumber (A+\lambda F)p_j(\lambda)=\alpha _j(\lambda)p_j(\lambda)
\end{equation}

\noindent and derivating this along $\lambda$, one gets

\begin{equation}\nonumber (A+\lambda F-\alpha_j(\lambda))p_j'(\lambda)+(F-\alpha _j'(\lambda))p_j(\lambda)=0
\end{equation}

\noindent Now let $\beta_1,\dots, \beta_n$ be such that 

\begin{equation}\nonumber p_j'(\lambda)=\sum_{l=1}^n \beta_l p_l(\lambda)
\end{equation}

\noindent Then

\begin{equation}\nonumber \sum_{l\neq j} \beta_l(A+\lambda F-\alpha_j(\lambda))p_l(\lambda)+(F-\alpha _j'(\lambda))p_j(\lambda)=0
\end{equation}

\noindent and taking the scalar product with $p_j(\lambda)$,

\begin{equation}\nonumber \langle (Fp_j(\lambda),p_j(\lambda)\rangle= \alpha _j'(\lambda) 
\end{equation}

\noindent Therefore

\begin{equation}\nonumber \mid \alpha_j(0)-\alpha_j(\lambda)\mid \leq \mid \lambda\mid \sup_{\lambda '} \mid \alpha _j'(\lambda') \mid
\leq \lambda \ \Box
\end{equation}

\subsection{A lemma on integer series with non-negative terms}

\noindent The following lemma was proven by D. Sauzin and is used in Section \ref{PRGev}.

\begin{lem}\label{sauzin}
Let $a>0$. For all $r\geq 0$, consider $E_a(r)=\sum_{k\geq 0} \frac{r^k}{k! ^a}$.
Then 

\begin{equation}
\left\{ \begin{array}{ccc}
K_1 e^{\lambda a r^{\frac{1}{\alpha}}} \leq E_a(r) \leq e^{ar^{\frac{1}{a}}} & \mathrm{if} & a>1, 0<\lambda < 1\\
e^{ar^{\frac{1}{a}}}\leq E_a(r) \leq K_1 e^{\lambda ar^{\frac{1}{a}}} & \mathrm{if} & 0<a<1, \lambda>1\\
\end{array}\right.
\end{equation}

\noindent where 

$$K_1=(1-\lambda^{\frac{a}{a-1}})^{a-1}<1$$

\end{lem}

\dem One uses the following inequalities: if $\alpha>1$ and $(X_k)_{k\in\mathbb{N}},(Y_k)_{k\in\mathbb{N}}$ are families of non-negative numbers,

\begin{equation}\label{1}\sum X_k^\alpha \leq (\sum X_k)^\alpha
\end{equation}

\noindent and

\begin{equation}\label{2}\sum X_k^\alpha \geq \frac{ (\sum X_k Y_k)^\alpha} {(\sum Y_k^\beta)^{\frac{\alpha}{\beta}} }
\end{equation}

\noindent for $\beta=\frac{\alpha}{\alpha-1}$, if $\sum X_kY_k <\infty$. Note that \eqref{1} is equivalent to

\begin{equation}\label{1''} \sum x_k^{\frac{1}{\alpha}} \geq (\sum x_k)^{\frac{1}{\alpha}}
\end{equation}

\noindent with $X_k=x_k^{\frac{1}{\alpha}}$. Also \eqref{2} is equivalent to

\begin{equation}\label{2''}\sum x_k ^{\frac{1}{\alpha}} \leq (\sum y_k ^{\frac{\alpha}{\alpha-1} })^{1-\frac{1}{\alpha}} (\sum \frac{x_k}{y_k})^{\frac{1}{\alpha}}
\end{equation}

\noindent with $X_k=(\frac{x_k}{y_k})^{\frac{1}{\alpha}}$ and $y_k=Y_k^\alpha$.

\bigskip
\begin{enumerate}
\item In the first case, apply \eqref{1} and \eqref{2} with $\alpha=a, X_k= \frac{r^{\frac{k}{\alpha}}} {k!}$ and $Y_k=\lambda^k$.

\item In the second case, apply \eqref{1''} and \eqref{2''} with $\alpha=\frac{1}{a},x_k=\frac{r^{\frac{k}{a} } }  {k!}$ and $y_k=\lambda^{-k}$. $\Box$
\end{enumerate}

\section*{Acknowledgments}

The author would like to thank H\aa kan Eliasson for his supervision, as well as David Sauzin, Raphaël Krikorian, Artur Avila and Xuanji Hou for suggestions and useful discussions.

\end{document}